\documentclass{NP-IV}
\usepackage[latin1]{inputenc}
\usepackage{amsmath}
\usepackage{amscd}
\usepackage{mathrsfs}
\usepackage[all]{xy}
\usepackage{graphicx}
\usepackage{pstricks} 
\usepackage{amsfonts}   
\usepackage{pict2e}   
\usepackage{makeidx}
\usepackage[colorlinks=true,linkcolor=red,citecolor=red]{hyperref}

\usepackage{amssymb}
\usepackage{latexsym}
\usepackage{bold-extra}

\newcommand{\Bs}{\ensuremath{\mathscr{B}}}
\newcommand{\bs}[1]{\boldsymbol{#1}}
\newcommand{\C}{\mathscr{C}}

\newcommand{\chidr}{\chi^{\phantom{1}}_{\mathrm{dR}}}
\newcommand{\Ds}{\mathscr{D}}

\newcommand{\RR}{\bs{\mathcal{R}}}

\renewcommand{\H}{\mathscr{H}}

\newcommand{\M}{\mathrm{M}}

\newcommand{\N}{\mathrm{N}}

\renewcommand{\O}{\mathscr{O}}

\renewcommand{\sp}{\mathrm{sp}}

\newcommand{\R}{\mathcal{R}}

\newcommand{\is}[1]{i_{#1}^{\mathrm{sol}}}

\newcommand{\simto}{\stackrel{\sim}{\to}}

\DeclareMathOperator\Int{Int}
\DeclareMathOperator\Card{Card}
\DeclareMathOperator\WC{Card_{w}}

\DeclareMathOperator\Irr{Irr}

\DeclareMathOperator\charac{char}

\newcommand\topo{\mathrm{top}}
\newcommand\an{\mathrm{an}}
\newcommand\alg{\mathrm{alg}}
\newcommand{\ho}{\ensuremath{\hat{\otimes}}}

\newcommand\wti{\widetilde}
\newcommand{\tot}{\mathrm{tot}}

\newcommand\emb{\mathrm{emb}}
\newcommand\lin{\mathrm{lin}}
\newcommand\vw{v_{\mathrm{w}}}

\newcommand\ew{e_{\mathrm{w}}}
\newcommand\rw{\mathrm{w}}

\newcommand{\Cs}{\mathscr{C}}

\newcommand\Es{\mathscr{E}}
\newcommand\Fs{\mathscr{F}}
\newcommand\Gs{\mathscr{G}}
\newcommand\Hs{\mathscr{H}}

\newcommand\Ps{\mathscr{P}}

\newcommand\Ic{\mathcal{I}}
\newcommand\Rc{\mathcal{R}}

\newcommand\NN{\mathbb{N}}
\newcommand\ERRE{\mathbb{R}}
\newcommand\Z{\mathbb{Z}}

\newcommand{\E}[2]{\ensuremath{\mathbb{A}^{#1,\mathrm{an}}_{#2}}}
\newcommand\of[3]{\mathopen#1 #2 \mathopen#3}

\newcommand\wc{{\mkern 2mu\cdot\mkern 2mu}}
\newcommand\va{|\wc|}

\newcommand\wKa{\widehat{K^\alg}}
\newcommand\wLa{\widehat{L^\alg}}

\newcommand{\comm}[1]{
\noindent{\magenta \framebox{\framebox{\framebox{
\begin{minipage}{450pt}#1
\end{minipage}}}}}
}

\usepackage{amsthm}
\swapnumbers

\makeatletter
\def\swappedhead#1#2#3{%
  \thmname{#1}\;%
  \thmnumber{\@upn{\the\thm@headfont#2\@ifnotempty{#1}}}%
  \thmnote{\,{\the\thm@notefont(#3)}}{.~}}
\makeatother

\newtheoremstyle{dotless-thm}
  {10pt}
  {10pt}
  {\itshape}
  {}
  {\bfseries}
  {}
  {.0em}
  {}
\theoremstyle{dotless-thm}

\newtheorem{theorem}{\textbf{Theorem}}[subsection]
\newtheorem{thm-intro}{\textbf{\textsc{Theorem}}}
\newtheorem{rk-intro}[thm-intro]{\textbf{\textsc{Remark}}}
\newtheorem{cor-intro}[thm-intro]{\textbf{\textsc{Corollary}}}
\newtheorem{proposition}[theorem]{\textbf{Proposition}}
\newtheorem{lemma}[theorem]{\textbf{Lemma}}
\newtheorem{corollary}[theorem]{\textbf{Corollary}}
\newtheorem{definition}[theorem]{\textbf{Definition}}
\newtheorem{remark}[theorem]{\textbf{Remark}}

\newtheorem{notation}[theorem]{\textbf{Notation}}
\numberwithin{equation}{section}

\title[Convergence Newton 
polygons IV : graphs]{The convergence Newton polygon 
of a $p$-adic differential equation IV : 
controlling graphs}

\author{J\'er\^ome Poineau}
\email{jerome.poineau@math.lnmo.fr}
\address{Laboratoire de mathématiques Nicolas Oresme
Université de Caen, BP 5186, F-14032, Caen Cedex}

\author{Andrea Pulita}
\email{andrea.pulita@univ-grenoble-alpes.fr}
\address{Univ. Grenoble Alpes, CNRS, IF, 38000 Grenoble, France.}

\date{\today}



\begin{abstract}
In our previous works we proved a finiteness property of 
the radii of convergence functions associated with 
a vector bundle with connection on $p$-adic analytic 
curves. We showed that the radii are locally constant 
functions outside a locally finite graph in the curve, called 
controlling graph.
In this paper we refine that finiteness results by giving a 
bound on the size of the controlling graph in terms of 
the geometry of the curve and the rank of the module. 
This is based on super-harmonicity properties of 
radii of convergence and partial heights of the Newton 
polygon. 
Under suitable assumptions, we relate the size of the 
controlling graph associated with 
the total height of the convergence 
Newton polygon to the Euler 
characteristic in the sense of de Rham cohomology.
\end{abstract}

\begin{document}
\maketitle

\begin{center}
Version of \today
\end{center}

\makeatletter
\renewcommand\tableofcontents{%
    \subsection*{\contentsname}%
    \@starttoc{toc}%
    }
\makeatother

\begin{small}
\setcounter{tocdepth}{3} \tableofcontents
\end{small}

\setcounter{section}{0}


\section*{\textsc{Introduction}}
\addcontentsline{toc}{section}{\textsc{Introduction}}


Since the inception of the $p$-adic theory of differential equations, initiated by 
the pioneering works of Bernard Dwork \cite{Dwork-zeta, Dw-hyp} and 
Philippe Robba \cite{Ro-I, Ro-II, RoIII, RoIV}, the radii of convergence of the 
Taylor series solutions to these equations has proven to be a powerful tool for their classification.
Although the radii are invariants of topological nature 
of the differential equation, they provide significant algebraic 
information about it. For example, they enable the factorization 
of the equation (cf. \cite{Ro-I, Dw-Robba, Ch-Me-III, 
Kedlaya-Book, Kedlaya-draft, NP-III}) 
and contain the information regarding the \emph{finite dimensionality 
of its de Rham cohomology} (cf. \cite{RoIV, Pons, NP-V, NP-VI}).
A crucial property of the radii is that they vary in a continuous way 
over the curve where the equation is defined. Due to the totally 
disconnected nature of the non-archimedean topology, the language 
of Berkovich geometry is essential to give meaning to this continuity 
property (cf. \cite{Ch-Dw, DV-Balda, Balda-Inventiones}). 

Recently, we have significantly strengthened the continuity of the 
radii functions by proving that they are locally constant outside a 
locally finite graph (cf. \cite{NP-I, NP-II, Kedlaya-draft}). 
In other words, the radii are determined by their restriction to the graph. 
This implies that they are locally controlled by a finite set of 
numerical invariants. Each graph associated with a radius function is a 
topological invariant of the equation, known as the \emph{controlling graph}. 
It contains important algebraic information about the equation; 
for instance, the relative position of these graphs indicates whether 
the decomposition by the radii is a direct sum \cite[Theorem 5.2.12]{NP-III}.
We will demonstrate in \cite{NP-V, NP-VI} that this enhanced continuity 
property has significant implications for the finite dimensionality of 
their de Rham cohomology. Additionally, K.S. Kedlaya proved that the 
endpoints of the controlling graphs are never points of type 4 
(cf. \cite[Theorem 4.5.15]{Kedlaya-draft}).

In the present paper, we reinforce the results of \cite{NP-I, NP-II} 
by providing information about the size of the controlling graphs. 
Specifically, we establish a 
remarkable link between the size of the controlling graph associated with 
the total height of the convergence Newton polygon and the 
\emph{index formula in de Rham cohomology}. 
Moreover, for projective curves, we obtain a bound on the size of the 
graph that depends only on the rank of the equation. In other 
words, if the rank is fixed, the bound does not depend on the specific 
equation.

\if{
This text is the fourth paper of our series devoted to 
differential equations on $p$-adic analytic curves (see 
\cite{NP-I,NP-II,NP-III} for the first ones).
}\fi
Let us now be more precise about the results of this paper. Let $K$ be a $p$-adic field or, more generally, a complete valued field of characteristic~0. Let $X$ be a quasi-smooth Berkovich $K$-analytic curve and let $(\Fs,\nabla)$ be a locally free $\O_X$-module with connection of finite rank~$r$ on~$X$. For simplicity, we will often remove~$\nabla$ from the notation. Recall that the curve~$X$ has the structure of a real graph. This structure may be encoded in a so-called triangulation~$S$, or in the associated skeleton~$\Gamma_{S}$. The latter is a locally finite graph on which the curve~$X$ retracts by deformation.

Following the work of F.~Baldassarri \cite{Balda-Inventiones}, one may associate to any point $x\in X$ a family of radii of convergence $\R_{S,1}(x,\Fs),\dotsc, \R_{S,r}(x,\Fs)$. 
The papers \cite{NP-I,NP-II} are devoted to a precise study of those radii of convergence. In particular, for each $i\in \{1,\dotsc,r\}$, we proved that the map $\R_{S,i}(-,\Fs)$ is continuous on~$X$ and factors through a locally finite graph~$\Gamma_{S,i}(\Fs)$ (containing~$\Gamma_{S}$), called controlling graph of~$\R_{S,i}(-,\Fs)$. 

In this paper, we investigate more precisely the controlling graphs $\Gamma_{S,i}(\Fs)$  and compute explicit upper bounds on their number of vertices and edges. 
Our best-looking result is obtained for the first radius of convergence, when the curve is compact (see Corollary~\ref{cor:RS1edge}). For any finite subgraph~$\Gamma$ of~$X$, denote by~$e(\Gamma)$ (resp. $v(\Gamma)$) its number of edges (resp. vertices).

\begin{thm-intro}\label{Thm-intro: key point}
Assume that $X$ is the analytification of a smooth geometrically connected projective curve of genus~$g\ge 1$ and that $S$ is minimal. Then, we have
\begin{equation}
\begin{cases}
e(\Gamma_{S,1}(\Fs)) \le e(\Gamma_{S}) + 4(g-1) 
\cdot r \cdot\max(r-1,2);\bigskip\\
v(\Gamma_{S,1}(\Fs)) \le v(\Gamma_{S}) + 
4(g-1)\cdot r\cdot\max(r-1,2).
\end{cases}
\end{equation}
\end{thm-intro}
Notice that 
the right hand bounds \emph{are the 
same for all differential equations of a given rank}.

When the curve~$X$ is not compact, a similar result holds. In order for the controlling graph to be finite, one needs to require that the first radius~$\R_{S,1}(-,\Fs)$ is $\log$-affine at the open boundary of~$X$, and the bounds then involve the corresponding slopes (see Corollary~\ref{cor:RS1edge}). 

A similar statement holds for the higher radii~$\R_{S,i}(-,\Fs)$, still under an assumption of $\log$-affinity at the open boundary, although the bounds are less explicit (see Corollary~\ref{cor:HSjedgetotalCS}).

As an illustration, we apply those results in the particular 
case of elliptic curves. In this situation, they imply that all 
radii of convergence are constant and that differential 
modules split as direct sums of modules with all radii 
equal (see Corollary~\ref{Cor Ell curves}).

\medbreak

It is worth mentioning a particular result that holds for the controlling graph $\Gamma_{S}(H_{S,r}(-,\Fs))$ of the total height~$H_{S,r}(-,\Fs)$ of the  convergence Newton polygon, under additional assumptions, 
see Corollary~\ref{cor:HrIndex}. Its proof relies on cohomological computations and an index formula, which will appear in the forecoming papers~\cite{NP-V, NP-VI}. The statement exhibits a surprising relationship between the size of the controlling graph $\Gamma_{S}(H_{S,r}(-,\Fs))$ and the Euler characteristic~$\chidr(X,\Fs)$ in the sense of de Rham cohomology with coefficients in~$\Fs$.

\begin{thm-intro}\label{thintro:Hr}
Under suitable assumptions on~$S$ and~$(\Fs,\nabla)$, we have
\begin{equation}
\begin{cases}
e(\Gamma_{S}(H_{S,r}(-,\Fs))) \le e(\Gamma_{S}) - 2\cdot r\cdot \max(r-1,2)\cdot\chidr(X,\Fs);\bigskip\\ 
v(\Gamma_{S}(H_{S,r}(-,\Fs))) \le v(\Gamma_{S}) - 2\cdot r\cdot \max(r-1,2)\cdot\chidr(X,\Fs). 
\end{cases}
\end{equation}
\end{thm-intro}

\medbreak

The key ingredient in the proof of Theorem~\ref{Thm-intro: key point} and its generalizations is a control of the Laplacian of the first radius $dd^c \R_{S,1}(x,\Fs)$ or, more generally, a partial height $dd^c H_{S,i}(x,\Fs)$, where $H_{S,i}(-,\Fs)$ denotes the product of the first $i$ radii. For interior points belonging to the skeleton, we obtain a precise upper-bound (see 
Theorem~\ref{thm:curveTR}). The bound holds under a 
technical assumption (TR) which is automatically satisfied 
if the residual characteristic of $K$ is different from $2$ 
(cf. Definition \ref{def:TR}). 

\begin{thm-intro}[\protect{Theorem \ref{thm:curveTR}}]\label{thintro:TR}
Let~$x$ be a point in $\Gamma_{S} \cap \Int(X)$. If it is of type~2, assume that it satisfies the condition~$(TR)$. Then, for each $i\in\{1,\dotsc,r\}$, we have
\begin{equation}
dd^c H_{S,i}(x,\Fs)\; \le\; 
(2g(x) - 2 \deg(x) + N_{S}(x)) \cdot i,
\end{equation}
where $N_{S}(x)$ denotes the arity of~$x$ in~$\Gamma_{S}$ (counted with suitable multiplicities).
\end{thm-intro}

A similar formula appears in a paper by K.~Kedlaya, see \cite[Theorem~5.3.6]{Kedlaya-draft}. Note that the technical assumption~(TR) has later been removed in~\cite{Push} (which appeared after the first version of this paper was released).

It is also important to be able to control the locus where the super-harmonicity of the partial heights fails. To this end, we introduce some rather explicit exceptional sets~$\Cs_{S,i}(\Fs)$, relying crucially on prior work from \cite{NP-I} in the case of the affine line. When bounding the size of the controlling graphs, an important part of the work consists in bounding the size of those sets. We also obtain general upper-bounds for the Laplacians of the partial heights (see Theorem~\ref{thm:Ei} and Proposition~\ref{prop:HiCi}).

\begin{thm-intro}\label{thintro:ddcHi}
Let $x \in X$ and let $i\in \{1,\dotsc,r\}$. We have 
\begin{equation}
dd^c H_{S,i}(x,\Fs)\;\le\; i-1.
\end{equation}
If $x\notin S \cup \C_{S,i}(\Fs)$, then we have
\begin{equation}
dd^c H_{S,i}(x,\Fs)\;\le\; 0.
\end{equation}
\end{thm-intro}

\subsection*{Structure of the paper.}

In Section~\ref{Notations - 1}, we collect some 
definitions and results about graphs and $K$-analytic 
curves. In particular, we introduce a notion of 
pseudo-triangulation, slightly generalizing that of triangulation. 
The notation defined there will be constantly used in this 
paper and the next ones in the series. The section is 
quite long but necessary to develop our theory on a 
sound basis. In particular, pseudo-triangulations 
constitute a crucial step to define finite curves, which 
permit to \emph{characterize curves without boundary 
with finite dimensional de Rham cohomology (cf. \cite{NP-VI})}.

In Section~\ref{sec:sh}, we investigate the Laplacians of the first radius of convergence and the partial heights and prove Theorems~\ref{thintro:TR} and~\ref{thintro:ddcHi}. 

In Section~\ref{section : explicit bound}, we prove the explicit bounds on the size of the controlling graphs: Theorem~\ref{Thm-intro: key point} and its generalizations and Theorem~\ref{thintro:Hr}.


\subsection*{Acknowledgments.}
Thanks to Yves Andr\'e, Francesco Baldassarri, Gilles Christol, 
Kiran S. Kedlaya, Adriano Marmora, Zoghman Mebkhout, and 
Nobuo Tsuzuki for helpful comments.

\subsection*{Notation}
For $a,b\in\mathbb{R}$, we use the notation 
$[a,b]:=\{r\in\mathbb{R}, a\leq r\leq b\}$, (resp. 
$]a,b[=\{r\in\mathbb{R}, a< r <b\}$, 
$[a,b[:=\{r\in\mathbb{R}, a\leq r <b\}$) 
to indicate a closed (resp. open, semi-open) and bounded real interval. 
We also use the terminology \emph{segment} to indicate a bounded 
real interval.

Let $(K,\va)$ be a ultrametric complete valued field of characteristic $0$. Denote by $\wti K :=\{x\in K,|x|\leq 1\}/\{x\in K,|x|< 1\}$ its residue field and by~$p$ be the characteristic exponent of the latter (either~1 or a prime number). Let~$K^\alg$ be an algebraic closure of~$K$. The absolute value~$\va$ on~$K$ extends uniquely to it. We denote by $(\wKa,\va)$ its completion.

\section{Curves and graphs}
\label{Notations - 1}

In this section, we introduce some definitions and notation about $K$-analytic curves and graphs. They will be constantly used in this paper, and in the next papers of the series. 

We will use of the theory analytic geometry over~$K$ as developed by V.~Berkovich in~\cite{Ber,bleu}. 

Let $Z$ be a $K$-analytic space. For every complete valued extension~$L$ of~$K$, and every complete valued extension~$M$ of~$L$, we set $Z_{L} := Z \ho_{K} L$,  $Z_{M} := Z \ho_{K} M$, and we denote by $\pi_{M/L} \colon Z_{M} \to Z_{L}$ the canonical projection morphism. We set $\pi_{M} := \pi_{M/K}$.


\subsection{General graphs.}\label{Section : Graphs.}\label{sec:graph}
We borrow the definition of graph from~\cite[1.3.1]{Duc}. It is defined as a 
topological space.


\begin{definition}\label{def:graph}
We say that a topological space~$\Gamma$ is a \emph{graph} if it satisfies the following properties:
\begin{enumerate}
\item $\Gamma$ is Hausdorff and locally compact\;;
\item $\Gamma$ has a basis of open sets~$U$ such that
\begin{enumerate}
\item for every pair $(x,y)$ of points of~$U$, there exists a unique closed 
subset $[x,y]$ in $U$ that is homeomorphic to a closed segment with 
boundary points~$x$ and~$y$\;;
\item the boundary of~$U$ in~$\Gamma$ is finite\;.
\end{enumerate}
\end{enumerate}
\end{definition}

By~\cite[Th\'eor\`eme~3.5.1]{Duc}, every $K$-analytic curve is a graph.

\begin{definition}\label{def:finitegraph}
We say that a graph~$\Gamma$ is \emph{finite at a point}~$x$ if there 
exists a positive integer~$n$, a neighborhood~$U$ of~$x$, and a 
homeomorphism
\begin{equation}\label{eq : star around x}
\varphi \colon U \longrightarrow \big\{r \exp(2 i k \pi/n) \in \mathbb C,\ r\in \of{[}{0,1}{[}, k\in\{0,\dotsc,n-1\}\big\}\
\end{equation}
such that $\varphi(x)=0$. The integer~$n$ is independent of the choices of 
$U$ and the homeomorphism $\varphi$.
 We call it the arity of~$\Gamma$ at~$x$ and denote it by~$a_{\Gamma}(x)$.

We say that a point~$x$ in a graph~$\Gamma$ is an end-point (resp. a bifurcation point) if $a_{\Gamma}(x)=1$ (resp. $a_{\Gamma}(x) \ge 3$.

We say that a graph is \emph{locally finite} if it is finite at each of its points. We say that it is \emph{finite} if it is moreover compact. We say that it is \emph{quasi-finite} if it is homeomorphic to the complement of a finite number of points in a finite graph.

\end{definition}


\begin{lemma}
Let~$\Gamma$ be a graph.
\begin{enumerate}
\item The graph~$\Gamma$ is finite if, and only if, it is a finite union of closed segments.
\item The graph~$\Gamma$ is quasi-finite if, and only if, it is a finite union of open or closed segments.
\end{enumerate}
\hfill$\Box$
\end{lemma}

\subsection{Affine line, disks and annuli.}

Let~$L$ be a complete valued extension of~$K$ and let $c\in L$. We denote by~$\E{1}{L}$ the Berkovich affine line in the sense of~\cite{Ber}. Let us fix a coordinate~$T$ on it. We set 
\begin{eqnarray}
D_{L}^+(c,R) &\;=\;& 
\big\{x\in \E{1}{L}\, \big|\, |(T-c)(x)|\le R\big\}\;,\quad\qquad\qquad 
R> 0\\
D_{L}^-(c,R) &\;=\;& \big\{x\in \E{1}{L}\, \big|\, |(T-c)(x)|<R\big\}
\;,\qquad\qquad \quad R>0\\
C_{L}^+(c;R_{1},R_{2}) &\;=\;& 
\big\{x\in \E{1}{L}\, \big|\, R_{1}\le |(T-c)(x)|\le R_{2}\big\}\;, 
\qquad R_2\geq R_1> 0\\
C_{L}^-(c;R_{1},R_{2}) &\;=\;& \big\{x\in \E{1}{L}\, \big|\, 
R_{1} < |(T-c)(x)| < R_{2}\big\}\;.\qquad R_2>R_1> 0
\label{eq : punctured disk is an annulus}
\end{eqnarray} 

\begin{definition}\label{def:modulusannulus}
Let $C$ be an open or closed annulus. If $C= C^\pm(c;R_{1},R_{2})$ in some coordinate, we define the \emph{modulus} of~$C$ to be
\[\mathrm{Mod}(C) := \frac{R_{2}}{R_{1}} \in [1,+\infty].\]
It is independent of the chosen coordinate.
\end{definition}

\subsection{Pseudo-triangulations.}

Let us recall the definition of virtual open disks and annuli (see~\cite[3.6.32 and 3.6.35]{Duc}). 

\begin{definition}
A non-empty connected $K$-analytic space is called a \emph{virtual open 
disk} (resp. \emph{open annulus}) if it becomes isomorphic to a disjoint 
union of open disks (resp. open annuli whose orientations are preserved by 
$\mathrm{Gal}(K^\alg/K))$ over~$\wKa$.
\end{definition}

\begin{definition}\label{def:modulusvirtualannulus}
The \emph{modulus} of a virtual open annulus~$C$ is the modulus of any connected component of~$C_{\wKa}$. We denote it by $\mathrm{Mod}(C)$.
\end{definition}

Following~\cite[5.1.8]{Duc}, we may now define the analytic skeleton of an analytic curve. 

\begin{definition}\label{def:analyticskeleton}
We call \emph{analytic skeleton} of an analytic curve 
$X$ the set of points that have no neighborhoods 
isomorphic to a virtual open disk. We usually denote it by $\Gamma_X$.
\end{definition}

The analytic skeleton of a virtual open disk is empty. That of a virtual open annulus~$C$ is homeomorphic to an open segment~$\Gamma_{C}$ and there is a canonical proper strong deformation retraction $r_{C} \colon C \to \Gamma_{C}$.

Let us now recall the definition of triangulation in the sense of~\cite[5.1.13]{Duc}.

\begin{definition}[Triangulation]\label{def:triangulation}
Let $X$ be a $K$-analytic curve.
A \emph{triangulation} of~$X$ is a locally finite subset 
$S\subseteq X$, formed by points of type~2 or~3, such 
that every connected component of $X-S$ is a relatively 
compact virtual open disk or a relatively compact virtual 
open annulus. 

The \emph{skeleton}~$\Gamma_{S}$ of a 
triangulation~$S$ is the union of~$S$ with the analytic 
skeletons of the connected components of~$X-S$ that 
are virtual open annuli.
\end{definition}

For each triangulation~$S$ of a $K$-analytic curve~$X$, consider the map $r_{S} \colon X \to \Gamma_{S}$ that induces the identity on~$\Gamma_{S}$ and sends each point of a connected component of $X-\Gamma_{S}$ (which is necessarily a relatively compact virtual open disk) to its boundary point (in~$\Gamma_{S})$. It is  a proper map and a strong deformation retraction. 

The 
main result of~\cite{Duc} claims that every 
quasi-smooth $K$-analytic curve may be endowed with a 
triangulation (see~\cite[Th\'eor\`eme~5.1.14]{Duc}). 
This is equivalent to the semistable reduction theorem.


In~\cite{NP-II}, we have used a weaker notion (called \textit{weak 
triangulation}) where we did not require relative compactness of the connected components of $X-S$. 
This allows open disks with empty weak triangulation, hence a 
more adapted framework to study the radii of convergence of differential 
equations. Nevertheless, in some situations this definition forces 
some ``\emph{topologically finite}'' skeletons~$\Gamma_S$ 
to have infinitely many edges 
(that is to say, $\Gamma_S$ contains infinitely many points of~$S$).
For instance, a weak triangulation of an 
open disk $D$ with a rational point $x$ removed is necessarily an  infinite set, 
even in the case where $\Gamma_S$ coincides with an open segment. 
This happens because $D-\{x\}$ is not an open annulus, hence its skeleton is not 
considered a possible edge, by definition. 

In order to have a better notion of edges, 
it will be convenient to use an even weaker notion of triangulation, 
that we now introduce. We will need a more general 
notion of disks and annuli.

Let us first give the definition of the genus of a point (see~\cite[4.4.1 and 4.5.2]{Duc}).
Recall that, for every point~$x$ of type~2 in a $K$-analytic curve, there exists a unique normal irreducible projective curve~$\Cs_{x}$ over~$\wti K$ such that $\wti{\Hs(x)}$ is isomorphic to the function field $\wti K(\Cs_{x})$ of~$\Cs_{x}$. It is called the residual curve of~$x$ (see~\cite[3.3.5.2]{Duc}). We denote by $g(\Cs_x)$ its genus.

\begin{definition}[Genus of a point]\label{Def : genus of a point}
Let $X$ be a $K$-analytic curve. Let~$x\in X$.

Assume that~$K$ is algebraically closed. We define the genus of~$x$ by 
\[g(x) = 
\begin{cases}
g(\Cs_{x}) \textrm{ if } x \textrm{ has type } 2;\medskip\\
0 \textrm{ otherwise.}
\end{cases}\]

If~$K$ is arbitrary, the genus~$g(x)$ of a point $x\in X$ is defined to be the sum of the genuses of the points of $X_{\wKa}$ above~$x$.
\end{definition}

\begin{lemma}
In the situation of Definition~\ref{Def : genus of a point}, let $x_{1},\dotsc,x_{d}$ be the preimages of~$x$ in~$X_{\wKa}$. We have
\begin{equation}
g(x_{1}) = \dotsb = g(x_{d}) = \frac1d \, g(x).
\end{equation}
\hfill$\Box$
\end{lemma}

\begin{definition}[Open pseudo-annulus]
\label{Def : Pseudo-annulus}
Assume that~$K$ is algebraically closed. We say that a connected quasi-smooth $K$-analytic curve~$C$ is an 
\emph{open pseudo-annulus} if
\begin{enumerate}
\item it has no boundary;
\item it contains no points of positive genus;
\item its analytic skeleton~$\Gamma_{C}$ is homeomorphic to 
a non-empty open segment.
\end{enumerate}  
We call~$\Gamma_{C}$ the skeleton of~$C$.

If~$K$ is arbitrary, we say that a connected quasi-smooth $K$-analytic curve~$C$ is an 
\emph{open pseudo-annulus} if $C \ho_{K} \widehat{K^\alg}$ is a (finite)\footnote{Recall that every connected $K$-analytic curve $X$ has the 
property that $X_{\wKa}$ has finitely many connected components. 
This is due to the fact that the Galois orbit of a point of type $2$ or $3$ is 
always finite.} disjoint union of open pseudo-annuli and if 
$\mathrm{Gal}(K^\alg/K)$ preserves the orientation of their skeletons. 
In this case, one can check that~$C$ has no boundary, contains no points 
of positive genus and that its analytic skeleton~$\Gamma_{C}$ is an open 
segment. We call $\Gamma_{C}$ the skeleton of~$C$.
\end{definition}

\begin{remark}\label{Rk : union of virtual ann = pseudo ann}
A increasing union of virtual open annuli with nested skeletons is an open 
pseudo-annulus. We will later prove that the converse holds (see 
Proposition~\ref{prop:pseudo-annulus}).
\end{remark}

\begin{remark}\label{Rk : pseudo-annuli}
By \cite[Proposition 3.2]{Liu},  if~$K$ is non-trivially valued, 
algebraically closed and maximally complete, any open 
pseudo-annulus~$C$ may be embedded into the affine 
line. We deduce that, in this case, either~$C$ is an 
open annulus or it may be written in the form 
$C=Y-\{y\}$, where $Y$ is either the affine line or an open 
disk and $y\in Y$ is a rational point.
\end{remark}

\begin{definition}
Let~$C$ be an open pseudo-annulus. We define the degree of~$C$ as the 
number of connected components of $C_{\wKa}$
\begin{equation}
\deg(C) \; := \; \Card(\pi_{0}(C_{\wKa}))\;.
\end{equation}
\end{definition}

\begin{lemma}\label{lem:ccpseudoannulus}
Let~$C$ be an open pseudo-annulus. Each connected component of $C - \Gamma_{C}$ is a virtual open disk.
\end{lemma}
\begin{proof}
A graph with the claimed 
property is called analytically admissible in \cite[(5.1.3)]{Duc}.
The result then follows from \cite[(5.1.11)]{Duc}.
\end{proof}


Similarly to the notion of open pseudo-annulus generalizing the notion of annulus, there is a notion of pseudo-disk generalizing that of a disk. We will not use it later.

\begin{definition}[Open pseudo-disk]
\label{Def : Pseudo-disk}
We say that a connected quasi-smooth $K$-analytic curve~$D$ is an \emph{open pseudo-disk} if it is not compact and if its analytic skeleton is empty.
\end{definition}

We note that if $D$ is a pseudo-disk, then $D_{\widehat{K^{\mathrm{alg}}}}$ is a finite disjoint union of open 
pseudo-disks over $\widehat{K^{\mathrm{alg}}}$.

\begin{remark}
A increasing union of virtual open disks is an open pseudo-disk. We will later prove that the converse holds (see Proposition~\ref{prop:pseudo-disk}).
\end{remark}

\begin{remark}\label{Rk : pseudo-disk}
By \cite[Prop. 3.2]{Liu},  if~$K$ is non-trivially valued, algebraically closed, and spherically complete, any open pseudo-disk~$D$ may be embedded into the affine line. We deduce that, in this case, $D$ is either an open disk or the affine line.
\end{remark}

\begin{definition}[Pseudo-triangulation]
\label{def:pseudo-triangulation}
Let $X$ be a $K$-analytic curve. 
A \emph{pseudo-triangulation} of~$X$ is a locally finite 
subset $S\subseteq X$, formed by points of type~2 
or~3, such that every connected component of $X-S$ is 
a virtual open disk or an open pseudo-annulus. 

The \emph{skeleton}~$\Gamma_{S}$ of a pseudo-triangulation~$S$ is the union of~$S$ with the 
skeletons of the connected components of~$X-S$ that are open pseudo-annuli. 
\end{definition}

\begin{remark}\label{Remark : pseudo = weak} 
Let $X$ be a connected $K$-analytic curve endowed with a pseudo-triangulation~$S$.
\begin{enumerate}
\item If $\Gamma_{S}=\emptyset$ (hence $S=\emptyset$), then~$X$ is a virtual open disk.
\item If $S=\emptyset$ and $\Gamma_{S}\ne\emptyset$, then $X$ is an open pseudo-annulus.
\end{enumerate}
\end{remark}

\begin{proposition}\label{prop:GammaSS'}
Let $X$ be a connected $K$-analytic curve endowed with a pseudo-triangulation~$S$. Assume that $\Gamma_{S}\ne\emptyset$. 
\begin{enumerate}
\item\label{prop:GammaSS'-i} Each connected component of $X-\Gamma_{S}$ is a virtual open disk $D$ which is relatively compact in $X$ and whose relative boundary $x_D$ in $X$ belongs to $\Gamma_S$. 
\item\label{prop:GammaSS'-ii} The map 
\begin{equation}\label{eq : retraction prop:GammaSS'}
r_{S} \;\colon\; X \to \Gamma_{S}
\end{equation}
 which is the identity on $\Gamma_S$ and maps the virtual open disks $D$ from item~\eqref{prop:GammaSS'-i} to their boundary points $x_D$ is a proper continuous retraction.
\item\label{prop:GammaSS'-iii} There exists a triangulation $S'$ (in the strong sense of Definition~\ref{def:triangulation}) such that $\Gamma_{S'} = \Gamma_{S}$.
\end{enumerate}
\end{proposition}
\begin{proof}
Item~\eqref{prop:GammaSS'-i} follows from Lemma~\ref{lem:ccpseudoannulus}.

Item~\eqref{prop:GammaSS'-ii} follows from \cite[(1.5.16)]{Duc}.

Item~\eqref{prop:GammaSS'-iii} follows from \cite[(5.1.14)]{Duc} since, by 
item~\eqref{prop:GammaSS'-i}, $\Gamma_{S}$ is analytically admissible (cf. \cite[(5.1.3)]{Duc}).
\end{proof}

We finally add a notion of controlling graph of a function (cf. \cite[Definition 2.4.2]{NP-III}).


\begin{definition}\label{def:controllinggraph}
Let $X$ be a connected $K$-analytic curve endowed with a 
pseudo-triangulation~$S$ and $f \colon X\to\ERRE$ be a function.
We call $S$-\emph{controlling graph} (or $S$-\emph{skeleton}) of~$f$ the 
set $\Gamma_{S}(f)$ of points $x\in X$ that admit no 
neighborhoods~$D_{x}$ in $X$ such that
\begin{enumerate}
\item\label{def:controllinggraph-i} $D_{x}$ is a virtual open disk;
\item\label{def:controllinggraph-ii} $f$ is constant on $D_{x}$;
\item\label{def:controllinggraph-iii} $D_{x}\cap \Gamma_S=\emptyset$ (or equivalently 
$D_{x}\cap S=\emptyset$).
\end{enumerate}
\end{definition}

\begin{lemma}\label{lem:controllinggraph}
Let $X$ be a connected $K$-analytic curve endowed with a pseudo-triangulation~$S$ and $f \colon X\to\ERRE$ be a function. We have $\Gamma_S\subseteq\Gamma_S(f)$.

Let $D$ be a connected component of~$X - \Gamma_{S}$ (hence $D$ is a virtual open disk). Let $x \in \Gamma_{S}(f) \cap D$. Then, $\Gamma_{S}(f)$ contains the half-open segment joining~$x$ to the boundary of~$D$. 

In particular, each connected component of $\Gamma_{S}(f) - \Gamma_{S}$ is a tree. 
\end{lemma}
\begin{proof}
The fact that $\Gamma_S\subseteq\Gamma_S(f)$ follows from item \eqref{def:controllinggraph-iii} of Definition~\ref{def:controllinggraph}.

Let us prove the second part of the statement. Let $y$ be a point of the segment joining~$x$ to the boundary of~$D$. If $y \notin \Gamma_{S}(f)$, then, there exists a virtual open disk~$D_{y}$ in $X$ containing~$y$ such that $D_{y}\cap \Gamma_S=\emptyset$ and $f$ is constant~$D_{y}$. Since $D_{y}\cap \Gamma_S=\emptyset$, we have $D_{y} \subseteq D$, hence $D_{y}$ contains~$x$, and we get a contradiction.
\end{proof}


\subsection{Segments and germs of segment.}
\label{section branch}

Let $X$ be a $K$-analytic curve. 

Recall that a map between topological spaces is called a \emph{topological embedding} if it induces a homeomorphism onto its image. 

We define a \emph{closed (resp. open, resp. semi-open) segment} in~$X$ to be the image of a topological embedding from a closed (resp. open, resp. semi-open) real interval to~$X$.

Assume given a topological embedding $\varphi \colon [0,1] \to X$ with $\varphi(0)=x$ and $\varphi(1)=y$. In this case, when no confusion may occur, we write $[x,y] = \varphi([0,1])$, $\of{[}{x,y}{[}\ = \varphi(\of{[}{0,1}{[})$, $\of{]}{x,y}{]}\ = \varphi(\of{]}{0,1}{]})$ and $\of{]}{x,y}{[}\ = \varphi(\of{]}{0,1}{[})$.



Let~$I$ be an open segment in $X$. Let~$\Ps_{I}$ be the set of topological embeddings $\varphi \colon \of{]}{0,1}{[} \to X$ whose image is~$I$. We 
define an equivalence relation~$\sim$ on~$\Ps_{I}$ by decreeing that 
$\varphi \sim \psi$ if $\varphi \circ \psi^{-1}$ is increasing. There are two 
equivalence classes, which we call \emph{orientations} of~$I$.

Consider the set~$\Ic_{oo}$ of oriented open segments in~$X$. We define an equivalence relation~$\sim_{oo}$ on~$\Ic_{oo}$ by decreeing that $(I_{1},\varphi_{1}) \sim_{oo} (I_{2},\varphi_{2})$ if there exists $\alpha_{1},\alpha_{2}\in \of{]}{0,1}{[}$ such that $\varphi_{1}(\of{]}{0,\alpha_{1}}{[}) = \varphi_{2}(\of{]}{0,\alpha_{2}}{[})$ and the map $\varphi_{2}^{-1} \circ \varphi_{1} \colon \of{]}{0,\alpha_{1}}{[} \to \of{]}{0,\alpha_{2}}{[}$ is increasing.

\begin{definition}
A \emph{germ of segment} is an equivalence class of oriented open segments, \textit{i.e.} an element of $\Ic_{oo}/\sim_{oo}$.

A germ of segment~$b$ is said to be \emph{relatively compact} if one of 
its representative is. Then, every representative has one or two boundary 
points and all representatives have exactly one in common, which is called 
\emph{boundary point} of the germ of segment. In this case, $b$ is said to 
be a \emph{germ of segment out of this boundary point}. 
This corresponds to A.~Ducros's notion of branch (see~Remark \ref{Rk : Branch vs germs}).

We call \emph{germ of segment at infinity} every non-relatively-compact 
germ of segment. 
We call 
\emph{open boundary} of~$X$, and denote by~
\begin{equation}
\partial^o X\;,
\end{equation} the set of those segments. We denote by~$N_{\infty}(X)$ the cardinal of the open boundary of~$X_{\wKa}$.
\end{definition}


For instance, if~$X$ is a virtual open disk, then 
$N_{\infty}(X) = \Card(\pi_{0}(X_{\wKa}))$, while 
if it is an virtual open annulus, then 
$N_{\infty}(X)=2\cdot \Card(\pi_{0}(X_{\wKa}))$.

\begin{remark}\label{Rk: Orientation}
It follows from the definition that if $b$~is a germ of 
segment out of a point (i.e. if $b$ is relatively compact), 
then $b$ is oriented away from that 
point. On the other hand, if $b$ lies at the open 
boundary of $X$, then it is oriented towards the interior of~$X$.   
\end{remark}

%

\begin{lemma}\label{Lemma : Top finite >> N(X)<inf}
Assume that~$X$ is connected. Let~$S$ be a pseudo-triangulation of~$X$ 
such that $\Gamma_{S} \ne \emptyset$. Then every germ of segment at 
infinity of~$X$ is contained in~$\Gamma_{S}$.
\end{lemma}
\begin{proof}
Let~$b$ be a germ of segment at infinity. Let $\varphi \colon \of{]}{0,1}{[} \to X$ be a topological embedding whose image~$I$ represents~$b$. 

Let us first assume that~$I$ does not meet~$\Gamma_{S}$. Then it belongs to some connected component~$D$ of $X - \Gamma_{S}$. By assumption, $D$ is a virtual open disk. If it had no boundary point in~$X$, then it would be closed in~$X$. Since~$X$ is connected, we would have~$D=X$, hence~$\Gamma_{S}=\emptyset$, which contradicts the assumptions of the lemma. We deduce that~$D$ is relatively compact in~$X$, hence so is~$I$, contradicting the fact that~$b$ is a germ of segment at infinity. 

We have just proved that there exists~$u \in \of{]}{0,1}{[}$ such that $\varphi(u) \in \Gamma_{S}$. Assume, by contradiction, that there exists $v<u$ such that $\varphi(v) \notin \Gamma_{S}$. Then~$\varphi(v)$ belongs to some connected component~$D$ of $X - \Gamma_{S}$. Set $w = \sup(\{w'\in \of{[}{v,u}{]} \mid \varphi(w')\in D\})$. Then~$\varphi(w)$ belongs to the boundary of~$D$ in~$X$. Since this boundary is reduced to a point and~$\varphi$ is injective, for every~$w'<w$, we have $\varphi(w')\in D$. 
%
It follows that we may find a representative of~$b$ whose image lies in $D$. It is then relatively compact, and we get a contradiction.
\end{proof}

\begin{lemma}\label{lem:bordvide}
Assume that~$X$ is connected. Let~$S$ be a pseudo-triangulation of~$X$. The following are equivalent:
\begin{enumerate}
\item $X$ is compact;
\item $\Gamma_{S}$ is non-empty and compact (\emph{i.e.} a finite graph);
\item $\partial^oX=\emptyset$.
\end{enumerate}
\end{lemma}
\begin{proof}
The implications $i) \implies ii)$ and $i) \implies iii)$ are clear (cf. Lemma \ref{Lemma : Top finite >> N(X)<inf})

The implication $ii) \implies i)$ follows from the properness of the retraction map $r_{S} \colon X \to \Gamma_{S}$ (see Proposition~\ref{prop:GammaSS'}).

It remains to prove $iii) \implies ii)$. Assume that $\partial^oX=\emptyset$. 

If $\Gamma_{S} = \emptyset$, then $X$ is a virtual open disk, which contradicts the assumption. It follows that $\Gamma_{S} \ne \emptyset$. By Proposition~\ref{prop:GammaSS'}, we may assume that $S$ is a triangulation of~$X$. In this case, $\Gamma_{S}$ is a graph in the usual sense (defined by vertices and edges between them), with vertex set~$S$.

Assume, by contradiction, that $\Gamma_{S}$ is not compact. 
Since~$\Gamma_{S}$ is locally finite, by K\"onig's lemma, $\Gamma_{S}$ 
contains an infinite path, that is to say a segment going through infinitely 
many elements of~$S$. Since~$S$ is locally finite, this segment may not 
be relatively compact, which contradicts the fact that 
$\partial^oX=\emptyset$.

%
\end{proof}

\begin{remark}\label{Rk : Branch vs germs}
The notion of branch in \cite[Section~1.7]{Duc} corresponds to that of 
relative compacts germs of segments.
Indeed, a Ducro's branch is by 
definition associated with a base point $x$:  
over $\wKa$, a branch can be seen as a connected 
components of $U-\{x\}$ where $U$ is an unspecified connected 
neighborhood of $x$ which is a tree (no loops). 
Our notion of germ of segment, instead, is not necessarily 
associated with a base point. 
In this sense, a Ducro's branch corresponds to a 
relatively compact germ of segments.
\end{remark}

\begin{remark}\label{rk : bad examples}
In Lemma \ref{lem:preimage} we will describe precisely the germs of 
segments that can be represented by the analytic skeleton of a 
pseudo-annulus. Let us mention here a couple of examples that can not be 
represented in that way. Let $X$ be the open unit disk.

\begin{enumerate}
\item\label{rk : bad examples-i} Assume $K$ algebraically closed. Let 
$Z\subseteq X$ the set of zeros of an unbounded analytic function on $X$. 
In this case the set $Z$ is infinite and locally finite. 
The germ of segment at infinity of $X$ produces a germ of segment in 
$X-Z$ which can not be represented as the skeleton of a pseudo-annulus 
in $X-Z$.
\item\label{rk : bad examples-ii} Assume that $K$ is not algebraically closed. Let $t\in \in X$ be a 
non-rigid point of type $1$. Since degree of $t$ is infinite 
(cf. Definition \ref{def:degree}), every segment out of $t$ 
has a pre-image in $X_{\wKa}$ with infinitely many bifurcations, while the 
analytic skeleton of a pseudo-annulus has no bifurcations.
\end{enumerate}
\end{remark}

\subsection{Slopes and Laplacian.}
Recall that it is possible to define a canonical length~$\ell(I)$ for any closed segment~$I$ inside a $K$-analytic curve (see~\cite[Proposition 4.5.7 and Corollaire 4.5.8]{Duc} or \cite[beginning of Section~3.1]{NP-II}). The fundamental example is that of the skeleton~$\Gamma_{C}$ of a virtual open annulus~$C$, in which case the length is defined as (cf. Definitions 
\ref{def:modulusannulus} and \ref{def:modulusvirtualannulus})
\[ \ell(\Gamma_{C}) = \log(\mathrm{Mod}(C)).\]




\medbreak

Let $X$ be a $K$-analytic curve. 

\begin{definition}\label{Def : partial_b}
Let $(I,\varphi)$ be an oriented open segment in~$X$. We say that a map $F \colon X \to \ERRE_{>0}$ is \emph{log-affine} on~$I$ if there exists $\lambda\in \ERRE$ such that, for each $u < v \in \of{]}{0,1}{[}$, we have 
\begin{equation}
\log(F(\varphi(v))) - \log(F(\varphi(u))) = \lambda \, \ell([\varphi(u),\varphi(v)]).
\end{equation}
%

In this case, we define the \emph{slope of~$F$ along~$I$} to be~$\lambda$. We denote it by~$\partial_I F$. The property of being log-affine and the slope are preserved by equivalence of oriented open segments.

Let~$b$ be a germ of segment. We say that a map $F \colon X \to \ERRE_{>0}$ is \emph{log-affine} on~$b$ if  it is log-affine on some representative~$I$ of~$b$ (or equivalently on every representative of~$b$ that is small enough). In this case, we define the \emph{slope of~$F$ along~$b$} to be the slope of~$F$ on~$I$ (or equivalently on every representative of~$b$ that is small enough). We denote it by~$\partial_b F$. 

%
\end{definition}


\begin{lemma}
Let~$L$ be a complete valued extension of~$K$. Let~$b$ be a germ of 
segment in~$X$ and $F \colon X \to \ERRE$ be a map that is log-affine 
on~$b$. 

Let~$c$ be a germ of segment in~$X_{L}$ above~$b$. Then the map 
$F_{L} := F \circ \pi_{L} \colon X_{L} \to \ERRE$ is log-affine on~$c$ and we 
have
\begin{equation} 
\partial_c F_{L} = \partial_{b} F\;.
\end{equation}
\hfill$\Box$
\end{lemma}

\begin{definition}\label{def:degreegerm}
Let~$b$ be a germ of segment in~$X$. Let $c$ be a germ of segment in~$X_{\wKa}$. We say that \emph{$c$ lies above~$b$} and we write $\pi(c) = b$ if there exists a representative $(I,\varphi)$ of~$c$ such that $(\pi_{\wKa}(I),\varphi \circ \pi_{\wKa})$ is an oriented open segment in~$X$ whose equivalence class is~$b$.

We denote by $\pi^{-1}_{\wKa}(b)$ the set of germs of segments in~$X_{\wKa}$ above~$b$. 

We define the \emph{degree} of~$b$ as
\begin{equation}
\deg(b) \; := \; \Card(\pi^{-1}_{\wKa}(b))\;.
\end{equation}

\if{
\comm{
$\pi:X_{\wKa}\to X$ induit bien une fonction envoyant un germe de 
segment dans $X_{\wKa}$ vers un germe dans $X$. 
Et donc ça a bien du sens de parler d'antécedent d'un germe. 
Et $\deg(b)$ est défini correctement.\\ 

Toutefois, il y peut y avoir ambiguité quand l'on travaille avec des 
antécédents en prenant des représentants. 
Et la notion de "disjoint union of germs of segments" 
a plus d'un sens possible.\\

Notamment , on a le problème suivant.\\

\textbf{PROBLEME :} On a deux possibilités pour 
"disjoint union" 

\begin{itemize}
\item pour tout $b$ et toute paire de germes $b',b''$ au dessus de $b$ 
on peut choisir des representants $I'$ et $I''$ de $b',b''$ de sorte que 
$I'\cap I'' = \emptyset$ 
\item  pour chaque $b'$ au dessus de $b$ on choisit un representant $I'$
de sorte que tous ces representants soient disjoints simultanément
\end{itemize}
Il y a différence ... dans le premier cas on choisit $b$ et 
$b'$ et on adapte ensuite le choix de $I'$ et $I''$, alors que dans le 
deuxième cas on choisit d'abord tous les representants des $b'$ et on 
affirme qu'on peut tous les choisir d'un seul coup 
de sorte que ce soit une réunion disjointe.\\

\textbf{VRAI PROBLEME DANS UN CAS CONCRET :} Soit $X$ le 
disque unité centré en $0$, et je considère un point de Berkovich de type $1$ donnée par une orbite infinie d'un point $t\in\wKa-K^{\mathrm{alg}}$, 
avec $|t|<1$. 

\begin{itemize}
\item Soit $\{t_i\}_{i}$ l'orbite de $t$; 

\item $D_i$ un disque ouvert centré en $t_i$; 

\item $C_i=D_i-\{t_i\}$ la pseudocouronne dont le squelette $I_i$ represente la branche sortante de $t_i$.
\end{itemize}

Pour chaque couple $t_1,t_2$ d'antécedents dans la fibre, je peux trouver 
des germes de segments disjoints, mais c'est pas sur qu'on puisse faire un 
choix global de sorte qu'ils soient tous disjoints.\\

Notamment, voici ce que ça donnerait : 
%
La réunion des $D_i$ recouvre notre orbite $\{t_i\}_i$, 
qui est compacte, et donc il existe un sous-recouvrement fini. 
Il en suit qu'il existe $i,j$ tels que $D_i\cap D_j\neq \emptyset$. 
Il faut donc faire le choix du rayon  $r_i$ de $D_i$ de sorte que 
$r_i<|t_i-t_j|$ autrement les segments $I_i$ et $I_j$ s'intersectent. 
Je crois que c'est possible si l'on suppose la fibre dénombrable (dans ce 
cas on établit un ordre sur l'ensemble des indices et on fait un choix de 
proche en proche), mais en général je ne sais pas faire.

\begin{remark}\label{Rk discr avec Ducros} On a peut être une discrépance de définitions entre nous et 
Ducros ici. Car pour nous un germe de segment est un truc vraiment 
homeomorphe à un segment, alors que pour Ducros c'est la composante 
connexe toute entière (ou une limite de comosantes connexes). 
Quelle serait la définition de Ducros dans l'exemple précédent ? 
Je crois que la "branche" de Ducros est $C_i$, modulo faire 
une classe d'équivalence qui consiste à s'autoriser à retrecir le rayon 
exterieur de $C_i$. Or, avec cette définition, quand $D_i\subseteq D_j$ on a 
forcement l'intersection $C_i\cap C_j$ non vide. 
Il en résulte qu'il est impossible de choisir une famille de representants 
$C_i$ qui soient tous disjoints. \\

\textbf{Donc la définition de Ducros ne permet pas de réunion disjointe.}\\

\textbf{Par contre, notre définition permet une réunion disjointe si la fibre 
est dénombrable.}
\end{remark}
}
}\fi
\end{definition}
\begin{remark}
The degree of a germ of segment~$b$ may be infinite. If~$b$ may be 
represented by the skeleton of an open pseudo-annulus~$C$, 
for instance if~$b$ is a germ of segment out of a point of type~2 or~3 or 
a rigid point, then its degree is finite (and equal to the degree of~$C$).
\end{remark}


\begin{remark}
It follows from the definitions that we have 
\begin{equation}
N_{\infty}(X) = \sum_{b \in \partial^o X} \deg(b).
\end{equation}
\end{remark}

\begin{definition}[Laplacian]
Let~$x\in X$. Let $F \colon X \to \ERRE$ be a map that is $\log$-affine on every germ of segment out of~$x$ and constant on all but finitely many of them. We define the Laplacian of~$F$ at~$x$ as
\begin{equation}
dd^c F(x) \; := \; \sum_{b} \deg(b)\, \partial_{b} F\;,
\end{equation}
where~$b$ runs through the family of germs of segment out of~$x$.
\end{definition}

%

\begin{definition}\label{def:degree}
Let~$x\in X$. 
We define the degree of~$x$ as
\begin{equation}
\deg(x) \; := \; \Card(\pi^{-1}_{\wKa}(x))\;.
\end{equation}
\end{definition}

\begin{remark}
The degree of a point may be infinite. It is automatically finite for points of 
type~2 and~3 and for rigid points. 

\end{remark}

\begin{remark}\label{rem:s}
Let $x\in X$. We may consider the separable (algebraic in our case) closure~$\mathfrak{s}(x)$ of~$K$ in~$\Hs(x)$ (see~\cite[3.1.1.5]{Duc}). We then have $[\mathfrak{s}(x) \colon K] = \deg(x)$.
\end{remark}

We record here some properties of the degrees in the case of virtual open annuli (see \cite[(3.6.35.5)]{Duc}) and virtual open disks (\cite[(3.6.36.2), (4.5.14.2)]{Duc}).

\begin{proposition}\label{prop:degannulus}
Let $C$ be a virtual open annulus. Then, for each $x\in \Gamma_{C}$, we have $\deg(x) = \deg(C)$. 

In particular, for each virtual open annulus~$C'$ contained in~$C$ with $\Gamma_{C'} \subseteq \Gamma_{C}$, we have $\deg(C') = \deg(C)$. 
\qed
\end{proposition}

\begin{proposition}\label{prop:degdisk}
Let $D$ be a virtual open disk. 
Let $I$ be a semi-open segment in $D$ 
with end-point~$x\in D$ of type~2 or~3. We assume that $I$ 
is not relatively-compact in~$D$, 
\textit{i.e.} $I - \{x\}$  represents the open boundary of~$D$. 
Then, the map $y \in I \mapsto \deg(y)$ is non-increasing in the direction of the open boundary of~$I$, and eventually constant with value $\deg(D)$.
\qed
\end{proposition}

%
%

\subsection{Properties of pseudo-disks and pseudo-annuli.}

We now give some more properties of pseudo-disks and pseudo-annuli.

%

\begin{lemma}\label{lem:X=Y}
Let $X$ be a connected analytic curve and let $Y$ be an analytic domain of~$X$. Assume that 
\begin{enumerate}
\item $Y$ has no boundary;
\item $Y$ is not relatively compact in~$X$;
\item $Y$ has exactly one germ of segment at infinity.
\end{enumerate} 
Then $Y=X$.
\end{lemma}
\begin{proof}
Let us first remark that $Y$ is connected. Indeed, otherwise $Y$ would have a connected component~$C$ with no germ of segment at infinity. By Lemma~\ref{lem:bordvide}, $C$ is compact, hence closed in~$X$. Since~$Y$ has no boundary, $C$ is also open in~$X$. By connectedness of~$X$, we deduce that $C=X$, which is a contradiction since $C$ is strictly contained in~$Y$.

Fix $y\in Y$. Let $x\in X$. Since $X$ is connected, it is path-connected, so there exists $\varphi \colon [0,1] \to X$ such that $\varphi(0) = y$ and $\varphi(1) = x$.

Let $S$ be the subset of elements~$s$ in~$[0,1]$ such that $\varphi([0,s]) \subseteq Y$. Set $s_{0} := \sup(S) \in [0,1]$. Since~$Y$ has no boundary, it contains a neighborhood of~$y$, hence we have $s_{0}>0$. In particular, $I_{0} := \varphi(\of{]}{0,s_{0}}{[})$ is an open segment.

Assume, by contradiction, that $\varphi(s_{0}) \notin Y$. Then, the segment $I_{0}$ is not relatively compact in~$Y$, so it represents the unique germ at infinity of~$Y$. Let $V$ be a connected affinoid neighborhood of~$\varphi(s_{0})$ in~$X$. Note that $Y\cup V$ is connected. We will prove that $Y\cup V$ is compact, which contradicts hypothesis~ii). By Lemma~\ref{lem:bordvide}, it is enough to prove that $\partial^\circ (Y\cup V) = \emptyset$. 

Let $b$ be a germ of segment in~$Y\cup V$. Remark that we have $b\subseteq V$ or $b\subseteq Y$. Indeed, if neither property holds, then any representative~$I_{b}$ of~$b$ would pass several times from~$Y$ to~$V$, hence $I_{b}$ would pass several times through the open boundary of~$Y$. Since $Y$ has only one germ of segment at infinity, this would contradict the injectivity of the map defining~$I_{b}$. 

Let us now prove that $b$ is relatively compact in $Y\cup V$. If $b\subseteq V$, this follows from the compactness of~$b$. If $b\subseteq Y$ and $b \notin \partial^\circ Y$, the result holds by definition. If $b \in \partial^\circ Y$, then $b$ is represented by~$I_{0}$, and the latter in relatively compact in $Y\cup V$, since it contains~$\varphi(s_{0})$. We deduce that $\partial^\circ (Y\cup V) = \emptyset$, as claimed.


%

We have proved that $\varphi(s_{0}) \in Y$. Since $Y$ has no boundary, it contains a neighborhood of~$\varphi(s_{0})$, hence $s_{0}$ cannot be the supremum of~$S$ unless $s_{0}=1$. We deduce that $x = \varphi(1) = \varphi(s_{0}) \in Y$.
\end{proof}

\begin{lemma}\label{lem:S-s}
Let $X$ be a connected analytic curve that is not compact. Let $S$ be a triangulation of~$X$ (in the strong sense of Definition \ref{def:triangulation}). Let $s \in S$ be an end-point of~$\Gamma_{S}$. If $s$ does not belong to the analytic skeleton of~$X$, then $S - \{s\}$ is still a triangulation of~$X$.
\end{lemma}
\begin{proof}
We use Ducros's criterion \cite[5.1.16.3]{Duc}. We need to check 3 points.

$\bullet$ Since~$s$ does not belong to the analytic skeleton of~$X$, it belongs to a virtual open disk. In particular, $s$ does not belong to the boundary of~$X$ and its genus is~0.

$\bullet$ Since $s$ is an end-point of~$\Gamma_S$, it belongs to a unique edge $\Gamma_0$ of~$\Gamma_{S}$, of the form $\of{[}{s,s'}{]}$. With Ducros's notation, we have $\Gamma_{0} = \of{[}{s,s'}{[}$, which is a semi-open interval.

$\bullet$ It remains to show that $\mathfrak{s}(s) = \mathfrak{s}(\beta)$, where~$\beta$ is the branch out of~$s$ along~$\Gamma_{0}$. Let~$C$ be the connected component of $X-\{s\}$ containing~$\Gamma_{0}$. Since $s$ is an end-point of~$\Gamma_{S}$, all the connected components of $X-\{s\}$ are virtual open disks, except possibly~$C$. Since~$X$ is not compact, $C$ is not relatively compact in~$X$. In particular, $C$ is not a virtual open disc.

By assumption, $s$ belongs to a virtual open disk~$D_{s}$. Let~$C'$ be the non-relatively compact connected component of $D_{s}-\{s\}$. Let~$C''$ be another connected component of~$D_{s} - \{s\}$. It is a virtual open disk. If $C'' \subseteq C$, then, by Lemma~\ref{lem:X=Y}, we have $C''=C$, and we get a contradiction. 
We deduce that $C' \subseteq C$. 
In particular, the branch~$\beta$ coincides with the branch~$\beta'$ out of~$s$ containing the boundary of~$D_{s}$. By \cite[5.1.16.2 (2nd paragraph)]{Duc}, we have $\mathfrak{s}(s) = \mathfrak{s}(\beta')$. The result follows.
\end{proof}

\begin{proposition}\label{prop:pseudo-disk}
Let~$D$ be an open pseudo-disk. Then the following properties hold.

\begin{enumerate}
\item \label{prop:pseudo-disk-i} The topological space~$D$ is a tree (\textit{i.e.} its topological genus is~0).
\item \label{prop:pseudo-disk-ii} The space~$D$ is a countable increasing union of virtual open disks.
\item \label{prop:pseudo-disk-iii} The space~$D$ has exactly one germ of segment at infinity~$b_{\infty}$ and we have 
\begin{equation}
\deg(b_{\infty}) = N_{\infty}(D) = \Card(\pi_{0}(D_{\wKa})).
\end{equation}
\item \label{prop:pseudo-disk-iv} For every representative~$I$ of~$b_{\infty}$ that it small enough, we have
\begin{equation}
\forall x\in I,\ \deg(x) = \Card(\pi_{0}(D_{\wKa})).
\end{equation} 
\end{enumerate}
\end{proposition}
\begin{proof}
By \cite[(1.5.10)]{Duc}, $D$ is a tree with at most one germ of segment at infinity. Since $D$ is not compact, it has exactly one germ of segment at infinity. 


\medbreak

For the rest of the proof, we fix a triangulation~$S$ of~$D$ (in the strong sense of Definition \ref{def:triangulation}) and consider its skeleton~$\Gamma_{S}$. Since~$D$ is not compact and~$r_{S}$ is proper, the locally finite 
graph~$\Gamma_{S}$ is not finite, hence it contains a 
non-relatively-compact segment~$I$. We may assume 
that~$I$ is semi-open and that its boundary point belongs to~$S$. Denote it by~$s_{0}$.

By \cite[Th\'eor\`eme (4.5.10)]{Duc}, $D$ is paracompact, hence~$S$ is countable. By Lemma~\ref{lem:S-s}, if $s$ is an end-point of~$\Gamma_{S}$, then $S-\{s\}$ is still a triangulation of~$D$. Recall also that, by item~\eqref{prop:pseudo-disk-i}, $\Gamma_{S}$ is a tree. Arguing by induction, we may now prove that there exists a subset~$S'$ of~$S$ containing~$s_{0}$ that is still a triangulation of~$D$ and such that~$\Gamma_{S'} = I$.

Identify~$I$ with $\of{[}{0,1}{[}$ and denote by $0 = s_{0} < s_{1} < \dotsb$ the points of~$S'$. 
Denote by $r\colon D \to I$ the canonical retraction. By Lemma~\ref{lem:S-s}, $S'_{0} :=S'-\{s_{0}\}$ is still a triangulation of~$D$. The preimage $D_{1} := r^{-1}(\of{[}{s_{0},s_{1}}{[})$ is a connected component of $S'_{0}-\{s_{1}\}$ that does not meet~$\Gamma_{S'_{0}}$, hence it is a virtual open disk. Arguing by induction, we prove that, for every $n\ge 0$, the preimage $D_{n} := r^{-1}(\of{[}{s_{0},s_{n}}{[})$ is a virtual open disk. This proves item \eqref{prop:pseudo-disk-ii}.

\medbreak

It remains to prove the statements about the degrees. By Proposition~\ref{prop:degdisk} applied to the~$D_{n}$'s, the map $y \in I \mapsto \deg(y)$ is non-increasing in the direction of the open boundary of~$I$. Since the degree only takes positive integral values, it is eventually constant, with value~$d$. By Proposition~\ref{prop:degdisk} again, for $n$ big enough, we have $\deg(D_{n})=d$. Since $D$ is the increasing union of the $D_{n}$'s, we deduce that $ \Card(\pi_{0}(D_{\wKa})) = d$. This finishes the proof.

\end{proof}

For a given pseudo-triangulation~$S$, we introduce an exceptional set~$E_{S}$ of~$\Gamma_{S}$. It contains the nodes of~$\Gamma_{S}$ in Ducros's terminology (see \cite[(5.1.12)]{Duc}).



\begin{definition}\label{def:exceptionalset}
Let $X$ be a $K$-analytic curve endowed with a pseudo-triangulation~$S$.
The \emph{exceptional set}~$E_{S}$ associated to~$S$ is the set of points~$x\in \Gamma_{S}$ where (at least) one of the following conditions is satisfied:
\begin{enumerate}
\item $x\in \partial X$;
\item $g(x) >0$;
\item $x$ is a bifurcation point of~$\Gamma_{S}$;
\item $x$ is an end-point of~$\Gamma_{S}$;
\item the map $y\mapsto \deg(y)$ is not locally constant on~$\Gamma_{S}$ at~$x$.\footnote{This cannot happen when~$K$ is algebraically closed.}
\end{enumerate}
This set is locally finite (see~\cite[(4.5.12)]{Duc} for the last condition). 

%
%
%
%
\end{definition}

\begin{lemma}\label{lem:preimage}
Let $X$ be a $K$-analytic curve endowed with a pseudo-triangulation~$S$.
Let~$I$ be an open segment in~$X$.
\begin{enumerate}
\item Assume that $I \subseteq X-\Gamma_{S}$. The connected 
component of~$X-\Gamma_{S}$ containing~$I$ is a virtual open disk and we denote by~$b_{\infty}$ its germ of segment at infinity. If $I$ 
and~$b_{\infty}$ are aligned (\textit{i.e.} belong to a common segment) 
and if the map $x\mapsto \deg(x)$ is constant on~$I$, then $I$ is the 
skeleton of  an open pseudo-annulus.


\item Assume that $I \subseteq \Gamma_{S} - E_{S}$. If~$I$ is relatively 
compact, then $r_{S}^{-1}(I)$ (cf. \eqref{eq : retraction prop:GammaSS'}) is 
a virtual open annulus whose skeleton is~$I$. In the general case, 
$r_{S}^{-1}(I)$ is an open pseudo-annulus whose skeleton is~$I$.  

\end{enumerate}
\end{lemma}
\begin{proof}
ii) Let $I \subseteq \Gamma_{S} - E_{S}$. If $I$ is relatively compact, then 
the result follows from \cite[(5.1.12.3)]{Duc}. The general case follows by 
writing~$I$ as an increasing union of relatively compact open segments.

%
%

%
%

Let~$D$ be the connected component of~$X-\Gamma_{S}$ containing~$I$. It is a virtual open disk. 
Let $S'$ be a triangulation of~$D$. Let~$\Delta$ denote the convex hull of $\Gamma_{S'} \cup I$. By \cite[(1.5.19), (5.1.6.2)]{Duc}, $\Delta$ is admissible and analytically admissible. By \cite[(5.1.14)]{Duc}, there exists a triangulation~$S''$ such that $\Gamma_{S''} = \Delta$. 

Let $s \in I$ that is a bifurcation point of~$\Gamma_{S''}$. Let $C$ be a connected component of~$D - \{s\}$ that meets~$\Gamma_{S''}$ but not~$I$. Since~$b_{\infty}$ and~$I$ are aligned, $C$ does not contain~$b_{\infty}$, hence $C$ is a virtual open disk. As a consequence, $\Gamma_{S''} \cap C$ is relatively compact, with boundary~$s$. Using Lemma~\ref{lem:S-s} repeatedly, we may remove all the points of $C \cap S''$ from~$S''$ and stil get a triangulation of~$D$. Using the same argument for each connected component of~$D - \{s\}$, we end up in a situation where $s$ is no longer a bifurcation point of~$\Gamma_{S''}$. Finally, we may assume that $I$ contains no bifurcation point of~$\Gamma_{S''}$.

The assumptions now imply that $I \cap E_{S''} = \emptyset$. The result then follows from point~ii).
\end{proof}

\begin{proposition}\label{prop:pseudo-annulus}
Let~$C$ be an open pseudo-annulus. 

\begin{enumerate}
\item $C$ is a countable increasing union of virtual open annuli with nested skeletons.
\item $C$ has exactly two germs of segment at infinity~$b_{-}$ and~$b_{+}$ and we have 
\begin{equation}
N_{\infty}(C) = 2 \Card(\pi_{0}(C_{\wKa})).
\end{equation}
\item For every $x\in \Gamma_{C}$, we have
\begin{equation}
\deg(x) = \deg(b_{-}) = \deg(b_{+}) = \Card(\pi_{0}(C_{\wKa})).
\end{equation}
\end{enumerate}
\end{proposition}
\begin{proof}
Denote by $r_{C} \colon C \to \Gamma_{C}$ the canonical retraction. We can write the open segment~$\Gamma_{C}$ as a countable increasing union~$(I_{n})_{n\ge 0}$ of relatively compact open segments. By Lemma~\ref{lem:preimage}, for every~$n$, $r_{C}^{-1}(I_{n})$ is a virtual open annulus with skeleton~$I_{n}$. Point~i) follows, as well as the fact that $C$ has exactly two germs of segment at infinity.

The statements about the degrees follow from Proposition~\ref{prop:degannulus}.

%

\end{proof}


\subsection{Extension of scalars.}\label{section:extensionofscalars}

Let $X$ be a $K$-analytic curve. 

%
%

Let~$L$ be an algebraically closed complete valued 
extension of~$K$. Let~$M$ be a complete valued 
extension of~$L$. In this case, 
by~\cite[Corollaire~3.14]{Angie}, every point 
of~$X_{L}$ may be canonically lifted to a point 
in~$X_{M}$. We denote the canonical lift of a 
point~$x \in X_{L}$ by $\sigma_{M/L}(x)$ or simply 
by~$x_{M}$ when no confusion may arise. 

By~\cite[Lemma~1.15]{Push}, if~$x$ is a point of~$X_{L}$ of type~2, then~$x_{M}$ is a point of type~2 and we have $g(x_{M})=g(x)$. The structure of the fiber~$\pi^{-1}_{M/L}(x)$ may be 
described precisely: it is connected and the connected 
components of $\pi^{-1}_{M/L}(x) - \{x_{M}\}$ are 
virtual open disks with boundary~$\{x_{M}\}$ that are 
open in~$X_{M}$ (see~\cite[Theorem~2.2.9]{NP-II}).

\bigbreak

We will need to understand the behavior of pseudo-triangulations with 
respect to extensions of scalars. The arguments that we used 
in~\cite[Section~2.2]{NP-II} for weak triangulations still apply here, because of proposition~\ref{prop:GammaSS'}.


Let $S$ be a pseudo-triangulation of~$X$. 
It is easy to check that the preimage $S_{\wKa} := \pi^{-1}_{\wKa}(S)$ is a pseudo-triangulation of~$X_{\wKa}$ and that $\Gamma_{S_{\wKa}} = \pi^{-1}_{\wKa}(\Gamma_{S})$. To go further, we use the canonical lifts mentioned above.


\begin{definition}\label{def:liftS}
Let~$L$ be a complete valued extension of~$K$. Let~$\wLa$ be the completion of an algebraic closure of~$L$, which we see as an extension of~$\wKa$. Set
\begin{equation}
S_{\wLa} := \{x_{\wLa} \mid x\in S_{\wKa}\}\quad ,\quad \Gamma_{S_{\wLa}} :=\{x_{\wLa} \mid x\in \Gamma_{S_{\wKa}}\}
\end{equation}
and
\begin{equation}
S_{L} :=  \pi_{\wLa/L}(S_{\wLa})\quad,\quad \Gamma_{S_{L}} :=  \pi_{\wLa/L}(\Gamma_{S_{\wLa}}).
\end{equation}
\end{definition}

The sets of the previous definition are well-defined and independent of the choices. 

\begin{proposition}\label{prop:liftS}
For every complete valued extension~$L$ of~$K$, the set~$S_{L}$ is a pseudo-triangulation of~$X_{L}$ whose skeleton is~$\Gamma_{S_{L}}$. \hfill $\Box$
\end{proposition}

The proof is based on the description of the structure of the fibers of the map~$\pi_{\wLa/\wKa}$ (see \cite[Corollary~2.2.10]{NP-II}).

\medbreak

Using pseudo-triangulations and extensions of scalars, it is possible to associate disks to the points of~$X$. 

\begin{notation}\label{nota:D(x)}
Let $x\in X$. Let~$L$ be an algebraically closed complete valued extension of~$K$ such that $X_{L}$ contains an $L$-rational point~$t_{x}$ over~$x$. 

We denote by $D(x)$ the biggest open disk with center~$t_{x}$ in~$\pi_{L}^{-1}(x)$. 

We denote by $D(x,S)$ the biggest open disk with center~$t_{x}$ in~$X_{L}$ that does not meet~$S_{L}$.

For each $r \in \of{]}{0,1}{]}$, we denote by $D(x,S,r)$ the open subdisk of~$D(x,S)$ with center~$t_{x}$ whose radius is $r$ times that of~$D(x,S)$.
\end{notation}

Remark that the notation for the disks does not keep track of~$L$ and~$t_{x}$. This is due to the fact that, over a big enough field (algebraically closed and spherically complete), all the rational points over~$x$, hence all the disks centered at those points, are permuted by the action of the continuous Galois group (see \cite[Corollary~2.20]{NP-II}). As a result, the different possible disks behave similarly and the choice of~$L$ and~$t_{x}$ is rather irrelevant for our purposes 
(cf. \cite[Sections 2.1 and 2.2]{NP-III}).

\medbreak

We may also investigate the behavior of some of the invariants introduce before under scalar extension. 

\begin{lemma}
Let~$L$ be a complete valued extension of~$K$. Then, we have
\begin{equation}
N_{\infty}(X_{L}) = N_{\infty}(X)\;.
\end{equation}
\end{lemma}
\begin{proof}
This follows from Lemma~\ref{Lemma : Top finite >> N(X)<inf} and the behavior of pseudo-triangulations under scalar extension.
\end{proof}

\begin{lemma}\label{lem:extensionddc}
Let~$x\in X$ with finitely many preimages in~$X_{\wKa}$. Let~$L$ be a complete valued extension of~$\wKa$. Let $x_{1},\dotsc,x_{\deg(x)}$ be the preimages of~$x$ in~$X_{\wKa}$. For each $i\in\{1,\dotsc,\deg(x)\}$, let~$y_{i}$ be the canonical lift of~$x_{i}$ in~$X_{L}$.

Let $F \colon X \to \ERRE$ be a map that is $\log$-affine on every germ of segment out of~$x$ and constant on all but finitely many of them. Then, for each $i\in\{1,\dotsc,\deg(x)\}$, the map $F_{L} := F \circ \pi_{L} \colon X_{L} \to \ERRE$ is log-affine on every germ of segment out of~$y_{i}$ and constant on all but finitely many of them. Moreover, we have
\begin{equation}
dd^c F_{L}(y_{1}) = \dotsc = dd^c F_{L}(y_{\deg(x)}) = \frac{1}{\deg(x)}\, dd^c F(x)\;. 
\end{equation}
\hfill$\Box$
\end{lemma}

\subsection{Graphs in curves.}

Let $X$ be a $K$-analytic curve~$X$. Recall that it has a graph structure, by~\cite[Th\'eor\`eme~3.5.1]{Duc}. 

\begin{definition}
We say that a germ of segment~$b$ in~$X$ \emph{belongs to
a graph $\Gamma \subseteq X$} if~$b$ may be represented by a segment belonging to~$\Gamma$. We then write $b\subseteq \Gamma$.
\end{definition}

We wiil use the following notation. 

\begin{notation}
Let~$\Gamma$ be a subgraph of~$X$. 
For every~$x\in X$, we denote by 
\begin{equation}\label{eq : B(x,Gamma)}
\Bs(x, \Gamma) \quad\textrm{(resp. } \Bs(x, -\Gamma)\textrm{)}
\end{equation} 
the set of germs of segment out of~$x$ that belong to~$\Gamma$ (resp. do not belong to~$\Gamma$). 
\end{notation}

\begin{notation}
Let~$x\in X$. Let $F \colon X \to \ERRE$ be a map that is $\log$-affine on every germ of segment out of~$x$ and constant on all but finitely many of them. Let~$\Gamma$ be a subgraph of~$X$. We set 
\begin{equation}
dd^c_{\subseteq \Gamma} F(x) \; := \; \sum_{b \in \Bs(x, \Gamma)} \deg(b)\, \partial_{b} F
\end{equation}
and
\begin{equation}
dd^c_{\not\subseteq \Gamma} F(x) \; := \; \sum_{b \in \Bs(x, -\Gamma)} \deg(b)\, \partial_{b} F\;.
\end{equation}
\end{notation}


\begin{lemma}\label{lem:locallybounded}
Let~$x\in X$. For every segment~$[x,y]$ in~$X$, there exists $z\in\, ]x,y]$ such that the map $\deg \colon t\mapsto \deg(t)$ is constant on~$]x,z]$.

In particular, if~$\Gamma$ is a subgraph of~$X$ such that $\Bs(x, \Gamma)$ is finite, then the map~$\deg$ is locally bounded at~$x$ on~$\Gamma$.
\end{lemma}
\begin{proof}
If~$x$ has type~1 or~4, then it has a neighborhood in~$X$ that is isomorphic to a virtual open disk~$D$. There exists $y'\in\, ]x,y]$ such that $[x,y'] \subset D$. By~\cite[1.9.16.3]{Duc}, the map~$\deg$ is non-decreasing on~$[x,y']$ towards~$x$, and the result follows. 

If~$x$ has type~2 or~3, then, by~\cite[4.5.12]{Duc}, every germ of segment~$b$ out of~$x$ has a representative~$I_{b}$ on which the map~$\deg$ is constant. The result follows.
\end{proof}

\begin{corollary}\label{cor:GammawKa-1}
Let~$\Gamma$ be a quasi-finite subgraph of~$X$. Assume that~$X$ may be embedded as an analytic domain in a quasi-smooth $K$-analytic curve~$X'$ such that~$\Gamma$ is relatively compact in~$X'$ (\textit{e.g.} if~$X$ may be embedded in the analytification of an algebraic curve). Then, $\pi^{-1}_{\wKa}(\Gamma)$ is quasi-finite.
\end{corollary}
\begin{proof}
We may assume that~$X'=X$. It is enough to prove that the map $\deg \colon x\mapsto \deg(x)$ is bounded on~$\Gamma$. The closure~$\overline \Gamma$ of~$\Gamma$ in~$X$ being a compact graph, it is enough to prove that it is locally bounded on~$\overline \Gamma$.

Since~$\Gamma$ is quasi-finite, its closure~$\overline\Gamma$ in the graph~$X$ is quasi-finite too, hence, for every $x\in \overline\Gamma$, there are only finitely many segments out of~$x$ belonging to~$\overline\Gamma$. The result now follows from Lemma~\ref{lem:locallybounded}.
\end{proof}

It is sometimes useful to be able to consider vertices and edges of graphs. We give a corresponding definition. Beware that it differs from the usual one since we allow edges with no end-points. We will use this definition in Section~\ref{section : explicit bound}.

\begin{definition}\label{def:markedgraph}
A \emph{marked graph} is a pair $(\Gamma,V)$ where~$\Gamma$ is a locally finite graph and~$V$ a locally finite subset of~$\Gamma$ such that every connected component of~$\Gamma-V$ is an open segment.

An element of~$V$ is called the a \emph{vertex} of the graph, a connected component of~$\Gamma-V$ is called an \emph{edge} of the graph.

By abuse of notation, we usually say that~$\Gamma$ is a marked graph. In this case, we denote by~$V(\Gamma)$ and~$E(\Gamma)$ its sets of vertices and edges. 
\end{definition}

\begin{lemma}
Let~$\Gamma$ be a locally finite graph.
\begin{enumerate}
\item The graph~$\Gamma$ is quasi-finite if, and only if, it may be endowed with the structure of a marked graph with finitely many vertices and edges.
\item The graph~$\Gamma$ is finite if, and only if, it may be endowed with the structure of a marked graph with finitely many vertices and edges such that all the edges are relatively compact.
\end{enumerate}
\hfill$\Box$
\end{lemma}

\begin{definition}\label{def:markedanalyticgraph}
A \emph{marked analytic subgraph} of~$X$ is a marked graph $(\Gamma,V)$ such that
\begin{enumerate}
\item $\Gamma$ is a subgraph of~$X$;
\item every edge of~$\Gamma$ is the skeleton of an open pseudo-annulus.
\end{enumerate}
\end{definition}


Remark that, if~$S$ is a pseudo-triangulation of the curve~$X$, then $(\Gamma_{S},S)$ is a marked analytic subgraph of~$X$.

\begin{lemma}\label{lem:markedanalyticgraph}
Let~$\Gamma$ be a locally finite subgraph of~$X$ that contains~$\Gamma_{S}$. Let~$V$ be the union of
\begin{enumerate}
\item the elements of~$S$;
\item the bifurcation points of~$\Gamma$;
\item the end-points of~$\Gamma$;
\item the points of~$\Gamma$ at which the map $x\mapsto \deg(x)$ is not locally constant.
\end{enumerate}
Then $(\Gamma,V)$ is a marked analytic subgraph of~$X$.

The set~$V$ is the smallest set containing~$S$ with this property.
\hfill$\Box$
\end{lemma}
\begin{proof}
It is clear that the connected components of~$\Gamma-V$ are open segments. By Lemma~\ref{lem:preimage}, they are skeletons of open pseudo-annuli.
\end{proof}

In the same way that we lifted triangulations (see Definition~\ref{def:liftS} and Proposition~\ref{prop:liftS}), we can lift marked analytic graphs to any complete valued extension of~$K$, starting with the case of~$\wKa$. The proofs are the same and we omit them. 

Let us first lift sets of points.

\begin{definition}\label{def:liftV}
Let~$V$ be a subset of~$X$. Let~$L$ be a complete valued extension of~$K$. Let~$\wLa$ be the completion of an algebraic closure of~$L$, which we see as an extension of~$\wKa$. Set
\begin{equation}
V_{L} := \pi_{\wLa/L}(\{x_{\wLa} \mid x\in \pi^{-1}_{\wKa}(V)\}).
\end{equation}
\end{definition}


This set is well-defined.

\begin{lemma}
Let $(\Gamma,V)$ be a marked analytic subgraph of~$X$. Then $(\Gamma_{\wKa} := \pi^{-1}_{\wKa}(\Gamma), V_{\wKa} := \pi^{-1}_{\wKa}(V))$ is a marked analytic subgraph of~$X_{\wKa}$ and we have
\begin{equation}
E(\Gamma_{\wKa}) = \{\textrm{connected components of } \pi^{-1}_{\wKa}(E),\ E \in E(\Gamma)\}.
\end{equation}
\hfill$\Box$
\end{lemma}


\begin{definition}\label{def:liftGamma}
Let $(\Gamma,V)$ be a marked analytic subgraph of~$X$. Let~$L$ be a complete valued extension of~$K$. Let~$\wLa$ be the completion of an algebraic closure of~$L$, which we see as an extension of~$\wKa$. Set
\begin{equation}
E_{\wLa} :=  \{\textrm{connected components of } \pi^{-1}_{\wLa}(E),\ E \in E(\Gamma_{\wKa})\}\quad,
\end{equation}
and
\begin{equation}
E_{L} :=  \{ \pi_{\wLa/L}(E),\ E \in E(\Gamma_{\wLa}) \}.
\end{equation}
\end{definition}

The sets of the previous definition are well-defined and independent of the choices. 

\begin{proposition}\label{prop:liftGamma}
Let $(\Gamma,V)$ be a marked analytic subgraph of~$X$. For every complete valued extension~$L$ of~$K$, the pair $(\Gamma_{L},V_{L})$  is a marked analytic subgraph of~$X_{L}$ and we have
\begin{equation}
E(\Gamma_{L}) = E_{L}.
\end{equation}
\hfill$\Box$
\end{proposition}

\subsection{Euler characteristic.}

We now define the Euler characteristic and the genus of a curve. Let~$X$ be a quasi-smooth $K$-analytic curve. 

\begin{definition}\label{Def : Finite genus}
Assume that~$K$ is algebraically closed and that~$X$ is 
connected. 
We denote by $\chi_{\topo}(X) \in \of{]}{-\infty,1}{]}$ 
the Euler characteristic (in the sense of singular 
homology) of the topological space underlying~$X$.

We define the genus of~$X$ by
\begin{equation}\label{eq : def of g(X)}
g(X)\;:=\;  1-\chi_{\topo}(X) +\sum_{x\in X} g(x) \;\in\; \of{[}{0,+\infty}{[} \;.
\end{equation}

We define the compactly supported Euler characteristic of~$X$ by
\begin{equation}\label{eq : def of chic(X)}
\chi_{c}(X)\;:=\;  2 - 2 g(X) - N_{\infty}(X) \;\in\; \of{]}{-\infty,2}{]} \;.
\end{equation}


If~$K$ is arbitrary and~$X$ has a finite number of connected components, we define the genus and compactly supported Euler characteristic of~$X$ to be the sum of those of the connected components of~$X_{\wKa}$. \end{definition}

For every open pseudo-disk~$D$, we have $g(D) = 0$ and $\chi_{c}(D) = \Card(\pi_{0}(D_{\wKa}))$ and for every open pseudo-annulus~$C$, we have $g(C)=\chi_{c}(C) = 0$.


\begin{remark}
Assume that $K$ is algebraically closed and that~$X$ is 
connected. The topological Euler characteristic 
$\chi_{\topo}(X)$ of~$X$ may be easily computed in 
this case. 

Recall that~$\chi_{\topo}$ is a homotopy invariant. In particular, if~$X$ is contractible (\textit{e.g.} any analytic 
domain of~$\mathbb{P}^{1,\an}_{K}$), then we have $\chi_{\topo}(X) = 1$.

If~$\Gamma_{S}\ne \emptyset$, then the curve~$X$ retracts by deformation onto~$\Gamma_{S}$, hence $\chi_{\topo}(X) = \chi_{\topo}(\Gamma_{S})$. If the graph~$\Gamma_{S}$ itself retracts by 
deformation onto a finite subgraph~$\Gamma'$ (in the sense of Definition~\ref{def:finitegraph}), then its Euler 
characteristic may be explicitly computed. To this end, endow~$\Gamma'$ with a structure of strict finite graph with~$v'$ 
vertices and~$e'$ edges. Then, we have
\begin{equation}
\chi_{\topo}(X) = \chi_{\topo}(\Gamma_{S}) = \chi_{\topo}(\Gamma') = v'-e'\;.
\end{equation}
\end{remark}

\begin{remark}\label{rem:genus}
Assume that~$K$ is algebraically closed and that~$X$ is connected.
\begin{enumerate}
\item If~$X$ is the analytification of an algebraic curve~$C$, then $g(X)=g(C)$ (see~\cite[p.82, before Thm. 4.3.1]{Ber} 
or~\cite[5.2.6]{Duc}).
\item The quantify~$\chi_{c}(X)$ could also be defined as a compactly supported \'etale Euler characteristic (see~\cite[Section~5]{bleu} for the definition of \'etale cohomology with compact support).
%
\end{enumerate}
\end{remark}



Curves of finite genus are close to algebraic ones thanks to the following result.

\begin{theorem}[\protect{\cite[Th\'eor\`eme 3.2]{Liu}}]
\label{Thm : Liu-plongement}
Assume that $K$ is non-trivially valued, algebraically closed and spherically complete. Assume that~$X$ is strictly $K$-analytic, connected and has finite genus~$g$. Then, there exists a smooth connected projective curve~$Y$ with the same genus~$g$ such that~$X$ is isomorphic to an analytic domain of $Y^\an$.\hfill$\Box$
\end{theorem}

\begin{lemma}
Assume that~$X$ has finitely many connected components. Then, for every complete valued extension~$L$ of~$K$, we have 
\begin{equation}
g(X_{L}) \;=\;  g(X) \textrm{ and }
\chi_{c}(X_{L}) \;=\;  \chi_{c}(X)\;.
\end{equation}
\hfill $\Box$
\end{lemma}

\begin{notation}\label{Def; chi(x,S)} 
Let~$\Gamma$ be a subgraph of~$X$. Let~$x$ be a point of~$\Gamma$. We set 
\begin{equation}
N_{\Gamma}(x) := \sum_{b\in\Bs(x,\Gamma)} \deg(b)
\end{equation}
and 
\begin{equation}
\chi(x,\Gamma) := 2\deg(x) - 2 g(x) - N_{\Gamma}(x).
\end{equation}

When~$\Gamma=\Gamma_{S}$, we set $N_{S}(x) := N_{\Gamma_{S}}(x)$ and $\chi(x,S) := \chi(x,\Gamma_{S})$ for short. 
\end{notation}

\begin{lemma}
Let~$\Gamma$ be a subgraph of~$X$. Let~$x$ be a point of~$\Gamma$. Let~$L$ be a complete valued extension of~$\wKa$. Let $x_{1},\dotsc,x_{\deg(x)}$ be the preimages of~$x$ in~$X_{\wKa}$. For each $i\in\{1,\dotsc,\deg(x)\}$, let~$y_{i}$ be the canonical lift of~$x_{i}$ in~$X_{L}$.
Then, we have
\begin{equation}
N_{\Gamma_{L}}(y_{1}) = \dotsb = N_{\Gamma_{L}}(y_{\deg(x)}) = \frac{1}{\deg(x)}\, N_{\Gamma}(x)
\end{equation}
and
\begin{equation}
\chi(y_{1},\Gamma_{L}) = \dotsb = \chi(y_{\deg(x)},\Gamma_{L})  = \frac{1}{\deg(x)}\, \chi(x,\Gamma).
\end{equation}
\hfill$\Box$
\end{lemma}

\begin{lemma}
\label{Lemma : chi(X)=sum chi(x,S)}

Assume that~$X$ has finitely many connected components and that~$S$ meets all of them. Moreover, assume that
\begin{enumerate}
\item the skeleton $\Gamma_S$ is quasi-finite and the map $x\mapsto \deg(x)$ is bounded on it;
\item the number of points with positive genus is finite;
\item the boundary~$\partial X$ is finite.
\end{enumerate}
Then, there exists a finite subset~$F$ of~$S$ such that, for every finite subset~$S'$ of~$S$ containing~$F$, we have
\begin{equation}\label{eq : hfsdtleoi}
\sum_{x\in S'}\chi(x,S') \;=\;  \chi_{c}(X).
\end{equation}
\end{lemma}
\begin{proof}
We may assume that~$K$ is algebraically closed. Indeed, assumption~i) ensures that~$\Gamma_{S_{\wKa}}$ is quasi-finite, and similarly for the others.

From the assumptions we deduce that there exists a marked finite subgraph~$\Gamma$ of~$\Gamma_{S}$ with vertices in~$S$ and intervals $I_{1},\dotsc,I_{N(X)}$ such that
\begin{enumerate}
\item every point with positive genus belongs to~$\Gamma$~;
\item every point of~$\partial X$  belongs to~$\Gamma$~;
\item for every~$i\in\{1,\dotsc,N(X)\}$, $I_{i}$ is a semi-open segment that is closed in~$X$ and whose end-point belongs to~$S$~;
\item $\Gamma_{S} = \Gamma \cup \bigcup_{1\le i\le N(X)} I_{i}$.
\end{enumerate}

It follows from these conditions that $\Gamma_{S}$ (hence~$X$) retracts onto~$\Gamma$.

Set $F= S \cap \Gamma$. Remark that every point of $S-F$ belongs to one of the~$I_{i}$'s.

\medbreak

Let~$S'$ be a finite subset of~$S$ containing~$F$. 

Let~$i\in\{1,\dotsc,N(X)\}$. Let~$a_{i} \in S$ be the end-point of~$I_{i}$. For every points $x,y \in I_{i}$, we denote by $[x,y]$ the image of an injective path in~$I_{i}$ joining~$x$ to~$y$. This is well-defined. We define an order relation~$\le$ on~$I_{i}$ by setting $x\le y$ if $x \in [a_{i},y]$. Let~$b_{i}$ be the biggest element of~$I_{i}\cap S'$.

Let us now consider the graph $\Gamma' = \Gamma \cup \bigcup_{1\le i\le N(X)} [a_{i}, b_{i}]$ with set of vertices~$S'$. Denote by~$E'$ the number of edges of~$\Gamma'$. 

For every $x \in S' -\{b_{1},\dotsc,b_{N(X)}\}$, we have $N_{S'}(x) = N_{S}(x)$ and, for every $i\in \{1,\dotsc,N(X)\}$, we have $N_{S}(b_{i}) = 2 = N_{S'}(b_{i})+1$. We deduce that
\begin{align}
\sum_{x\in S'} (2-2g(x)-N_S(x)) & = \sum_{x\in S'} (2-2g(x)-N_{S'}(x)) - N(X)\\
&= -2 \sum_{x\in S'} g(x) + 2 \Card(S') - 2 E' - N(X)\;.
\end{align}
The curve~$X$ retracts onto~$\Gamma'$, hence
\begin{equation}
\chi_{\topo}(X) = \chi_{\topo}(\Gamma') = \Card(S') - E'
\end{equation}
and we now conclude by~\eqref{eq : def of g(X)}.
\end{proof}

\begin{lemma}\label{lem:Ninfty}
Let~$U$ and~$V$ be open subsets of~$X$ that have finite open boundaries. Then $U\cup V$ and $U\cap V$ have finite open boundaries and we have
\begin{equation}
N_\infty(U\cup V) = N_\infty(U) + N_\infty(V) - N_\infty(U\cap V)\;.
\end{equation}
\end{lemma}
\begin{proof}
We may assume that~$K$ is algebraically closed. For every $W = U, V, U\cup V$ or $U\cap V$, we let $\Gs_{nrc}(W)$ (resp. $\Gs_{rc}(W)$, resp. $\Gs_{\notin}(W)$) be the set of germs of segment of~$X$ that belong to~$W$ and are not relatively compact in~$W$ (resp. belong to~$W$ and are relatively compact in~$W$, resp. do not belong to~$W$).

Remark that an open subset contains all the germs of segment out of its points. We deduce that 
\begin{enumerate}
\item $\Gs_{nrc}(U\cup V) \cap \Gs_{nrc}(U\cap V) = \Gs_{nrc}(U) \cap \Gs_{nrc}(V)$\;;
\item $\Gs_{nrc}(U\cup V) \cap \Gs_{\notin}(U\cap V) = \big(\Gs_{nrc}(U) \cap \Gs_{\notin}(V)\big) \sqcup \big(\Gs_{\notin}(U) \cap \Gs_{nrc}(V)\big)$\;;
\item $\Gs_{rc}(U\cup V) \cap \Gs_{nrc}(U\cap V) = \big(\Gs_{nrc}(U) \cap \Gs_{rc}(V)\big) \sqcup \big(\Gs_{rc}(U) \cap \Gs_{nrc}(V)\big)$\;.
\end{enumerate}
Denote by~$N_{1}$, $N_{2}$ and~$N_{3}$ the cardinals of the sets in the lines i), ii) and iii) respectively. We have 
\begin{equation}
N_\infty(U\cup V) + N_\infty(U\cap V)= 2N_{1}+N_{2} + N_{3} = N_\infty(U)+N_\infty(V)\;.
\end{equation}
\end{proof}

\begin{corollary}
\label{Cor : chi(X)=chi(U)+chi(V)-chi(U cap V)}
Let~$U$ and~$V$ be open subsets of~$X$ that both have finite genus and finite open boundary. Then $U\cup V$ and $U\cap V$ have finite genus and finite open boundary and we have
\begin{equation}
\chi_{c}(U\cup V) = \chi_{c}(U) + \chi_{c}(V) - \chi_{c}(U\cap V)\;.
\end{equation}
\end{corollary}
\begin{proof}
We may assume that~$K$ is algebraically closed. By Mayer-Vietoris' formula for singular homology, we have 
\begin{equation}
\chi_{\topo}(U\cup V) = \chi_{\topo}(U) + \chi_{\topo}(V) - \chi_{\topo}(U\cap V)\;.
\end{equation}
Moreover, we have 
\begin{equation}
\sum_{x\in U\cup V} g(x) = \sum_{x\in U} g(x) + \sum_{x\in V} g(x) - \sum_{x\in U\cap V} g(x)\;.
\end{equation}
Using the Lemma~\ref{lem:Ninfty}, the result follows.
\end{proof}

\subsection{Finite curves.}
\label{Section : finite curves}

Let $X$ be a $K$-analytic curve.

\begin{definition}[Finite curve]\label{Def : finite curve}
We say that $X$ is a \emph{finite curve} 
if it has finitely many connected components and admits a finite pseudo-triangulation. 
\if{such 
that the pull-back $\Gamma_{S_{\wKa}}
\subseteq X_{\wKa}$ of its skeleton 
$\Gamma_S\subseteq X$ is topologically finite.}\fi
%
\end{definition}


We have the following characterization of finite curves.

%


\begin{lemma}\label{Lemma : finite ssi Gama_S top finite}
The curve~$X$ is finite if, and only if, 
\begin{enumerate}
\item $X$ has finitely many connected components;
\item $X$ has finite genus;
\item $X$ has finite boundary $\partial X$;
\item there exists a pseudo-triangulation~$S$ of~$X$ such that
\begin{enumerate}
\item $\Gamma_{S}$ is a quasi-finite graph;
\item the map $x \in\Gamma_{S} \mapsto \deg(x)$ is continuous outside finitely many points.
\end{enumerate}
\end{enumerate}
\hfill$\Box$
\end{lemma}

\begin{remark}\label{rem:degfini}
Condition iii) (b) of the previous lemma is satisfied in the following cases:
\begin{enumerate}
\item $K$ is algebraically closed;
\item $\Gamma_{S}$ is a finite graph, see Lemma~\ref{lem:locallybounded};
\item $\Gamma_{S}$ is quasi-finite and $X$ may be embedded as an analytic domain in a quasi-smooth $K$-analytic curve~$X'$ such that~$\Gamma_{S}$ is relatively compact in~$X'$ (\textit{e.g.} if~$X$ may be embedded in the analytification of an algebraic curve), see Corollary~\ref{cor:GammawKa-1}.
\end{enumerate}
\end{remark}

\begin{lemma}
Let $Y$ be a quasi-smooth $K$-analytic curve.
Let $U,V\subseteq Y$ be two analytic domain in $Y$.
If $U$ and $V$ are finite curves, then so are $U\cup V$ 
and $U\cap V$.
\end{lemma}
%
%
%
%
%
%
%


\begin{proposition}\label{Prop : smallest}
Let~$Y$ be a quasi-smooth $K$-analytic curve. If no connected component of~$Y$ is isomorphic to~$\mathbb{P}^{1,\an}_{K}$ or a Tate curve, then~$Y$ admits a smallest 
pseudo-triangulation~$S_{Y}$. 

Moreover, the skeleton of~$S_{Y}$ coincides with the analytic skeleton~$\Gamma_{Y}$ of~$Y$ and, for each $y\in S_{Y}$, we have $\chi(y,S_{Y})\le 0$.
\end{proposition}

\begin{proposition}\label{Prop : smallest-1}
Let~$Y$ be a quasi-smooth $K$-analytic curve. Then
$Y$ admits a smallest pseudo-triangulation~$S_{Y}$. 

Moreover, if no connected component of~$Y$ is 
isomorphic to~$\mathbb{P}^{1,\an}_{K}$ or a Tate 
curve, then the following hold:
\begin{enumerate}
\item $S_Y$ is unique;
\item the skeleton of~$S_{Y}$ coincides with the 
analytic skeleton~$\Gamma_{Y}$ of~$Y$;
\item for each $y\in S_{Y}$, 
we have $\chi(y,S_{Y})\le 0$.
\end{enumerate}
\end{proposition}
%

\section{Super-harmonicity of partials heights}\label{sec:sh}

Let $X$ be a quasi-smooth $K$-analytic curve endowed with a pseudo-triangulation~$S$. Let $(\Fs,\nabla)$ be a module with connection of rank~$r$ on~$X$. For simplicity, we will often remove~$\nabla$ from the notation.

In this section, we are interested in super-harmonicity properties of the partial heights $H_{S,i}(-,\Fs)$. We provide conditions ensuring that $dd^c H_{S,i}(x,\Fs)\le 0$ at a given point~$x$ and, in general, compute upper bounds for $dd^c H_{S,i}(x,\Fs)$.


%

\subsection{Partial heights and their behavior.}
\label{A criterion for the controlling graphs-ref}

In this section, we recall the definitions of radii of convergence and partial heights and state some properties. 

Let $x\in X$. Recall Notation~\ref{nota:D(x)} for the disks attached to~$x$. For $i\in \{1,\dotsc,r\}$, we denote by 
\begin{equation}
\R_{S,i}(x,\Fs)
\end{equation}
the \emph{$i^\textrm{th}$ radius of convergence} of~$\Fs$ at~$x$, that is to say the supremum of the set of $R \in \of{]}{0,1}{]}$ such that $\Fs$ admits $r-i+1$ linearly independent solutions on~$D(x,S,R)$ (see \cite[Definition~2.36]{NP-II}).

For $i\in \{1,\dotsc,r\}$, we define the \emph{$i^\textrm{th}$ partial height} of~$\Fs$ at~$x$ by
\begin{equation}
H_{S,i}(x,\Fs) := \prod_{j=1}^i \R_{S,i}(x,\Fs).
\end{equation}

We say that the radius $\R_{S,i}(x,\Fs)$ (or the index~$i$) is \emph{spectral} (resp. \emph{solvable}, resp. \emph{over-solvable}) if we have $D(x,S,\R_{S,i}(x,\Fs))\subseteq D(x)$ (resp. $=D(x)$, resp. $\supset D(x)$).


We denote by $i_x^{\mathrm{sp}}, \is{x} \in \{0,\dotsc,r\}$ the indexes such that $\R_{S,i}(x,\Fs)$ is spectral non-solvable for $i\leq i_x^{\mathrm{sp}}$, solvable for 
$i_x^{\mathrm{sp}}<i\le \is{x}$ and over-solvable for $i>\is{x}$.

We say that $i$ is a \emph{vertex} at $x$ if $i=r$, or if $i <r$ and $\R_{S,i}(x,\Fs) < \R_{S,i+1}(x,\Fs)$. We say that $i$ is a \emph{vertex free of solvability} at $x$ if, moreover, $i_x^{\mathrm{sp}} = \is{x}$.

%
%

By \cite[Theorem~3.6]{NP-II}, for each $i\in \{1,\dotsc,r\}$, the function $\R_{S,i}(-,\Fs)$ is continuous on~$X$ and factors through the retraction on a locally finite subgraph of~$X$. In particular, the controlling graph (see Definition~\ref{def:controllinggraph})
\begin{equation} 
\Gamma_{S,i}(\Fs) := \Gamma_{S}(\R_{S,i}(-,\Fs))
\end{equation}
is locally finite.

Set $\Gamma_{S,0}^\tot(\Fs):=\Gamma_S$ and, for $i \in \{1,\dotsc,r\}$, 
\begin{equation}\label{eq : Gamma'_S,i(Fs)}
\Gamma_{S,i}^\tot(\Fs)\;:=\;
\Gamma_{S,1}(\Fs)\cup\cdots\cup\Gamma_{S,i}(\Fs)\;.
\end{equation}

%
%
%

%

\medbreak

For later use, we record here a result about the form of the slopes of the radii of convergence (see \cite[Theorem~3.6 (iii))]{NP-II}).


\begin{proposition}\label{prop:slopesm/i}
Let $i \in\{1,\dotsc,r\}$. On any segment in~$X$, the slopes of~$\Rc_{S,i}(-,\Fs)$ are of the form $m/j$ with $m\in\NN$ and $1\le j\le r$.
\qed
\end{proposition}

\begin{notation}\label{nota:kappar}
Set
\begin{equation}
\kappa_{r} := \max(r(r-1),1) = 
\begin{cases}
1 & \textrm{if } r=1;\\
r(r-1) & \textrm{if } r\ge 2.
\end{cases}
\end{equation}
\end{notation}

\begin{corollary}\label{cor:penteminoree}
Let $x$ and $y$ be points of~$X$ and let $b$ and $c$ be germs of segment out of $x$ and $y$ respectively. Let $i,j \in \{1,\dotsc,r\}$. If $\partial_{b} H_{S,i}(x,\Fs) \ne 0$, then we have
\begin{equation}
|\partial_{b} H_{S,i}(x,\Fs)| \ge \frac1r
\end{equation}
and, if $\partial_{c} H_{S,j}(y,\Fs) \ne\partial_{b} H_{S,i}(x,\Fs)$, then we have
\begin{equation}
|\partial_{c} H_{S,j}(y,\Fs) -\partial_{b} H_{S,i}(x,\Fs)| \ge \frac{1}{\kappa_{r}}.
\end{equation}
%
\qed
\end{corollary}
\begin{proof}
This is an easy consequence of the fact that $\partial_{b} H_{S,i}(x,\Fs)$ and $\partial_{c} H_{S,j}(y,\Fs)$ belong to the set $\bigcup_{1\le s\le r} \frac{1}{s}\,\Z$ (see \cite[Theorem~3.6 (iii))]{NP-II}). 
\end{proof}

\begin{lemma}[\protect{\cite[Thm.4.7]{NP-I}}]
\label{Lemma: concave alonga spectral annulus}
Let $i\leq r$.
Let $I\subseteq\Gamma_S$ be the skeleton of an annulus in $X$. 
Then the index $i$ is spectral on $I$, 
and $H_{S,i}(-,\Fs)$ is $\log$-concave on $I$.\hfill$\Box$
\end{lemma}

\begin{lemma}[\protect{\cite[11.3.2]{Kedlaya-Book}}] 
Let $i\leq r$.
Let $D\subseteq X$ be a disk such that $D\cap S=\emptyset$. 
Let $]x,y[$ be an open segment in $D$ oriented towards the exterior of $D$. 
If the index $i$ is spectral non solvable at each point of $]x,y[$, 
then $H_{S,i}(-,\Fs)$ is decreasing on it.
\hfill$\Box$
\end{lemma}
  
\begin{remark}[From \cite{NP-I}]
\label{remark : Linear properties of Gamma'}
The graphs $\Gamma_{S,i}'(\Fs)$ satisfy the following properties.
\begin{enumerate}
\item For every $i=0,\ldots,r$, the topological space 
$X-\Gamma_{S,i}'(\Fs)$ is a 
disjoint union of virtual open disks of the form 
$D(y,\Gamma_{S,i}'(\Fs))$, where 
$y\notin\Gamma_{S,i}'(\Fs)$ is a rigid point.
\item The radii $\R_{S,1}(-,\Fs),\ldots,\R_{S,i-1}(-,\Fs)$ 
are constant functions on the disk 
$D(x,\Gamma_{S,i-1}'(\Fs))$, for all $x\in X$. 
In particular the ratio $\R_{S,i}(-,\Fs)/H_{S,i}(-,\Fs)$ is constant on  
$D(x,\Gamma_{S,i-1}'(\Fs))$. This implies that 
the controlling graphs and the $\log$-slopes of 
 $\R_{S,i}(-,\Fs)$ and $H_{S,i}(-,\Fs)$
coincide on $D(x,\Gamma_{S,i-1}'(\Fs))$ :
\begin{equation}
\Gamma_{S,i}(\Fs)\cap(X-\Gamma_{S,i-1}'(\Fs))\;=\;
\Gamma_{S}(H_{S,i}(-,\Fs))\cap(X-\Gamma_{S,i-1}'(\Fs))\;.
\end{equation}
Hence 
\begin{equation}
\Gamma_{S,i}'(\Fs)\;=\;\bigcup_{j=1}^i
\Gamma_{S}(H_{S,j}(-,\Fs))\;.
\end{equation}
\item Either  
$\R_{S,i-1}(-,\Fs)=\R_{S,i}(-,\Fs)$ as functions over 
$D(x,\Gamma_{S,i-1}'(\Fs))$, or they are separated at each point of it 
by \cite[Prop.  7.5]{NP-I}. 
This follows from the fact that the restriction 
of $\Fs$ to $D:=D(x,\Gamma_{S,i-1}'(\Fs))$ 
decomposes by \cite[12.4.1]{Kedlaya-Book} as a direct sum
$\Fs_{|D}=(\Fs_{|D})_{\geq i}\oplus (\Fs_{|D})_{<i}$. 
Hence $\R_{S,i}(-,\Fs)$ and $H_{S,i}(-,\Fs)$ behave as  
first radii of convergence outside $\Gamma_{S,i-1}(\Fs)$. 
So they have the concavity property of point iv) below. 
\end{enumerate}
Point iii) is used in \cite[Section 7.4]{NP-I} to prove that the radii are 
separated over $D(x,\Gamma_{S,i-1}'(\Fs))$. %
As a consequence they have the following property:
\begin{itemize}
\item[iv)] 
Let $]z,y[$ be a segment in $X$ such that 
$\Gamma_{S,i-1}'(\Fs)\cap]z,y[=\emptyset$. 
We consider $]z,y[$ as oriented towards the exterior of 
the disk $D(x,\Gamma_{S,i-1}'(\Fs))$ containing it.
Then the functions $\R_{S,i}(-,\Fs)$ and $H_{S,i}(-,\Fs)$ 
are log-concave and decreasing on $]z,y[$
(cf. \cite[sections 3.1,3.2]{NP-I}). 
We refer to this property by saying that $\R_{S,i}(-,\Fs)$ and 
$H_{S,i}(-,\Fs)$ have \emph{the concavity property} outside 
$\Gamma_{S,i-1}'(\Fs)$.
\end{itemize}
In particular let $D\subset X$ be a virtual disk with boundary $x\in X$, 
such that $D\cap\Gamma_{S,i-1}'(\Fs)=\emptyset$, and let
$b$ be the germ of segment out of $x$ inside $D$. Then 
:
\begin{equation}\label{eq : partial R_b =0 .jky fd}
\textrm{$b\in\Gamma_{S,i}(\Fs)$\; if and only if } \;\;
\partial_b\R_{S,i}(x,\Fs)\;\neq \;0\;.
\end{equation}
Namely by iv) one has
\begin{equation}\label{eq : R_i constant on D iff slope not O}
\textrm{$\partial_b\R_{S,i}(x,\Fs)=0$\; if and only if\; $\R_{S,i}(-,\Fs)$
is constant on $D$\;.}
\end{equation}
The same is true replacing $\R_{S,i}$ by $H_{S,i}$ in 
\eqref{eq : partial R_b =0 .jky fd} and 
\eqref{eq : R_i constant on D iff slope not O}.
\end{remark}
The following lemma studies the structure of the controlling graphs in 
the solvable case. It is somehow a key lemma in what follows.

\begin{lemma}[\protect{\cite[Lemma 7.7]{NP-I}}]
\label{Lemma : solvable radii at boundary}
Assume that the index $i\in\{1,\ldots,r\}$ is solvable at 
$x\in X$, and that $x\in\Gamma_{S,i}(\Fs)$. 
Then the following holds:
\begin{enumerate}
\item If $x\in\Gamma_{S,i}(\Fs)-\Gamma_{S,i-1}'(\Fs)$, 
then $x$ is an end point of $\Gamma_{S,i}(\Fs)$; 
\item If $x\in\Gamma_{S,i-1}'(\Fs)$, then 
$\Gamma_{S,i}(\Fs)\subseteq \Gamma_{S,i-1}'(\Fs)$ 
around $x$ (i.e. if 
$[x,y[\subseteq \Gamma_{S,i}(\Fs)$, then 
$[x,y[\subseteq \Gamma_{S,i-1}'(\Fs)$ if $y$ is close enough to $x$). 
\end{enumerate}
\end{lemma}
\begin{proof}
Let $D$ be a virtual disk in $X-\Gamma_{S,i-1}'(\Fs)$ 
with boundary $x\in\Gamma_{S,i}(\Fs)$. 
With the notations of \eqref{eq : partial R_b =0 .jky fd},
solvability at $x$ implies $\partial_b\R_{S,i}(x,\Fs)\leq 0$, 
while the concavity property iv) of Remark 
\ref{remark : Linear properties of Gamma'} 
implies $\partial_b\R_{S,i}(x,\Fs)\geq 0$. 
So $\partial_b\R_{S,i}(x,\Fs)=0$ and 
$D\cap\Gamma_{S,i}(\Fs)=\emptyset$. 
\end{proof}

Reasoning as in \cite[Remark 7.1]{NP-I} one proves the following 
characterization of $\Gamma_{S,i}'(\Fs)$. 
It does not follow directly from from the Definition.
\begin{proposition}[\protect{\cite[Thm.4.7 iii), iv)]{NP-I}}]
\label{Prop. Characterization of the skeleton by the non constancy of the partial heights}
Let $\Gamma_i$ be the union of the closed 
segments $[x,y]\subseteq X$ on which at least one of the 
partial heights
$H_{S,1}(-,\Fs),\ldots,H_{S,i}(-,\Fs)$, or equivalently one of 
the radii $\R_{S,1}(-,\Fs),\ldots,\R_{S,i}(-,\Fs)$, is never constant on 
$[x,y]$. Then
\begin{equation}\label{eq : gamma'=Galla cum Gamma_i}
\Gamma_{S,i}'(\Fs)\;=\;\Gamma_S\cup\Gamma_i\;.
\end{equation}
\end{proposition}
\begin{proof}
\if{
We can assume $K$ is arbitrary large since the radii are 
insensitive to scalar extensions of $K$.  It is enough to prove 
that a closed segment $[x,y]$ in a maximal open disk 
$D(x,S)$ belongs to $\Gamma_{S,i}'(\Fs)$ if and only if at least one of 
the $H_{S,j}$'s has non zero slope for $j=1,\ldots,i$. 
So by localization we can assume that $X=D(x,S)$ is a disk with 
empty triangulation (cf. \eqref{Localization of sol}). 
This is done in \cite[Thm.4.7 iii), iv)]{NP-I}. 
}\fi
The proof is an induction on $i$. 
Namely $\R_{S,1}(-,\Fs)=H_{S,1}(-,\Fs)$ has the 
concavity property iv) of Remark 
\ref{remark : Linear properties of Gamma'} 
outside $\Gamma_S$, so 
\eqref{eq : gamma'=Galla cum Gamma_i} holds for $i=1$. 
Assume inductively that \eqref{eq : gamma'=Galla cum Gamma_i} 
holds for $i-1$. Again $H_{S,i}$ has the 
concavity property  
outside $\Gamma_{S,i-1}'(\Fs)$, so 
\eqref{eq : gamma'=Galla cum Gamma_i} holds for $i$.
\end{proof}

\begin{remark}
In section \ref{An operative description of Gamma'}
we give another, more operative, 
description of $\Gamma_{S,i}'(\Fs)$.
\end{remark}

\subsection{Weak super-harmonicity of partial heights: on the skeleton}
\label{Weak super-harmonicity of partial heights: on the skeleton}

In this section and the next one, we are interested in super-harmonicity properties of the 
partial heights $H_{S,i}(-,\Fs)$. Here, for a point~$x$ belonging to the skeleton of the curve, we prove that the Laplacians are bounded above by quantities involving the rank of~$\Fs$ and the geometry of the curve.

Let us begin with a preparatory lemma.

\begin{lemma}\label{lem:invarS}
Let $S$ and $S'$ be pseudo-triangulations of $X$.
Let~$x\in \Gamma_{S}\cap\Gamma_{S'}$. Let $i \in \{1,\dotsc,r\}$. Then, the index~$i$ is spectral non-solvable at~$x$ with respect to~$S$ if, and only if, it is spectral non-solvable at~$x$ with respect to~$S'$. Moreover, in this case, we have
\begin{equation}\label{eq:wshA1}
dd^c \Rc_{S',i}(x,\Fs) -  N_{S'}(x) \;=\; dd^c \Rc_{S,i}(x,\Fs)  - N_{S}(x)\;.
\end{equation}
\end{lemma}
\begin{proof}
Since $x\in \Gamma_{S}\cap\Gamma_{S'}$, the equivalence follows directly from the definitions.

Let us now prove \eqref{eq:wshA1}. We may assume that~$K$ is algebraically closed. Let~$b$ be a branch out 
of~$x$. If~$b$ belongs to the complement of 
$\Gamma_S\cup\Gamma_{S'}$, or to $\Gamma_{S}\cap\Gamma_{S'}$,
then we have $\partial_{b} \Rc_{S',i}(x,\Fs) = 
\partial_{b} \Rc_{S,i}(x,\Fs)$.

Assume that~$b$ belongs to~$\Gamma_{S'}$, but not 
to~$\Gamma_{S}$. Remark that we are reduced to computing radii on 
an annulus. By~\cite[Lemma 3.15 (b)]{NP-II}, we have 
$\partial_{b} \Rc_{S',i}(x,\Fs) = \partial_{b} \Rc_{S,i}(x,\Fs)+1$. 

Analogously if $b$ belongs to~$\Gamma_{S}$, but not 
to~$\Gamma_{S'}$, then 
$\partial_{b} \Rc_{S',i}(x,\Fs) = \partial_{b} \Rc_{S,i}(x,\Fs)-1$. 

The result follows by summing up all the contributions.
\end{proof}



We first study the points that lie on the skeleton of the curve.

\begin{proposition}\label{prop:wshA1}
Assume that~$X$ is an analytic domain of~$\E{1}{K}$. Let~$S$ be a 
pseudo-triangulation of $X$. 
Let $i\in\{1,\dotsc,r\}$. For every~$x\in \Gamma_{S} \cap \Int(X)$, we have
\begin{equation}\label{eq:wshA1}
dd^c H_{S,i}(x,\Fs)\;\le\; (N_{S}(x)-2\deg(x)) \cdot \min(i,i_x^{\mathrm{sp}})\;.
\end{equation}
Moreover equality holds if $i$ is a vertex free of solvability at $x$.
\end{proposition}
\begin{proof}
We may assume that~$K$ is algebraically closed. 

\medbreak

\noindent $\bullet$ Assume that $i\le i_x^{\mathrm{sp}}$. 

Then, for all $j\le i$, $\Rc_{S,j}(x,\Fs)$ is spectral non-solvable. As a consequence, by Lemma~\ref{lem:invarS}, we may replace~$S$ by any other pseudo-triangulation whose skeleton contains~$x$. In particular, we may enlarge~$S$ and assume that~$x$ belongs to~$S$ and is not an end-point of~$\Gamma_{S}$. 

In this situation, there exists an affinoid neighborhood~$Y$ of~$x$ in~$X$ such that $\Gamma_{S} \cap Y$ is the analytic skeleton of~$Y$, associated to its smallest triangulation~$S_{Y}$. (To construct~$Y$, one may proceed as follows. Let $[x,x_{1}],\dotsc,[x,x_{n}]$ be the edges of~$\Gamma_{S}$ containing~$x$. The affinoid domain~$Y$ may be chosen as the connected component of $X - \{x_{1},\dotsc,x_{n}\}$ containing~$x$. We then have $S_{Y}= \{x_{1},\dotsc,x_{n}\}$.) 

By \cite[Proposition~2.8.2]{NP-III}, the radii of~$\Fs$ remain unchanged by localization to~$Y$, hence we may assume that~$X$ is an affinoid domain of the affine line and that~$S$ is its smallest triangulation. The result then follows from~\cite[Theorem~3.9]{NP-I}.


%

\medbreak

\noindent $\bullet$ Assume that $i>i_x^{\mathrm{sp}}$. 

Since $x\in\Gamma_S$, we have~$\Rc_{S,i}(x,\Fs)=1$, 
hence $\partial_{b} \Rc_{S,i}(x,\Fs)\le 0$ for every branch~$b$ out 
of~$x$, and $dd^c \Rc_{S,i}(x,\Fs) \le 0$. The result for~$i$ now follows from the result for~$i_x^{\mathrm{sp}}$.

\end{proof}


We will now extend this result to arbitrary curves. To do so, we will map the 
curve to the affine line by a finite \'etale morphism and use the previous 
proposition. We will need to understand how the radii of convergence change 
in the process. This is easily described when the degree of the morphism is prime to the residual exponent, as we already investigated in~\cite{NP-II} (for weak triangulations, but the results and the proofs immediately adapt to pseudo-triangulations).


\begin{lemma}[\protect{\cite[Lemma 3.23]{NP-II}}]\label{lem:tameradii}
Let~$Y$ and~$Z$ be quasi-smooth $K$-analytic curves with pseudo-triangulations~$S$ and~$T$ respectively. Let $\varphi\colon Y\to Z$ be a finite \'etale morphism. Let $\Fs_{Y}$ be a module with connection of rank~$r$ on~$Y$. Let $y\in\Gamma_S \cap \varphi^{-1}(\Gamma_T)$. 
Assume that $d_{y} := [\Hs(y)\colon\Hs(\varphi(y))]$ 
is prime to $p$. 
Then for each $i\in \{1, \ldots, r\}$ and $j\in\{1,\ldots,d_y\}$, we have
$\R_{T,d_y(i-1)+j}(\varphi(y),\varphi_*\Fs_{Y})=\R_{S,i}(y,\Fs_{Y})$:
\begin{equation}\label{eq : nice formula of radii tame case}
\RR_{T}(\varphi(y),\varphi_*\Fs_{Y})\;=\;
\bigl(\underbrace{\R_{S,1}(y,\Fs_{Y}),\ldots,
\R_{S,1}(y,\Fs_{Y})}_{d_y \textrm{ times}}\;,\;\dotsc\;,\;
\underbrace{\R_{S,r}(y,\Fs_{Y}),\ldots,
\R_{S,r}(y,\Fs_{Y})}_{d_y \textrm{ times}}\bigr)\;.\quad\Box
\end{equation}
\end{lemma}

\begin{corollary}\label{cor:radiiCC'}
Let $\psi \colon C \to C'$ be a finite \'etale morphism between open annuli. Endow~$C$ and~$C'$ with the empty pseudo-triangulations. Let $\Fs_{C}$ be a module with connection of rank~$r$ on~$C$ and assume that its radii of convergence are $\log$-linear on~$\Gamma_{C}$ with slopes $\lambda_{1},\dotsc,\lambda_{r}$. Assume that $\deg(\psi)$ is prime to~$p$. Then, the radii of convergence of $\psi_{*} \Fs_{C}$ are $\log$-linear on~$\Gamma_{C'}$ with slopes 
\begin{equation}\label{eq : nice formula of slopes tame case}
\Bigl(\underbrace{\frac{1}{\deg(\psi)}\,\lambda_{1},\dotsc,
\frac{1}{\deg(\psi)}\,\lambda_{1}}_{\deg(\psi) \textrm{ times}}\;,\;
\ldots\;,\;
\underbrace{\frac{1}{\deg(\psi)}\,\lambda_{r},\ldots,
\frac{1}{\deg(\psi)}\,\lambda_{r}}_{\deg(\psi) \textrm{ times}}\Bigr)\;.
\end{equation}
\end{corollary}
\begin{proof}
It follows from the assumption that, for each $y \in \Gamma_{C}$, we have $\psi(y)\in \Gamma_{C'}$ and $[\Hs(y)\colon\Hs(\varphi(y))] = \deg(\psi)$. The result now follows from Lemma~\ref{lem:tameradii}, combined with~\cite[Thm~4.4.33]{Duc} in order to take into account the dilatation of distances induced by~$\psi$. 

%
%
\end{proof}

\begin{corollary}\label{cor:tameslopes}
Let~$Y$ and~$Z$ be quasi-smooth 
$K$-analytic curves with pseudo-triangulations~$S$ and~$T$ respectively. Let $\varphi\colon Y\to Z$ be a finite \'etale morphism. Let $\Fs_{Y}$ be a module with connection of rank~$r$ on~$Y$. Let $y\in\Gamma_S \cap \varphi^{-1}(\Gamma_T)$ such that $\varphi^{-1}(\varphi(y))=\{y\}$. Let~$c$ be a branch out of~$\varphi(y)$ and let $b_{1},\dotsc,b_{t}$ be the branches out of~$y$ over~$c$. Assume that $c \subset \Gamma_{T}$ and $b_{1},\dotsc,b_{t} \subset \Gamma_{S}$. Assume that, for each $i \in \{1,\dots,t\}$, $\deg(\varphi_{\vert b_{i}})$ is prime to~$p$. Then, we have 
\begin{equation}
\partial_{c}H_{T,r \cdot \deg(\varphi)}(\varphi(y),\varphi_*\Fs_{Y})\;=\; \sum_{i=1}^t \partial_{b_{i}}H_{S,r}(y,\Fs_{Y})
\;.
\end{equation}
\qed
\end{corollary}
\begin{proof}
We may assume that~$K$ is algebraically closed. 
There exists a section~$C_{c}$ of~$c$ and sections~$C_{i}$'s of the~$b_{i}$'s such that~$\varphi$ induces finite \'etale morphisms $\psi_{i} \colon C_{i} \to C_{c}$ of degree~$\deg(\varphi_{\vert b_{i}})$. We may assume that all those sections are open annuli and that the radii of convergence of~$\Fs_{Y}$ and $\varphi_{*} \Fs_{Y}$ are $\log$-linear on the corresponding skeletons. Let~$U$ denote the union of the~$C_{i}$'s. We have $\varphi_{*}(\Fs_{Y})_{|U} = \bigoplus_{1\le i\le t} (\psi_{i})_{*}(\Fs_{Y})_{|C_{i}}$ and the result now follows from Corollary~\ref{cor:radiiCC'}.
\end{proof}

\if{\begin{corollary}\label{cor:tameslopes}
Assume that~$K$ is algebraically closed. Let~$Y$ and~$Z$ be quasi-smooth 
$K$-analytic curves with weak triangulations~$S$ and~$T$ respectively. Let 
$f\colon Y\to Z$ be a finite \'etale morphism. 
Let $y\in\Gamma_S \cap f^{-1}(\Gamma_T)$ such that $f^{-1}(f(y))=\{y\}$.
Let~$c$ be a branch out of~$f(y)$ that belongs to~$\Gamma_{T}$. 
Let $b_{1},\dotsc,b_{t}$ be the branches out of~$y$ over~$c$. 
Assume that they all belong to~$\Gamma_{S}$. 
Let $j\in\{1,\dotsc,t\}$. 
Among the slopes of $\{\partial_{c}\R_{T,i}(f(y),f_*\Fs)\}_{i=1,
\dotsb,\deg(f)}$, there are~$r\cdot d_{b_j}$ 
slopes (counted with multiplicity) 
coming from the branch $b_j$.


If $d_{b_{j}}=\deg(f_{|b_j})$ is prime to $p$, then those slopes are
\begin{equation}\label{eq : nice formula of slopes tame case}
\Bigl\{\underbrace{\frac{1}{d_{b_{j}}}\,\partial_{b_{j}}\R_{S,1}(y,\Fs),\dotsc,
\frac{1}{d_{b_{j}}}\,\partial_{b_{j}}\R_{S,1}(y,\Fs)}_{d_{b_{j}} \textrm{ times}}\;,\;
\ldots\;,\;
\underbrace{\frac{1}{d_{b_{j}}}\,\partial_{b_{j}}\R_{S,r}(y,\Fs),\ldots,
\frac{1}{d_{b_{j}}}\,\partial_{b_{j}}\R_{S,r}(y,\Fs)}_{d_{b_{j}} \textrm{ times}}\Bigr\}\;.
\end{equation}
In particular, if all the $d_{b_{j}}$'s are prime to~$p$, then one has 
\begin{equation}
\partial_bH_{S,r}(y,\Fs)\;=\;
\partial_{f(b)}H_{S,r \cdot \deg(f)}(f(y),f_*\Fs)\;.
\end{equation}
\end{corollary}
\begin{proof}
By the above property iv), the result follows from 
Lemma~\ref{lem:tameradii} and~\cite[Thm~4.4.33]{Duc} in order to take 
into account the dilatation of distances induced by the finite map~$f$. 
%
%
\end{proof}
}\fi

We need to find conditions ensuring that we can apply the previous results.

\begin{definition}\label{def:TR}
Let~$x$ be a point of~$X$ of type~2. If~$K$ is algebraically closed, we say 
that the point~$x$ satisfies the condition~$(TR)$ if there exists a 
finite morphism 
\begin{equation}
f\;\colon\; \Cs_{x} \xrightarrow{\quad}\mathbb{P}^1_{\wti K}
\end{equation} 
that is tamely ramified everywhere and unramified almost 
everywhere, \textit{i.e.} the 
degree of~$f$ is prime to the characteristic exponent of~$\wti K$  at a finite 
number of 
closed point of~$\Cs_{x}$ and equal to~$1$ at every other. 

In general, we say that the point~$x$ satisfies the condition~$(TR)$ if 
one (or equivalently all) of its preimages on~$X_{\widehat{K^\alg}}$ satisfies the 
condition~$(TR)$.





\end{definition}

\begin{proposition}\label{prop:CSTR}
Let~$x$ be a point of~$X$ of type~2. 
If~$g(x)=0$ or if $\charac(\wti K)\ne 2$, then~$x$ satisfies the condition~$(TR)$.

\end{proposition}
\begin{proof}
We may assume that~$K$ is algebraically closed. If~$g(x)=0$, then~$\Cs_x$ is isomorphic to~$\mathbb{P}^1_{\wti K}$ and the result is obvious. If~$\charac(\wti K)\ne 2$, then the result follows from~\cite[Prop.~8.1]{Fulton}.
\end{proof}

\begin{proposition}\label{prop:ddcHr}
Let~$x$ be a point in $\Gamma_{S} \cap \Int(X)$ of type~2 that satisfies the 
condition~$(TR)$. Assume that all the radii are spectral 
non-solvable at $x$. Then, we have
\begin{equation}
dd^c H_{S,r}(x,\Fs)\; =\; (2g(x) - 2\deg(x) + N_{S}(x)) \cdot r\;.
\end{equation}
\end{proposition}
\begin{proof}
We may assume that~$K$ is algebraically closed. 
By assumption, there exists 
a finite morphism $f\colon \Cs_{x} \to \mathbb{P}^1_{\wti K}$ that is 
tamely ramified everywhere and unramified almost everywhere. We may lift it 
to a morphism $\varphi\colon Y\to W$, where~$Y$ is an affinoid 
neighborhood of~$x$ in~$X$ and~$W$ an affinoid domain 
of~$\mathbb{P}^{1,\an}_{K}$. 
By restricting~$Y$, we may assume that 
$\varphi^{-1}(\varphi(x))=\{x\}$.
By~\cite[Thm.~4.3.15]{Duc}, the degree of 
the restriction of~$\varphi$ to any branch of~$Y$ out of~$x$ is prime 
to~$p$ and it is equal to~$1$ for almost all of them.

By Lemma~\ref{lem:invarS}, the result does not depend on the chosen 
triangulation on~$X$. We may endow~$Y$ and~$W$ with 
triangulations~$S$ and~$T$ respectively such that 
\begin{enumerate}
\item $\varphi^{-1}(\Gamma_{T})=\Gamma_{S}$ ;
\item $\Gamma_{S}$ contains all the branches~$b$ out of~$x$ with $\deg(\varphi_{\vert b}) \ne 1$ ;
\item the radii are locally constant outside~$\Gamma_{S}$ 
and~$\Gamma_{T}$. 
\end{enumerate}

By applying Corollary~\ref{cor:tameslopes}
to each branch of~$\Gamma_{S}$ out of~$x$, we find
\begin{equation}
dd^c H_{S,r}(x,\Fs) \;=\; dd^c H_{T,dr}(\varphi(x),\varphi_{\ast}\Fs), 
\end{equation}
where $d=[\H(x):\H(\varphi(x))]$. Moreover, by Proposition~\ref{prop:wshA1}, we have
\begin{equation}
dd^c H_{T,dr}(\varphi(x),\varphi_{\ast}\Fs)\; =
\; (N_{T}(\varphi(x))-2\deg(x)) \cdot d\cdot r.
\end{equation}
Recall that we denote by $\Bs(x,\Gamma_{S})$ the set of branches of~$\Gamma_{S}$ out of~$x$, see~\eqref{eq : B(x,Gamma)}. We have
\begin{align}
N_{T}(\varphi(x)) \cdot d & = \sum_{b\in \Bs(x,\Gamma_{S})} \deg(\varphi_{|b})\\
& = N_{S}(x) +  \sum_{b\in \Bs(x,\Gamma_{S})} (\deg(\varphi_{|b})-1)\\
& = N_{S}(x) +  \sum_{b\in \Bs(x,X)} (\deg(\varphi_{|b})-1),
\end{align}
since we chose~$\Gamma_{S}$ so that it contains every branch out of~$x$ where the degree of~$\varphi$ is not~1.

Recall that there is a bijection 
\[\begin{array}{ccc}
\Bs(x,X) & \longrightarrow & \Cs_{x}(\tilde K)\\
b & \longmapsto & P_{b}
\end{array}\]
between the set of branches out of~$x$ and the set of closed points of the residual curve~$\Cs_{x}$ (see \cite[4.2.11.1]{Duc}). Moreover, for each branch~$b$ out of~$x$, we have
\begin{equation}
\deg(\varphi_{\vert b}) = e_{P_{b}}(\tilde \varphi),
\end{equation}
where $e_{P_{b}}(\tilde \varphi)$ denotes the ramification index of the residual morphism $\tilde \varphi \colon \Cs_{x} \to \Cs_{\varphi(x)}$ at~$P_{b}$ (see \cite[4.3.13]{Duc}).
We deduce that 
\begin{equation}
N_{T}(\varphi(x)) \cdot d \;=\; N_{S}(x) +  \sum_{P\in \Cs_{x}(\tilde K)} (e_{P}(\tilde \varphi)-1) \;=\; 2d+2g(x)-2,
\end{equation}
by the Riemann-Hurwitz formula (cf. \cite[Cor. 2.4]{Hartshorne}).
The result follows.
\end{proof}

\begin{theorem}\label{thm:curveTR}
Let~$x$ be a point in $\Gamma_{S} \cap \Int(X)$. If it is of type~2, assume that it satisfies the condition~$(TR)$. Then, for each $i\in\{1,\dotsc,r\}$, we have
\begin{equation}\label{eq:wsh}
dd^c H_{S,i}(x,\Fs)\; \le\; 
(2g(x) - 2 \deg(x) + N_{S}(x)) \cdot \min(i,i_x^{\mathrm{sp}}).
\end{equation}
Moreover equality holds if $i$ is a vertex free of solvability at $x$.
\end{theorem}
\begin{proof}
We may assume that $K$ is algebraically closed. If~$x$ is of type~3, then it has a neighborhood that is isomorphic to an annulus, and the result follows from Proposition~\ref{prop:wshA1}. As a consequence, we may assume that $x$ is of type~2.

As in the proof of Proposition~\ref{prop:wshA1}, the result for $i>i_x^{\mathrm{sp}}$ follows from that for $i=i_x^{\mathrm{sp}}$. As a result, we may assume that $i\le i_x^{\mathrm{sp}}$, \textit{i.e.} the radius 
$\Rc_{S,i}(x,\Fs)$ is spectral non-solvable. 

As in the proof of 
Proposition~\ref{prop:wshA1}, we may localize in the neighborhood 
of~$x$. Recall that, by Dwork-Robba's theorem~\cite[Theorem~4.1.1]{NP-III}, the 
differential module~$\Fs_{x}$ may be written as a direct sum $\Fs_{1,x} 
\oplus\cdots \oplus \Fs_{s,x} \oplus \Fs_{s+1,x}$ where, for each $k\in\{1,\dotsc,s\}$, 
the radii of $\Fs_{k,x}$ at~$x$ are spectral non-solvable and equal, and the 
radii of~$\Fs_{s+1,x}$ at~$x$ are all equal to~1. The result for~$\Fs_{x}$ follows from the result for the~$\Fs_{k,x}$'s, hence we may assume 
that all the radii $\Rc_{S,i}(x,\Fs)$'s are equal.


Let $j \le k\in\{1,\dotsc,r\}$. Since $\Rc_{S,j}(x,\Fs) = \Rc_{S,k}(x,\Fs)$, for 
every branch~$b$ out of~$x$, we have $\partial_{b} \Rc_{S,j}(x,\Fs) \le 
\partial_{b} \Rc_{S,k}(x,\Fs)$, hence $dd^c \Rc_{S,j}(x,\Fs)\le dd^c 
\Rc_{S,k}(x,\Fs)$. Consider the polygon~$\Ps$ whose vertices are the points 
\begin{equation}
\Big(i,dd^c H_{S,i}(x,\Fs) = \sum_{j=1}^i dd^c \Rc_{S,j}(x,\Fs)\Big),
\end{equation} 
for $i\in\{1,\dotsc,r\}$. The previous inequalities 
show that~$\Ps$ is concave, hence we have
$dd^c H_{S,i}(x,\Fs)/i \le dd^c H_{S,r}(x,\Fs)/r$. 
By Proposition~\ref{prop:ddcHr}, we have $
dd^c H_{S,r}(x,\Fs)/r \le 2g(x) - 2\deg(x) + N_{S}(x)$, and the result follows.
\end{proof}


\begin{remark}\label{rem:noTR}
Formula~\eqref{eq:wsh} first appears in \cite[Theorem~5.3.6]{Kedlaya-draft}, with a rather terse proof, and the authors felt the need to provide a more complete argument, under the additional $(TR)$ assumption. After a first version of this text was made available, the V.~Bojkovi\'c and the first-named author undertook a systematic study of the behaviour of radii of convergence under pushforward by finite \'etale morphisms, building on M.~Temkin's theory of ramification for morphisms of curves (see~\cite{TemkinMetricUniformization}). As a consequence, they manage to remove the assumption~$(TR)$ from Theorem~\ref{thm:curveTR}, thus proving Formula~\eqref{eq:wsh} in full generality.
\end{remark}

\subsection{Weak super-harmonicity of partial heights: outside the skeleton}
\label{Weak super-harmonicity of partial heights: outside the skeleton}

We pursue our study of super-harmonicity properties of the partial heights $H_{S,i}(-,\Fs)$ by investigating the situation at points outside the skeleton. 




\begin{definition}\label{def:Ei}
For each $i\in\{1,\dotsc,r\}$, set
\begin{equation}
\Es_{S,i}(\Fs) \;:=\; \{x\in X \mid dd^c H_{S,i}(x,\Fs)>0\}\;.
\end{equation}
In the sequel, we may write $\Es_{S,i}:=\Es_{S,i}(\Fs)$ for short.
\end{definition}


\begin{definition}\label{Def. C_i}
We define inductively a sequence of locally finite sets 
\begin{equation}\label{eq : C_i}
\C_{S,1}(\Fs)\;\subseteq \;\ldots\;\subseteq\; 
\C_{S,r}(\Fs)\;\subseteq \;X
\end{equation} 
as follows. Set $\C_{S,1}(\Fs) := \aleph_{S,1}(\Fs):= \emptyset$. For $i \in \{2,\dotsc,r\}$, let $\aleph_{S,i}(\Fs)$ be the locally 
finite set of points $x \in X -\Gamma_{S}$ satisfying 
\begin{enumerate}
\item $\R_{S,i}(-,\Fs)$ is solvable at $x$;
\item $x$ is an end-point of $\Gamma_{S,i}(\Fs)$;
\item $x\in\Gamma_{S,i-1}'(\Fs)\cap \Gamma_{S,i}(\Fs)\cap 
\Gamma_S(H_{S,i}(-,\Fs))$
\end{enumerate}
and let
\begin{equation}
\C_{S,i}(\Fs)\;:=\; \Cs_{S,i-1}(\Fs) \cup \aleph_{S,i}(\Fs)\;.
\end{equation}
In the sequel, we may write $\aleph_{S,i} := \aleph_{S,i}(\Fs)$ and
$\C_{S,i}:=\C_{S,i}(\Fs)$ for short.
\end{definition}


\begin{remark}\label{Rk : C_Sis position}\label{rem:cardinalCi}
Let $i\in\{1,\dotsc,r\}$. By definition, we have 
\begin{equation}
\C_{S,i}(\Fs)\cap \Gamma_S\;=\;\emptyset\;\textrm{ and }\;
\C_{S,i}(\Fs)\;\subseteq\;\Gamma_{S,i-1}'(\Fs)\;.
\end{equation}
The graph $\Gamma_{S,i}(\Fs) \cap 
\Gamma'_{S,i-1}(\Fs)$ contains $\Gamma_{S}$ and the points 
of~$\aleph_{S,i}(\Fs)$ are some of its end-points. We deduce that the 
cardinal of~$\aleph_{S,i}(\Fs)$ is at most the number of end-points 
of~$\Gamma'_{S,i-1}(\Fs)$ that do not belong to~$\Gamma_{S}$.
\end{remark}

\begin{remark}\label{rem:CSiwKa}
Let $i\in \{1,\dotsc,r\}$. If follows from the behaviour of triangulations under scalar extension (see Section~\ref{section:extensionofscalars}) and the invariance of radii that we have $\aleph_{S_{\wKa},i}(\Fs_{\wKa}) = \pi_{\wKa}^{-1}(\aleph_{S,i}(\Fs))$ and $\Cs_{S_{\wKa},i}(\Fs_{\wKa}) = \pi_{\wKa}^{-1}(\Cs_{S,i}(\Fs))$.
\end{remark}


The sets~$\aleph_{S,i}$ and~$\Cs_{S,i}$ have been precisely investigated by the second-named author in the case of the affine line (see~\cite[Theorem~3.9 (iv)]{NP-I}). We generalize his results to arbitrary curves.

\begin{theorem}\label{thm:Ei}
Let $i\in\{1,\dotsc,r\}$. For each $x \notin (S \cup \C_{S,i}(\Fs))$, we have
\begin{equation}\label{eq:sh-SCi}
dd^c H_{S,i}(x,\Fs)\;\le\; 0.
\end{equation}
Moreover, equality holds in \eqref{eq:sh-SCi} if $i$ is a vertex free of 
solvability at $x$. 
\end{theorem}
\begin{proof}
By Remark~\ref{rem:CSiwKa}, we may assume that $K$ is algebraically closed. Let~$x\in X\setminus (S \cup \C_{S,i})$. 

Let us first assume that~$x \in \Gamma_{S}$. Since~$x\notin S$, we have $g(x)=0$ and $N_{S}(x)=2$, hence the result follows 
from Theorem~\ref{thm:curveTR} and Proposition~\ref{prop:CSTR} (or Remark~\ref{rem:noTR}).

Let us now assume that~$x\notin \Gamma_{S}$. 
Let~$D$ be the connected component 
of~$X\setminus \Gamma_{S}$ containing~$x$. It is an open disk. Let us 
identify it with~$D^-(0,R)$ for some~$R>0$ and endow it with the empty 
triangulation. Remark that, for each~$i\in\{1,\dotsc,r\}$, we have $\Rc_{S,i}(-,\Fs)_{|D} = \Rc_{\emptyset,i}(-,\Fs_{|D})$. 

We will now consider the embedded radii in the 
sense of~\cite[Section~2.4.2]{NP-II}. 
By \cite[Formula~(2.3)]{NP-II}, for every $i\in\{1,\dotsc,r\}$, we have
\begin{equation}\label{eq:radiusembS0}
\Rc_{\emptyset,i}(x,\Fs_{|D})\; = \;
\frac{\Rc_{i}^\textrm{emb}(x,\Fs_{|D})}{R}.
\end{equation}

Let $i_{0} \in \{0,r\}$ such that, for each $i \in \{1,\dotsc,i_{0}\}$, we have $\Rc_{i}^\emb(x,\Fs_{|D}) < R$. There exists $R' \in \of{]}{0,R}{[}$ such that, for each $i \in \{1,\dotsc,i_{0}\}$, we have $\Rc_{i}^\emb(x,\Fs_{|D}) < R'$. Set $D' := D^+(0,R')$. Up to enlarging~$R'$, we may assume that $x \in D^-(0,R')$. Note that, for each $z\in D^-(0,R')$ and each $i \in \{1,\dotsc,i_{0}\}$, we have $\Rc_{i}^\emb(z,\Fs_{|D})<R'$ (as otherwise, we would have $\Rc_{i}^\emb(x,\Fs_{|D}) = \Rc_{i}^\emb(z,\Fs_{|D})\ge R'$), hence
\begin{equation}
\Rc_{i}^\emb(z,\Fs_{|D}) \;=\; \Rc_{i}^\emb(z,\Fs_{|D'}).
\end{equation}
As a consequence, for each $i \in \{1,\dotsc,i_{0}\}$, on $D^-(0,R')$, the radius $\Rc_{i}^\emb(-,\Fs_{|D'})$ coincides, up to a constant multiplicative factor, with the radii $\Rc_{\emptyset,i}(-,\Fs_{|D})$ and $\Rc_{S,i}(-,\Fs)_{|D}$. In particular, the associated Laplacian coincide, as well as the sets $\C_{S,i}(\Fs)\cap D^-(0,R')$ and $\C_{i} \cap D^-(0,R')$, where $\C_{i}$ is computed using the radii $\Rc_{j}^\emb(z,\Fs_{|D'})$, as in~\cite{NP-I} (which deals with the case of an affinoid domain of $\E{1}{K}$). The result now follows from~\cite[Theorem~3.9]{NP-I}.

Let $j\in \{i_{0}+1,\dotsc,r\}$. By definition, we have $\Rc_{j}^\emb(x,\Fs_{|D}) = R$, hence the map $\Rc_{j}^\emb(-,\Fs_{|D})$ is constant on~$D$ (with value~$R$). As a consequence, the map $\Rc_{S,j}(-,\Fs)$ is constant on~$D$ too, and we have $dd^c \Rc_{S,j}(x,\Fs)=0$ and $\aleph_{S,j}(\Fs)=\emptyset$. 

If $i_{0}=0$, we deduce that $dd^c H_{S,i}(x,\Fs) = 0$, for each $i\in \{i_{0}+1,\dotsc,r\}$.

If $i_{0}\ge 1$, we deduce that, for each $i\in \{i_{0}+1,\dotsc,r\}$, we have $\Cs_{S,i}(\Fs) =  \Cs_{S,i_{0}}(\Fs)$ and $dd^c H_{S,i}(x,\Fs) = dd^c H_{S,i_{0}}(x,\Fs)$, hence the result follows from the previous case.
\end{proof}

\begin{corollary}\label{cor:R1sh}
For each $x \notin S$, we have
\begin{equation}\label{eq:shR1}
dd^c \Rc_{S,1}(x,\Fs)\;\le\; 0.
\end{equation}

For each $i\in\{1,\dotsc,r\}$, we have $\Es_{S,i}(\Fs) \subseteq S \cup \C_{S,i}(\Fs)$. In particular, the set~$\Es_{S,i}$ is locally finite.
\qed
\end{corollary}

To complete the results obtained so far, we provide upper bounds for the Laplacians of the partial heights. 

\begin{proposition}\label{prop:HiCi}
Let $x\in X -\Gamma_{S}$. We have
\begin{equation}\label{eq:ddcR1}
dd^c \Rc_{S,1}(x,\Fs) \le 0
\end{equation}
and, for each $i\in\{2,\dotsc,r\}$ such that $\Rc_{S,i}(-,\Fs)$ is solvable at~$x$, we have
\begin{equation}\label{eq:ddcRisolvable}
dd^c \Rc_{S,i}(x,\Fs) \le 1.
\end{equation}
For each $i\in\{1,\dotsc,r\}$, we have
\begin{equation}
dd^c H_{S,i}(x,\Fs) \le i-1\;.
\end{equation}
\end{proposition}
\begin{proof}
We may assume that~$K$ is algebraically closed. Since~$x\notin \Gamma_{S}$, the connected component of $X-\Gamma_{S}$ containing~$x$ is an open disk. Up to replacing~$X$ by this connected component, we may assume that~$X$ is an open disk endowed with the empty pseudo-triangulation.

$\bullet$ Equation \eqref{eq:ddcR1} follows from Corollary~\ref{cor:R1sh}.


\medbreak

$\bullet$ Let $i\in\{2,\dotsc,r\}$ such that $\Rc_{\emptyset,i}(-,\Fs)$ is solvable at~$x$. 

Let~$b$ be a germ of segment out of~$x$. If $b$ points towards the interior of 
the disk, then $\partial_{b} \Rc_{\emptyset,i}(x,\Fs)\leq 0$. Indeed, otherwise, there would exist a point~$y$ in the connected component of $X-\{x\}$ containing~$b$ such that $\Rc_{\emptyset,i}(y,\Fs) > \Rc_{\emptyset,i}(x,\Fs)$. As a consequence, $\Fs$ admits at least $r-i+1$ linearly independent solutions at~$y$ with radii of convergence bigger than~$\Rc_{\emptyset,i}(x,\Fs)$. Since $\Rc_{\emptyset,i}(-,\Fs)$ is solvable at~$x$, those solutions convergence in the neighborhood of~$x$, hence $\Rc_{\emptyset,i}(-,\Fs)$ is over-solvable at~$x$, and we get a contradiction.


If $b$ is the germ of segment pointing towards the exterior of the disk, then 
$\partial_{b} \Rc_{\emptyset,i}(x,\Fs) \leq 1$. Indeed, otherwise, as above, we would have as over-solvability at~$x$. 

Summing up the contributions, we deduce that $dd^c \Rc_{\emptyset,i}(x,\Fs) \le 1$.

\medbreak

$\bullet$ Let $i\in\{1,\dotsc,r\}$. 

Assume that $i=1$. Then the result follows from~\eqref{eq:ddcR1}.


Assume that $i \le i_{x}^\sp$. By definition, for each $j\le i_{x}^\sp$, we have $x \notin \aleph_{S,j}(\Fs)$, hence $x\notin \Cs_{S,i}(\Fs)$. By Theorem~\ref{thm:Ei}, we then have $dd^c H_{S,i}(x,\Fs) \le 0$. 

Assume that $i> \max(i_{x}^\sp,1) =: i_{0}$. 
We have just proved that  $dd^c H_{S,i_{0}}(x,\Fs) \le 0$, so it is enough to prove that, for all $j > i_{0}$, we have 
$dd^c \Rc_{\emptyset,j}(x,\Fs) \le 1$. If $\Rc_{S,j}(-,\Fs)$ is solvable at~$x$, this follows from~\eqref{eq:ddcRisolvable}. If $\Rc_{S,j}(-,\Fs)$ is over-solvable at~$x$, then we even have $dd^c \Rc_{\emptyset,j}(x,\Fs) =0$.
\end{proof}

\subsection{An operative description of $\Gamma_{S,i}'(\Fs)$}
\label{An operative description of Gamma'}


\begin{notation}\label{notation : D_b - first}
Let $\Gamma\subseteq X$ be a locally finite 
graph.\footnote{Recall that $X-\Gamma$ is a disjoint union of 
virtual disks, cf. \ref{Section : Graphs.}.} 
Let $x\in\Gamma$, and let $b\notin\Gamma$ be a germ of segment 
out of $x$. We denote by $D_b\subset X$ the virtual open disk with 
boundary $x$ containing $b$. 
\end{notation}
\begin{proposition}\label{prop psedudo-harm}
Let $i\leq r$, let $\Gamma\subseteq X$ be a locally finite graph 
containing $\Gamma_S$. Let $x\in\Gamma$. 
The following conditions are equivalent:
\begin{enumerate}
\item One has $\Gamma_{S,i}'(\Fs)\cap D_b=\emptyset$, for all 
direction $b$ out of $x$ such that $b\notin\Gamma$ 
(i.e. $\Gamma_{S,i}'(\Fs)\subseteq\Gamma$ around $x$).
\item
For all $j=1,\ldots,i$, and for all germs of 
segment $b$ out of $x$ such that $b\notin\Gamma$, one has 
\begin{equation}\label{eq : notion fff}
\partial_b\R_{S,j}(x,\Fs)\;=\;0 \qquad
(\textrm{resp. }\partial_bH_{S,j}(x,\Fs)\;=\;0)\;.
\end{equation}
\item For all $j=1,\ldots,i$ the function $\R_{S,j}(-,\Fs)$ 
(resp. $H_{S,j}(-,\Fs)$) 
verifies 
\begin{equation}
dd^c_{\not\subseteq\Gamma}
\R_{S,j}(x,\Fs)\; =\; 0 
\qquad (\textrm{resp. } 
dd^c_{\not\subseteq\Gamma}
H_{S,j}(x,\Fs)\; 
=\; 0)\;.
\end{equation}
\item Same as iii), replacing the equalities by $\leq$.
\item Same as ii), replacing the equalities by $\leq$.
%
%
%
\end{enumerate}
In particular these equivalent conditions holds at each 
$x\in\Gamma$ if, and only if, one has 
\begin{equation}
\Gamma_{S,i}'(\Fs)\;\subseteq\;\Gamma\;.
\end{equation}
\end{proposition}
\begin{proof}
We can assume $K$ algebraically closed.
Clearly i) $\Rightarrow$ ii) $\Rightarrow$ iii) $\Rightarrow$ iv), 
and ii) $\Rightarrow$ v) $\Rightarrow$ iv).
We now prove that  iv)  imply i). 
We proceed by induction on $i$. 

If $i=1$, the first radius has the concavity property on each disk $D'$ 
such that $D'\cap\Gamma_S=\emptyset$. 
So both iv) and v) imply ii) for $\R_{S,1}(x,\Fs)$. 
By \eqref{eq : partial R_b =0 .jky fd} one has i). 

Let now $i>0$. Assume inductively that all conditions hold 
for $i-1$. Firstly notice that, for all $b\notin\Gamma$, one has
$\partial_b\R_{S,i}(x,\Fs)=\partial_bH_{S,i}(x,\Fs)$ by point ii) 
of Remark \ref{remark : Linear properties of Gamma'}.
In fact %
by i) one has 
$\Gamma_{S,i-1}(\Fs)\cap D_b=\emptyset$. 
Now, since $\R_{S,i}(-,\Fs)$ has the concavity property on $D_b$ 
(cf. point iv) of Remark \ref{remark : Linear properties of Gamma'}), 
then $\partial_b\R_{S,i}(x,\Fs)\geq 0$ for all $b\notin\Gamma$. 
Hence both v) and iv) imply ii), 
and by \eqref{eq : partial R_b =0 .jky fd} one has 
$\Gamma_{S,i}(\Fs)\cap D_b=\emptyset$, for all $b\notin\Gamma$. 
This implies i) since 
$\Gamma_{S,i}'(\Fs)=\Gamma_{S,i-1}'(\Fs)\cup\Gamma_{S,i}(\Fs)$.
\end{proof}

\begin{lemma}
If $i_x^{\mathrm{sp}}=0$ (i.e. if all the radii are solvable or 
over-solvable at $x$),  the conditions of Proposition 
\ref{prop psedudo-harm} are automatically fulfilled at $x\in \Gamma$.
\end{lemma}
\begin{proof}
Let $b\notin \Gamma$ be a germ of segment out of $x$. 
Since the first radius satisfies the concavity 
property outside $\Gamma_S$ 
(cf. point iv) of Remark \ref{remark : Linear properties of Gamma'}).
Then $\R_{S,1}(-,\Fs)$ is constant and over-solvable at all point of $D_b$. 
This will imply the same property for $\R_{S,i}(-,\Fs)$, for 
all $i\geq 1$. 
\end{proof}

\begin{proposition}\label{Prop : equi 3 cond f f^*}
Let $\Gamma\subseteq X$ be a locally finite graph containing 
$\Gamma_{S,i}'(\Fs)$. 
If the index $i$ is spectral non solvable at each point of $\Gamma$, 
then for all $j=1,\ldots,i$ one has 
$\R_{S,j}(x,\Fs)=\R_{S,j}(x,\Fs^*)$ for all $x\in X$, and 
$\Gamma_{S,j}(\Fs)=\Gamma_{S,j}(\Fs^*)$.
\end{proposition}
\begin{proof}
Since spectral radii are stable by duality, for all $j\leq i$ and all  
$x\in\Gamma_{S,i}'(\Fs)$, we have 
$\R_{S,j}(x,\Fs)=\R_{S,j}(x,\Fs^*)$. Moreover, if 
$x\in\Gamma_{S,i}'(\Fs)$, by continuity $\R_{S,j}(-,\Fs)$ 
remains spectral over all germ of segment out of $x$. 
So $\Fs^*$ verifies the assumptions of Proposition 
\ref{prop psedudo-harm} and one has 
$\Gamma_{S,j}'(\Fs^*)\subseteq \Gamma_{S,j}'(\Fs)\subseteq 
\Gamma$, for all $j=1,\ldots,i$. 
Repeating the argument for $\Fs^*$ we also have the converse of 
this inclusion. 
\end{proof}
\if{
\begin{lemma}
Let $B$ be a finite set of germs of segment out of $x\in X$, and let 
$B_{S,i}(x,\Fs)$ be the set of germs of segment out of $x$ belonging 
to $\Gamma_{S,i}(\Fs)$. Assume that
\begin{equation}
\sum_{b\in B}\partial_bH_{S,i}(x,\Fs) \;\geq\; 
(2g(x)-2+c(x))\cdot \min(i,i_x^{\mathrm{sp}})\;.
\end{equation}
Then ...
\end{lemma}
\begin{proof}

Questo Lemma discende dal corollario seguente rimpiazzando $(X,S)$ 
con un intorno $(U,S_U)$ di $x$ tale che 
$\Gamma_{S_U}=\Gamma_S\cap U$ ...

...........

????????????

???????????????

?????????

\end{proof}
}\fi

Combining 
Lemma \ref{Lemma : solvable radii at boundary}, and 
Prop. \ref{prop psedudo-harm},  one finds the following result.

\begin{definition}
Let $\Gamma\subseteq X$ be a locally finite graph containing $\Gamma_S$.
We denote by $B(\Gamma)$ %
the set of points  $x\in\Gamma$ such that either $x\in\partial X$, 
or $x$ does not satisfy the $(TR)$ condition.
\end{definition}

\begin{corollary}\label{Cor : effective computation of Gamma}
Let $i\leq r$. 
Let $\Gamma\subseteq X$ be a locally finite graph containing 
$\Gamma_S$. Assume that the following conditions hold:
\if{
Let $B(\Gamma):=\Bigl(\Gamma\cap\partial(X)\Bigr)\cup(TR)$. 
Let $B^+(\Gamma)\subset\Gamma$ be the union of 
$B(\Gamma)$ with the set of points $x\in\Gamma$ 
such that $g(x)>0$.
Assume that $p\neq 2$, and that the following 
conditions hold:
}\fi
\begin{enumerate}
\item For all $x\in B(\Gamma)$, one of the equivalent conditions of 
Proposition \ref{prop psedudo-harm} holds at $x$;
\item For all $x\in\Gamma-B(\Gamma)$, either the conditions of 
Proposition \ref{prop psedudo-harm} holds at $x$, 
or $\Gamma\neq\{x\}$ and for all 
$j=1,\ldots,\min(i,i_x^{\mathrm{sp}})$, one has 
\begin{equation}
\label{eq : condition on pt to ensure graph in Gamma}
dd^c_{\subseteq \Gamma}H_{S,j}(x,\Fs) \;\geq\; 
\left\{
\begin{array}{lcl}
0&\textrm{ if }&x\notin\Gamma_S\\
(2g(x)-2+N_S(x))\cdot j
&\textrm{ if }&x\in\Gamma_S\;.
\end{array}
\right.
\end{equation}
\end{enumerate}
Then 
\begin{equation}
\Gamma_{S,i}'(\Fs)\;\subseteq\;\Gamma\;.
\end{equation}
%
\end{corollary}

\begin{proof}
We can assume $K$ algebraically closed.
We fix $x\in\Gamma-B(\Gamma)$, and we prove that 
\eqref{eq : condition on pt to ensure graph in Gamma} 
implies that the equivalent conditions of 
Proposition \ref{prop psedudo-harm} at $x$. 

Assume that $i\leq i_x^{\mathrm{sp}}$. 
Then the indexes $j=1,\ldots,i$ are all spectral-non solvable, 
hence $x\notin\C_{S,i}$. 
Condition 
\eqref{eq : condition on pt to ensure graph in Gamma} 
implies $\sum_{b\notin\Gamma}\partial_bH_{S,j}(x,\Fs)\leq 0$, 
for all $j=1,\ldots,i$. 
So we are done by iv) of Proposition \ref{prop psedudo-harm}. 

Assume now that $i>i_x^{\mathrm{sp}}$. 
%
If $j\leq\min(i_x^{\mathrm{sp}},i)$, we 
proceed as above to prove that $\Gamma_{S,j}'(\Fs)\subseteq\Gamma$ 
around $x$. 
So we consider $i_x^{\mathrm{sp}}<j\leq i$, and, inductively, 
we assume that for all $k=1,\ldots,j-1$, 
and all $b\notin\Gamma$ one has $\partial_b\R_{S,k}(x,\Fs)=0$ (i.e. 
$\Gamma_{S,j-1}'(\Fs)\subseteq\Gamma$ around $x$). 
If $x\notin\Gamma_{S,j}(\Fs)$, then $\partial_b\R_{S,j}(x,\Fs)=0$ for 
all $b$, so we are done. 
Assume then that $x\in\Gamma_{S,j}(\Fs)$. 
In this case, since $j>i_x^{\mathrm{sp}}$ is solvable or over-solvable at 
$x$, Lemma \ref{Lemma : solvable radii at boundary} shows 
that $\partial_b\R_{S,j}(x,\Fs)=0$, for 
all $b\notin\Gamma_{S,j-1}'(\Fs)\subseteq\Gamma$. 
In particular this holds if $b\notin\Gamma$, and Proposition 
\ref{prop psedudo-harm} implies 
$\Gamma_{S,j}'(\Fs)\subseteq\Gamma$ around $x$. 
%
%
%
%
%
\if{
We fix $x\in\Gamma$, and we prove that 
\eqref{eq : condition on pt to ensure graph in Gamma} 
implies the equivalent conditions of 
Proposition \ref{prop psedudo-harm} hold at $x$. 

Assume that $i\leq i_x^{\mathrm{sp}}$. 
For all $j=1,\ldots,i$, 
Theorem \ref{Prop : Laplacian} implies 
\begin{equation}\label{eq : sdfghjkiuytgv}
\sum_{b\notin\Gamma}\partial_bH_{S,j}(x,\Fs)\;\leq \;0\;.
\end{equation}
We now prove, by induction, that
for all $j=1,\ldots,i$, and all $b\notin\Gamma$, one has 
$\partial_bH_{S,i}(x,\Fs)\geq 0$. 
For $j=1$, the radius $\R_{S,1}(-,\Fs)$ has the concavity property 
outside $\Gamma$ 
(cf. Remark \ref{remark : Linear properties of Gamma'}). 
So $\partial_bH_{S,1}(x,\Fs)\geq 0$, for all $b\notin\Gamma$. 
Together with \eqref{eq : sdfghjkiuytgv} this implies 
$\partial_bH_{S,1}(x,\Fs)= 0$, for all $b\notin\Gamma$. 
Hence, by point iii) of Prop. \ref{prop psedudo-harm}, for all 
$b\notin\Gamma$  one has 
$\Gamma_{S,1}'(\Fs)\cap D_b=\emptyset$. 
In particular $D_b\cap\mathscr{C}_2=\emptyset$.

Now let $2\leq j\leq i$, and assume inductively that for all 
$k=1,\ldots,j-1$, and all $b\notin\Gamma$ one has 
$\partial_bH_{S,k}(x,\Fs)= 0$, and 
$\Gamma_{S,k}'(\Fs)\cap D_b=\emptyset$. 
In particular $\mathscr{C}_j\subset\Gamma$. 
By Remark \ref{remark : Linear properties of Gamma'} 
the function $H_{S,j}(-,\Fs)$ behave 
as a first radius outside $\Gamma$ so an argument as above  
proves that for all $b\notin\Gamma$, one has 
$\partial_bH_{S,j}(x,\Fs)= 0$, and 
$\Gamma_{S,j}'(\Fs)\cap D_b=\emptyset$. 
In particular $\mathscr{C}_{j+1}\subseteq \Gamma$. 
This concludes the induction for $i\leq i_x^{\mathrm{sp}}$.

Assume now that $i>i_x^{\mathrm{sp}}$. If $j=1$, or if $j\leq\min(i_x^{\mathrm{sp}},i)$, we 
proceed as above. So we assume $i_x^{\mathrm{sp}}<j\leq i$, and, inductively, 
that for all $k=1,\ldots,j-1$, 
and all $b\notin\Gamma$ one has $\partial_bH_{S,k}(x,\Fs)=0$, and 
$\Gamma_{S,k}'(\Fs)\cap D_b=\emptyset$. In particular 
$\mathscr{C}_{j}\subset \Gamma$. 
The index $j>i_x^{\mathrm{sp}}$ is solvable or over-solvable at $x$. 
If $x\notin\Gamma_{S,j}(\Fs)$, then $\partial_b\R_{S,j}(x,\Fs)=0$ for 
all $b$, so we are done. 
Assume then that $x\in\Gamma_{S,j}(\Fs)$. 
In this case Lemma \ref{Lemma : solvable radii at boundary} shows 
that $\partial_b\R_{S,j}(x,\Fs)=0$, for 
all $b\notin\Gamma_{S,j-1}(\Fs)$. In particular this holds if 
$b\notin\Gamma$. This concludes the proof in the case $i>i_x^{\mathrm{sp}}$.
}\fi
\end{proof}

\begin{remark}[Annuli in $\Gamma$]
\label{Remark : segments inside Gamma}
Let $\Gamma$ be a locally finite graph containing $\Gamma_S$.
Assume that $]y,z[\subset\Gamma$ is the skeleton of a virtual 
open annulus in $X$, such that no bifurcation point of $\Gamma$ lie in 
$]y,z[$. Condition 
\eqref{eq : condition on pt to ensure graph in Gamma} 
means, in this case, that the radii that are spectral non solvable 
at $x$ are all $\log$-affine over an open interval in $]y,z[$ containing $x$, 
while radii that are solvable or over-solvable at $x$, 
are allowed to have a break at $x$.
Indeed
\eqref{eq : condition on pt to ensure graph in Gamma} implies 
convexity of $H_{S,j}(-,\Fs)$ along $]y,z[$, and it is known that 
$H_{S,j}(-,\Fs)$ is concave if $j$ is spectral non solvable at $x$ 
\cite[Thm. 4.7, i)]{NP-I}.
\end{remark}

\begin{corollary}[Annuli]\label{cor : Condition on an annulus on graphs}
Let $X$ be an open annulus with empty triangulation. 
Let $I$ be the skeleton of the annulus. 
Let $i\leq r=\mathrm{rank}(\Fs)$. Assume that at each point 
$x$ of $I$, and for all 
$j\in\{1,\ldots,i\}$ one of the following conditions holds:
\begin{enumerate}
\item there exists an open subinterval $J\subseteq I$ containing $x$ 
such that 
the partial height $H_{S,j}(-,\Fs)$ 
is a $\log$-affine map on $J$ 
(cf. \cite[Def. 3.1.1]{Potentiel})
\item $\R_{S,j}(x,\Fs)$ is solvable or over-solvable at $x$.
\end{enumerate}
Then 
\begin{equation}
\Gamma_{S,j}(\Fs)\;=\;\Gamma_{S,j}'(\Fs)\;=\; I\;,
\quad \textrm{for all }j=1,\ldots,i\;.
\end{equation}
\end{corollary}
\begin{proof}
Apply Corollary \ref{Cor : effective computation of Gamma} 
to $\Gamma=I$.
\end{proof}

\begin{remark}
\label{Cor : Gamma' - Gamma - Gamma'}
Let $1\leq i_1\leq i_2\leq r$. Let $\Gamma$ be a locally finite graph 
containing $\Gamma_S$. Assume that 
\begin{enumerate}
\item For all germ of segment $b$ in $\Gamma-\Gamma_S$ at least 
one of the radii $\R_{S,j}(-,\Fs)$, $j\in\{1,\ldots,i_2\}$, 
has a non zero slope on $b$.
\item The conditions of Corollary 
\ref{Cor : effective computation of Gamma}, are 
fulfilled for $i=i_1$.
\end{enumerate}
Then by Prop. 
\ref{Prop. Characterization of the skeleton by the non constancy of the partial heights} we have
\begin{equation}
\Gamma_{S,i_1}'(\Fs)\;\subseteq\;
\Gamma\;\subseteq\; 
\Gamma_{S,i_2}'(\Fs)\;.
\end{equation}

If we replace condition i) by
\begin{enumerate}
\item[i$'$)] For each end point $x$ of $\Gamma$ such that 
$x\notin\Gamma_S$, 
one has $D_{S,i_2}^c(x,\Fs)=D(x)$.
\end{enumerate}
Then 
\begin{equation}
\Gamma_{S,i_1}'(\Fs)\;\subseteq\;
\Gamma\;\subseteq\; 
\Gamma_{S,i_2}(\Fs)\;.
\end{equation}
Indeed the $X-\Gamma_{S,i_2}(\Fs)$ is a disjoint union of virtual open disks. 
So if an end point $x$ of $\Gamma$ belongs to $\Gamma_{S,i_2}(\Fs)$, the 
whole segment joining $x$ to $\Gamma_S$ is also included in 
$\Gamma_{S,i_2}(\Fs)$. 
\end{remark}

\section{Explicit bounds on the size of the controlling graphs}
\label{section : explicit bound}


In this section, we use our former results on super-harmonicity to provide explicit bounds on the size of the controlling graphs~$\Gamma^\tot_{S,j}(\Fs)$. Such results hold when the curve is compact or, more generally, when the radii are $\log$-linear along the germs of segment at infinity. The bounds involve the geometry of the curve (which translates into bounds on the size of the skeleton~$\Gamma_{S}$), the rank~$r$ of the module~$\Fs$, and also the slopes of the radii of convergence at the open boundary of the curve, in the non-compact case.

The simpler bounds are obtained under super-harmonicity conditions on the partial heights $H_{S,i}(-,\Fs)$. However, using the properties of the exceptional set~$\Cs_{S,i}(\Fs)$, we are able to obtain bounds in general. 

At the end of the section, we provide specific bounds for the controlling graph of the total height $H_{S,r}(-,\Fs)$, under additional hypotheses, that are verified in geometric situations. Using results from the forecoming papers \cite{NP-V,NP-VI}, we are able to bound the size of this graph by invariants coming from the de Rham cohomology of the module.

\subsection{Marked analytic graph structures}

Let $X$ be a quasi-smooth $K$-analytic curve endowed with a pseudo-triangulation~$S$. Recall that~$\Gamma_{S}$ may be given the structure of a marked analytic subgraph of~$X$ with set of vertices~$S$ (see Definition~\ref{def:markedgraph}). By Lemma~\ref{lem:markedanalyticgraph}, it induces a structure of marked analytic graph on~$\Gamma^\tot_{S,j}(\Fs)$, for each $j\in\{1,\dotsc,r\}$. Let us introduce a refinement.

\begin{definition}
Let $j\in\{1,\dotsc,r\}$. For each edge~$e$ of~$\Gamma^\tot_{S,j}$, let~$V_{e}$ be the set of points of~$e$ where at least one of the maps $\Rc_{S,1}(-,\Fs),\dotsc,\Rc_{S,j}(-,\Fs)$ fails to be $\log$-affine. We denote by~$\Gamma^{\tot,\lin}_{S,j}(\Fs)$ the marked analytic subgraph of~$X$ with underlying graph~$\Gamma^\tot_{S,j}(\Fs)$ and set of vertices
\begin{equation}
V(\Gamma^\tot_{S,j}(\Fs)) \cup \bigcup_{e \in E(\Gamma^\tot_{S,j}(\Fs))} V_{e}.
\end{equation}

We also set $\Gamma^{\tot,\lin}_{S,0}(\Fs) := \Gamma_{S}$ as marked analytic subgraphs of~$X$.

\medbreak

We similarly define marked analytic graphs $\Gamma^{\lin}_{S,j}(\Fs) = \Gamma^{\lin}_{S}(\Rc_{S,j}(-,\Fs))$, $\Gamma^{\lin}_{S}(H_{S,j}(-,\Fs))$, etc.
\end{definition}




%




We introduce some terminology  to count the number of vertices and edges of graphs with some 
degrees taken into account. 


\begin{definition}
Let~$V$ be a finite subset of~$X$. We define the \emph{weighted cardinal} of~$V$ as 
\begin{equation}
\WC(V) := \sum_{v\in V} \deg(v) = \Card(\pi^{-1}_{\wKa}(V)).
\end{equation}
\end{definition}

\begin{notation}
Let~$\Gamma$ be a quasi-finite marked graph. We denote by~$v(\Gamma)$ the cardinal of its set of vertices. For each $n \in \NN$, we denote by $V_{n}(\Gamma)$ the set of vertices of~$\Gamma$ of arity~$n$, and by~$v_{n}(\Gamma)$ its cardinal. In particular, $v_{1}(\Gamma)$ is the number of end-points of~$\Gamma$.

If~$\Gamma$ is a subgraph of~$X$, we define similarly $\vw(\Gamma)$ and $v_{n,\rw}(\Gamma)$, for $n\in \NN$, using weighted cardinals instead of cardinals.
%
\end{notation}


\begin{definition}
Let $J$ be a segment in~$X$ on which the map $x\mapsto \deg(x)$ is constant. We define the \emph{degree} of~$J$ to be 
\begin{equation}
\deg(J) := \mathrm{Card}(\pi_{0}(\pi_{\widehat{K^\mathrm{alg}}}^{-1}(J))) = \deg(x),
\end{equation}
for any $x\in J$.
\end{definition}


\begin{notation}
Let~$\Gamma$ be a quasi-finite marked graph. We denote by~$e(\Gamma)$ the cardinal of its set of edges.

If $\Gamma$ is a subgraph of~$X$ and if the map $x\mapsto \deg(x)$ is constant on every edge, we denote by~$\ew(\Gamma)$ the sum of the degrees of its edges.
%

\end{notation}


These notions behave well with respect to extensions of scalars (see Definitions~\ref{def:liftV} and~\ref{def:liftGamma}).


\begin{lemma}\label{lem:graphK'K}
Let $L$ be a complete valued field extension of~$K$. For each finite subset~$V$ of~$X$, we have
\begin{equation}
\WC(V) = \WC(V_{L}).
\end{equation}
For each quasi-finite marked analytic subgraph~$\Gamma$ of~$X$ and each $n\in \N$, we have
\begin{equation}
\vw(\Gamma) = \vw(\Gamma_{L}),\ v_{n,\rw}(\Gamma) = v_{n,\rw}(\Gamma_{L}).
\end{equation}
The map $x\mapsto \deg(x)$ is constant on every edge of~$\Gamma$ if, and only if, it is constant on every edge of~$\Gamma_{L}$. In this case, we have 
\begin{equation}
\ew(\Gamma) = \ew(\Gamma_{L}).
\end{equation}
\hfill$\Box$
\end{lemma}

%

\begin{lemma}\label{lem:finitetree}
For every marked finite tree~$\Gamma$, we have
\begin{equation}
e(\Gamma) + 1 = v(\Gamma) \le 2 v_{1}(\Gamma) + v_{2}(\Gamma) -2.
\end{equation}
%
\end{lemma}
\begin{proof}
Let~$\Gamma$ be a marked finite tree. For a connected planar graph~$G$ of genus~$g$, we have $v(G) - e(G) = 1-g$. Since $\Gamma$ has genus~0, we find $e(\Gamma) + 1 = v(\Gamma)$.

Denote by~$E^\pm(\Gamma)$ the set of oriented edges of~$\Gamma$. We have
\begin{equation}
\Card (E^\pm(\Gamma)) = 2  e(\Gamma).
\end{equation}
For each vertex~$x$ of~$\Gamma$, the number of oriented edges of~$\Gamma$ with endpoint~$x$ is equal to the arity~$a_{x}$ of~$x$. It follows that 
\begin{equation}
\Card (E^\pm(\Gamma)) = \sum_{x\in V(\Gamma)} a_{x} = \sum_{n\ge 1} n \, v_{n}(\Gamma).
\end{equation}
We deduce that
\begin{align}
2 v(\Gamma) - 2 & = v_{1}(\Gamma) + 2 v_{2}(\Gamma) +  \sum_{n\ge 3} n \, v_{n}(\Gamma)\\
& \ge v_{1}(\Gamma) + 2 v_{2}(\Gamma) + 3 \sum_{n\ge 3}  v_{n}(\Gamma)\\
& \ge  - 2 v_{1}(\Gamma) - v_{2}(\Gamma) + 3 v(\Gamma).
\end{align}
The result follows.
\end{proof}

\begin{proposition}\label{prop:v1v2}
Let~$\Gamma$ be a quasi-finite marked graph. Let~$\Gamma_{0}$ be a marked subgraph of~$\Gamma$ such that every connected component of $\Gamma - \Gamma_{0}$ is isomorphic to a marked finite tree with root removed. Then, we have
\begin{equation}
v(\Gamma) \le v(\Gamma_{0}) + 2 v_{1}(\Gamma- \Gamma_{0}) + v_{2}(\Gamma- \Gamma_{0}) - \Card(\pi_{0}(\Gamma-\Gamma_{0}))
\end{equation}
and 
\begin{equation}
e(\Gamma) \le v(\Gamma_{0}) + 2 v_{1}(\Gamma- \Gamma_{0}) + v_{2}(\Gamma- \Gamma_{0}) - \Card(\pi_{0}(\Gamma-\Gamma_{0})).
\end{equation}
\end{proposition}
\begin{proof}
Let~$\Gamma'$ be a connected component of $\Gamma - \Gamma_{0}$. By Lemma~\ref{lem:finitetree}, we have
\begin{equation}
e(\Gamma') +1 = v(\Gamma')+1 \le 2 (v_{1}(\Gamma')+1) + v_{2}(\Gamma') -2, 
\end{equation}
hence 
\begin{equation}
e(\Gamma') = v(\Gamma') \le 2 v_{1}(\Gamma') + v_{2}(\Gamma') -1.
\end{equation}
The result follows by summing up over all connected components of $\Gamma - \Gamma_{0}$.
\end{proof}

\subsection{Bounds on disks and annuli}\label{sec:boundsondisks}


In this section, we investigate the graphs~$\Gamma^\tot_{S,j}(\Fs)$ on virtual open disks and open pseudo-annuli and bound their number of edges and vertices. 

We begin with disks. Recall Notation~\ref{nota:kappar} for~$\kappa_{r}$.

\begin{notation}
Set
\begin{equation}
\lambda_{r} := r \max(r-1,2) =
\begin{cases}
2\kappa_{r} = 2r & \textrm{if } r \in \{1,2\};\\
\kappa_{r} = r(r-1) & \textrm{if } r\ge 3.
\end{cases}
\end{equation}
\end{notation}


%
%
%


\begin{proposition}\label{prop:edgesdisksR1}
Let~$D$ be a virtual open disk endowed with the empty pseudo-triangulation. Let $\Fs$ be a module with connection on~$D$. Assume that the map $\Rc_{\emptyset,1}(-,\Fs)$ is $\log$-affine on the germ of segment at infinity~$b_{\infty}$ of~$D$ and let~$\sigma$ be its slope. 
Then, we have
\begin{equation}\label{eq:sigmadisk}
\frac{1}{r} \,v_{1,\rw}(\Gamma^{\lin}_{\emptyset,1}(\Fs)) + \frac{1}{\kappa_{r}} \, v_{2,\rw}(\Gamma^{\lin}_{\emptyset,1}(\Fs)) \le \deg(b_{\infty})\, \sigma
\end{equation}
and 
\begin{equation}\label{eq:sigmadisk2}
v_{\rw}(\Gamma^{\lin}_{\emptyset,1}(\Fs)) = e_{\rw}(\Gamma^{\lin}_{\emptyset,1}(\Fs)) \le \deg(b_{\infty}) \lambda_{r} \sigma.
\end{equation}
In particular, $\sigma \geq 0$ and we have $\sigma=0$ if, and only if, $\Gamma_{\emptyset,1}(\Fs) = \emptyset$. 

%
%
\end{proposition}
\begin{proof}
We may assume that $K$ is algebraically closed.
%

Let $x$ be a vertex of $\Gamma_{\emptyset,1}^\lin(\Fs)$. Let~$b_{x,\infty}$ be the germ of segment out of~$x$ pointing towards the boundary of~$D$. Denote by~$T_{x}$ the subtree of~$\Gamma_{\emptyset,1}^\lin(\Fs)$ with root~$x$, i.e. the subtree containing all vertices~$y$ such that $x$~is between~$y$ and the boundary of~$D$. Let~$v_{x,1}$ (resp. $v_{x,2}$) be the number of vertices of~$T_{x}$ that are end-points (resp. points of arity~2) in~$\Gamma_{\emptyset,1}^\lin(\Fs)$. We claim that we have
\begin{equation}
\partial_{b_{x,\infty}}\R_{\emptyset,1}(-,\Fs) \le - \frac{v_{x,1}}{r}  - \frac{v_{x,2}}{\kappa_{r}}.
\end{equation} 
Let us prove the result by induction on~$v_{x,1} + v_{x,2}$.

\noindent $\bullet$ Assume that $v_{x,1} + v_{x,2} =1$.

In this case, $x$ is an end-point of~$\Gamma_{\emptyset,1}^\lin(\Fs)$. It follows that the graph~$\Gamma_{\emptyset,1}^\lin(\Fs)$ contains exactly one germ of segment out of~$x$, necessarily~$b_{x,\infty}$. As a result, we have
\begin{equation}
\partial_{b_{x,\infty}}\R_{\emptyset,1}(-,\Fs) = dd^c \R_{\emptyset,1}(x,\Fs) \le 0,
\end{equation}
by Corollary~\ref{cor:R1sh}. By definition of~$\Gamma_{\emptyset,1}^\lin(\Fs)$, we have $\partial_{b_{x,\infty}}\R_{\emptyset,1}(-,\Fs) \ne 0$, hence 
\begin{equation}
\partial_{b_{x,\infty}}\R_{\emptyset,1}(-,\Fs) \le -\frac1r,
\end{equation}
by Proposition~\ref{prop:slopesm/i}.

\noindent $\bullet$ Assume that $v_{x,1} + v_{x,2}  > 1$ and that the result holds for smaller values.

Let $b_{1},\dotsc,b_{t}$ be the germs of segments emanating from~$x$ that belong to~$\Gamma_{\emptyset,1}^\lin(\Fs)$ and are different from~$b_{x,\infty}$. For each $j\in \{1,\dotsc,t\}$, denote by~$y_{j}$ the point of~$D$ such that $[x,y_{j}]$ is the edge of~$\Gamma_{\emptyset,1}^\lin(\Fs)$ containing~$b_{j}$. Denote by~$b_{y_{j},\infty}$ the germ of segment out of~$x$ pointing towards the boundary of~$D$. 

\noindent $\bullet\bullet$ If $t=1$, then $x$ is of arity~2 in~$\Gamma_{\emptyset,1}^\lin(\Fs)$. We have $v_{x,1}=v_{y_{1},e}$, $v_{x,2} = v_{y_{1},2}+1$ and
\begin{equation}
- \partial_{b_{1}}\R_{\emptyset,1}(-,\Fs) = \partial_{b_{y_{1},\infty}}\R_{\emptyset,1}(-,\Fs) \le - \frac{v_{y_{1},e}}{r} - \frac{v_{y_{1},2}}{\kappa_{r}} ,
\end{equation}
by induction. Moreover, since~$x$ is a vertex of~$\Gamma_{\emptyset,1}^\lin(\Fs)$, we have 
\begin{equation}
\partial_{b_{1}}\R_{\emptyset,1}(-,\Fs) + \partial_{b_{x,\infty}}\R_{\emptyset,1}(-,\Fs) \ne 0.
\end{equation}
By Corollary~\ref{cor:R1sh}, we have $dd^c \R_{\emptyset,1}(x,\Fs) \le 0$, hence, by Corollary~\ref{cor:penteminoree}, 
\begin{equation}
\partial_{b_{1}}\R_{\emptyset,1}(-,\Fs) + \partial_{b_{x,\infty}}\R_{\emptyset,1}(-,\Fs) \le - \frac1{\kappa_{r}}.
\end{equation}
The result follows.

\noindent $\bullet\bullet$ If $t>1$, then, for each $j\in \{1,\dotsc,t\}$, we have $v_{y_{j},e} < v_{x,1}$, hence 
\begin{equation}
- \partial_{b_{j}}\R_{\emptyset,1}(-,\Fs) = \partial_{b_{y_{j},\infty}}\R_{\emptyset,1}(-,\Fs) \le - \frac{v_{y_{j},e}}{r} - \frac{v_{y_{j},2}}{\kappa_{r}} ,
\end{equation}
by induction. By Corollary~\ref{cor:R1sh}, we have 
\begin{equation}
dd^c \R_{\emptyset,1}(x,\Fs) =  \partial_{b_{x,\infty}}\R_{\emptyset,1}(-,\Fs) + \sum_{j=1}^t  \partial_{b_{j}}\R_{\emptyset,1}(-,\Fs) \le 0.
\end{equation}
The result follows.

\medbreak

Let us now prove \eqref{eq:sigmadisk}. Assume that $\Gamma_{\emptyset,1}(\Fs) = \emptyset$. Then, the map $\Rc_{\emptyset,1}(-,\Fs)$ is constant on~$D$, so we have $\sigma=0$, and the result holds.

Assume that $\Gamma_{\emptyset,1}(\Fs) \ne \emptyset$. Let us start from the boundary of~$D$ and consider the first edge of $\Gamma^{\lin}_{\emptyset,1}(\Fs)$, which is half-open. Let~$x$ be its endpoint. We have $\sigma = - \partial_{b_{x,\infty}}\R_{\emptyset,1}(-,\Fs)$, hence, by the result above, 
\begin{equation}
\sigma \ge \frac{v_{x,1}}{r}  + \frac{v_{x,2}}{\kappa_{r}} =  \frac{1}{r} \,v_{1,\rw}(\Gamma^{\tot,\lin}_{\emptyset,i}(\Fs)) + \frac{1}{\kappa_{r}} \, v_{2,\rw}(\Gamma^{\tot,\lin}_{\emptyset,i}(\Fs)).
\end{equation}

Let us finally prove \eqref{eq:sigmadisk2}. We have 
\begin{equation}
 \sigma \ge \frac{1}{\lambda_{r}}\, \bigl(2 v_{1,\rw}(\Gamma^{\tot,\lin}_{\emptyset,i}(\Fs)) + v_{2,\rw}(\Gamma^{\tot,\lin}_{\emptyset,i}(\Fs))\bigr),
\end{equation}
since $\lambda_{r} = \kappa_{r} \ge 2r$ as soon as $r\ge 3$. The result now follows from Proposition~\ref{prop:v1v2}, applied with $\Gamma_{0} = \emptyset$.

\end{proof}

\begin{remark}
In the setting of Proposition~\ref{prop:edgesdisksR1}, when $\Gamma_{\emptyset,1}(\Fs) \ne \emptyset$, the same proof shows that 
\begin{equation}
v_{\rw}(\Gamma^{\lin}_{\emptyset,1}(\Fs)) = e_{\rw}(\Gamma^{\lin}_{\emptyset,1}(\Fs)) \le \deg(b_{\infty}) (\lambda_{r} \sigma -1).
\end{equation}
We will not need this more precise result.
\end{remark}

\begin{corollary}\label{cor:edgesdisks}
Let~$D$ be a virtual open disk endowed with the empty pseudo-triangulation. Let~$\Fs$ be a module with connection of rank~$r$ on~$D$. Let $i\in\{1,\dotsc,r\}$. Assume that $\Gamma^\tot_{\emptyset,i-1}(\Fs) = \emptyset$. Assume that the map $H_{\emptyset,i}(-,\Fs)$ is $\log$-affine on the germ of segment at infinity~$b_{\infty}$ of~$D$ and let~$\sigma$ be its slope. Then, we have
\begin{equation}
\frac{1}{r} \,v_{1,\rw}(\Gamma^{\tot,\lin}_{\emptyset,i}(\Fs)) + \frac{1}{\kappa_{r}} \, v_{2,\rw}(\Gamma^{\tot,\lin}_{\emptyset,i}(\Fs)) \le \deg(b_{\infty})\, \sigma
\end{equation}
and 
\begin{equation}
v_{\rw}(\Gamma^{\tot,\lin}_{\emptyset,i}(\Fs)) = e_{\rw}(\Gamma^{\tot,\lin}_{\emptyset,i}(\Fs)) \le \deg(b_{\infty}) \lambda_{r} \sigma.
\end{equation}
In particular, $\sigma \geq 0$ and we have $\sigma=0$ if, and only if, $\Gamma_{\emptyset,i}^\tot(\Fs) = \emptyset$. 
\end{corollary}
\begin{proof}
By extending the scalars, we may assume that $K$ is algebraically closed. By assumption, the radii $\Rc_{\emptyset,1}(-,\Fs),\dotsc,\Rc_{\emptyset,i-1}(-,\Fs)$ are constant on~$D$, hence we have
\begin{equation}
\Gamma^{\tot,\lin}_{\emptyset,i}(\Fs) = \Gamma^{\lin}_{\emptyset}(\Rc_{\emptyset,i}(-,\Fs)).
\end{equation}
As a result, it is enough to study $\Rc_{\emptyset,i}(-,\Fs)$. 

If $\Rc_{\emptyset,i}(-,\Fs)$ is constant on~$D$, then the result holds. Otherwise, by \cite[Proposition~6.2]{NP-I}, which generalizes easily to open disks, $\Rc_{\emptyset,i}(-,\Fs_{\vert D})$ may be identified to a first radius of convergence. The result then follows from Proposition~\ref{prop:edgesdisksR1}.


\end{proof}

\begin{corollary}\label{cor:edgesdisks2}
Let $X$ be a quasi-smooth $K$-analytic curve endowed with a pseudo-triangulation~$S$. Let~$\Fs$ be a module with connection of rank~$r$ on~$X$. Let~$x\in X$. Let $i\in\{1,\dotsc,r\}$. Let~$\Cs^i_{x}$ be a set of connected components of~$X\setminus \{x\}$ that are virtual open discs and do not meet~$\Gamma^\tot_{S,i-1}(\Fs)$. Set $D^i_{x} := \bigcup_{D\in \Cs^i_{x}} D$ and $\bar D^i_{x} := D^i_{x} \cup \{x\}$. Let
\begin{equation}
\sigma_x^d := dd^c_{\subseteq \bar D^i_{x}}H_{S,i}(x,\Fs).
\end{equation}
Then, we have 
\begin{equation}
\frac{1}{r} \,v_{1,\rw}(\Gamma^{\tot,\lin}_{\emptyset,i}(\Fs) \cap D_{x}^i) + \frac{1}{\kappa_{r}} \, v_{2,\rw}(\Gamma^{\tot,\lin}_{\emptyset,i}(\Fs) \cap D_{x}^i) \le \sigma_{x}^d
\end{equation}
and 
\begin{equation}
e_{\rw}(\Gamma^{\tot,\lin}_{\emptyset,i}(\Fs) \cap D_{x}^i) = v_{\rw}(\Gamma^{\tot,\lin}_{\emptyset,i}(\Fs) \cap D_{x}^i) \le \lambda_{r}\, \sigma_{x}^d.
\end{equation}
In particular, $\sigma_{x}^d \geq 0$ and we have $\sigma_{x}^d=0$ if, and only if, $\Gamma_{\emptyset,i}^\tot(\Fs) \cap D_{x}^i = \emptyset$. 
\qed
\end{corollary}

Let us now consider the case of annuli.

\begin{corollary}\label{cor:Hjconcave}
Let~$C$ be a pseudo-annulus endowed with the empty pseudo-triangulation. Let~$\Fs$ be a module with connection of rank~$r$ on~$C$. Let $i\in\{1,\dotsc,r\}$. Assume that $C\cap \Gamma^\tot_{\emptyset,i-1}(\Fs) \subseteq \Gamma_{C}$. 
Then the map $H_{\emptyset,i}(-,\Fs)$ is $\log$-concave on~$\Gamma_{C}$.

In particular, if $H_{\emptyset,i-1}(-,\Fs)$ is $\log$-affine on~$\Gamma_{C}$, then $\Rc_{\emptyset,i}(-,\Fs)$ is $\log$-concave on~$\Gamma_{C}$.
\end{corollary}
\begin{proof}
Let~$x\in \Gamma_{C}$. Let~$\sigma_{-}$ and~$\sigma_{+}$ be the slopes of the 
map $H_{\emptyset,i}(-,\Fs)$ at~$x$ along the two germs of segment that belong to~$\Gamma_{C}$. 
Let~$\Ds^i_{x}$ be the set of other germs of segment. Each of them 
corresponds to a virtual open disk with boundary~$x$ that does not 
meet~$\Gamma^\tot_{\emptyset,i-1}$. By Corollary~\ref{cor:edgesdisks2}, we have $dd^c_{\subseteq \bar D_{x}^i} H_{\emptyset,i}
(x,\Fs) \ge 0$. 
By Theorem~\ref{thm:curveTR}, we have $dd^c H_{\emptyset,i}
(x,\Fs) \le 0$, hence
\begin{equation}
\deg(C) \, (\sigma_{-}+\sigma_{+}) \;=\; 
dd^c H_{\emptyset,i}(x,\Fs) - dd^c_{\subseteq \bar D_{x}^i} H_{\emptyset,i}
(x,\Fs)\;\le \;0\;.
\end{equation}
The result follows.
\end{proof}

\begin{corollary}\label{cor:edgesJ}
Let~$C$ be a pseudo-annulus endowed with the empty pseudo-triangulation. Let~$\Fs$ be a module with connection of rank~$r$ on~$C$. Let $i\in\{1,\dotsc,r\}$. Assume that $C\cap \Gamma^\tot_{\emptyset,i-1}(\Fs) = \Gamma_{C}$. 
Assume that $H_{\emptyset,i}(-,\Fs)$ is $\log$-affine on the two germs of segment~$b_{-}$ and~$b_{+}$ at infinity in~$C$. Let~$\sigma_{-}$ and~$\sigma_{+}$ be the associated slopes. Then, we have
\begin{equation}\label{eq:GammaC1}
 \frac{1}{\kappa_{r}} \, v_{\rw}(\Gamma^{\tot,\lin}_{\emptyset,i}(\Fs) \cap \Gamma_{C}) =  \frac{1}{\kappa_{r}} \, \big(e_{\rw}(\Gamma^{\tot,\lin}_{\emptyset,i}(\Fs) \cap \Gamma_{C}) -1\big)  \le \deg(C)\,(\sigma_{-} + \sigma_{+}),
\end{equation}
\begin{equation}\label{eq:GammaC2}
\frac{1}{r} \,v_{1,\rw}(\Gamma^{\tot,\lin}_{\emptyset,i}(\Fs) \cap (C -\Gamma_{C}) ) + \frac{1}{\kappa_{r}} \, v_{2,\rw}(\Gamma^{\tot,\lin}_{\emptyset,i}(\Fs) \cap (C -\Gamma_{C}))\le \deg(C)\,(\sigma_{-} + \sigma_{+})
\end{equation}
and 
\begin{equation}\label{eq:GammaC3}
\begin{cases}
v_{\rw}(\Gamma^{\tot,\lin}_{\emptyset,i}(\Fs)) \le \deg(C)\, (2\lambda_{r} (\sigma_{-} + \sigma_{+})) ;\\ e_{\rw}(\Gamma^{\tot,\lin}_{\emptyset,i}(\Fs)) \le \deg(C)\, (2\lambda_{r} (\sigma_{-} + \sigma_{+})+1).
\end{cases}
\end{equation}
\end{corollary}
\begin{proof}
We may assume that~$K$ is algebraically closed. 

Consider the map~$H_{\emptyset,i}(-,\Fs)$ along~$\Gamma_{C}$ in the direction going from~$b_{-}$ to~$b_{+}$. 
Let $x_{1},\dotsc,x_{s}$ be its break-points.
For each $j \in \{1,\dotsc,s\}$, let $\tau_{j}$ be the difference between the slope going out 
and the slope going in at the point~$x_{j}$. By Corollaries~\ref{cor:Hjconcave} and~\ref{cor:penteminoree}, we have $\tau_{j} \le - \frac1{\kappa_{r}}$. 

We have $\sum_{j=1}^s \tau_{j} = -(\sigma_{-}+\sigma_{+})$, hence $\frac{s}{\kappa_{r}} \le \sigma_{-}+\sigma_{+}$. Since $s = v_{\rw}(\Gamma^{\tot,\lin}_{\emptyset,i}(\Fs) \cap \Gamma_{C}) = e_{\rw}(\Gamma^{\tot,\lin}_{\emptyset,i}(\Fs) \cap \Gamma_{C}) -1$, \eqref{eq:GammaC1} follows.


\medbreak

Let $z\in \Gamma_{C}$. 
By Theorem~\ref{thm:curveTR}, we have $dd^c H_{\emptyset,i}(z,\Fs) \le 0$. Denote by~$C_{z}$ the preimage of~$z$ by the retraction map $C \to \Gamma_{C}$.

If $z$ is not a break-point, then the sum of the two slopes in the directions 
that belong to~$\Gamma_{C}$ is zero, hence $dd^c_{\nsubseteq \Gamma_{C}} H_{\emptyset,i}(z,\Fs)\le 0$. 
By Corollary~\ref{cor:edgesdisks2}, we deduce that $\Gamma^{\tot,\lin}_{\emptyset,i}(\Fs) \cap C_{z} = \emptyset$.

If $z=x_{j}$ for some $j \in \{1,\dotsc,s\}$, we have 
$dd^c_{\nsubseteq \Gamma_{C}} H_{\emptyset,i}(x_{j},\Fs) \le -\tau_{j}$. By 
Corollary~\ref{cor:edgesdisks2}, we have 
\begin{equation}
\frac{1}{r} \,v_{1,\rw}(\Gamma^{\tot,\lin}_{\emptyset,i}(\Fs) \cap C_{z}) + \frac{1}{\kappa_{r}} \, v_{2,\rw}(\Gamma^{\tot,\lin}_{\emptyset,i}(\Fs) \cap C_{z}) \le - \tau_{j}.
\end{equation}

By summing up, \eqref{eq:GammaC2} follows.

\medbreak

The inequalities~\eqref{eq:GammaC3} follow by using Proposition~\ref{prop:v1v2} with $\Gamma_{0}=\Gamma^{\tot,\lin}_{\emptyset,i}(\Fs) \cap \Gamma_{C}$, as in the proof of Proposition~\ref{prop:edgesdisksR1}.
\end{proof}


\begin{corollary}\label{cor:annulusloglinear}
Let $C$ be an open pseudo-annulus endowed with the 
empty pseudo-triangulation. Let $\Fs$ be a module with connection of rank~$r$ over $C$. Let $i\in \{1,\dotsc,r\}$ such that, for each $j\in\{1,\dotsc,i\}$, $\Rc_{\emptyset,j}(-,\Fs_{C})$ is $\log$-affine along~$\Gamma_C$. Then, we have $\Gamma^\tot_{\emptyset,i}(\Fs_{C}) =\Gamma_{C}$.
\end{corollary}
\begin{proof}
We argue by induction, using Corollary~\ref{cor:edgesJ} at each step.
\end{proof}

\begin{remark}
From the above Lemmas it is possible to derive the following simple 
criterion. 

Assume that the following conditions are fulfilled :
\begin{enumerate}
\item $\R_{S,1}(s,\Fs)=1$ for all $s\in S$,
\item Let $C$ be an open virtual annulus such that
\begin{enumerate}
\item $C$ is a connected component of $X-S$. 
\item The topological closure $I$ of $\Gamma_C$ in $X$ is an 
open or semi-open interval (not a loop).
\end{enumerate}  
Identify $\Gamma_C$ with $]0,1[$. 
If $I=\Gamma_C$ is open, then assume that 
\begin{equation}
\lim_{x\to 0^+}\R_{S,1}(x,\Fs)\;=
\;\lim_{x\to 1^-}\R_{S,1}(x,\Fs)\;=\;1\;.
\end{equation}
If $I-\Gamma_C$ is a point, say $\{0\}$, then assume that
\begin{equation}
\lim_{x\in\Gamma_C,\;x\to 1^-}\R_{S,1}(x,\Fs)\;=\;1\;.
\end{equation}
\end{enumerate} 
Then 
\begin{equation}\label{eq : RoBBA COND}
\R_{S,i}(x,\Fs)=1,\;\textrm{ for all }i=1,\ldots,r, \textrm{ and all }
x\in X\;.
\end{equation}

In alternative, consider a triangulation $S'$, i.e. a weak 
triangulation such that each connected component of $X-S'$ is 
relatively compact in $X$. Then \eqref{eq : RoBBA COND} 
holds if $\R_{S',1}(x,\Fs)=1$ for all $x\in S'$.
This is a slight generalization of  
\cite[Theorem 0.1.8]{Balda-Inventiones}.
\end{remark}

\subsection{Bounds on curves}\label{sec:globalbounds}

In this section, we bound the number of edges and vertices of the graphs~$\Gamma^{\tot,\lin}_{S,i}(\Fs)$ on a arbitrary curve. 

Let $X$ be a quasi-smooth $K$-analytic curve endowed with a pseudo-triangulation. Let~$\Fs$ be a module with connection of rank~$r$ on~$X$. Recall that, for each $i\in \{1,\dotsc,r\}$, we introduced, in Definition~\ref{def:Ei}, a subset~$\Es_{S,i}(\Fs)$ of~$X$ outside which the map $H_{S,i}(-,\Fs)$ is super-harmonic. By Corollary~\ref{cor:R1sh}, it is contained in $\Gamma^\tot_{S,i-1}(\Fs)$ and locally finite.

Let us begin with the compact case.


\begin{proposition}\label{prop:branchesXcompact}
Assume that~$X$ is compact. Let $i\in\{1,\dotsc,r\}$. Set $\Es^0_{S,i}(\Fs) := \Es_{S,i}(\Fs) - V(\Gamma^{\tot,\lin}_{S,i-1}(\Fs))$ and, for each $z\in \Es^0_{S,i}(\Fs)$, let~$d_{z}$ be the degree of the edge of~$\Gamma^{\tot,\lin}_{S,i-1}(\Fs)$ containing~$z$. Let~$W_{i}$ be a set such that $\Es_{S,i}(\Fs)\subseteq W_{i}\subseteq V(\Gamma^{\tot,\lin}_{S,i-1}(\Fs))\cup\Es_{S,i}(\Fs)$. Then, we have
\begin{equation}
\begin{cases}
\ew(\Gamma^{\tot,\lin}_{S,i}(\Fs)) \le \ew(\Gamma^{\tot,\lin}_{S,i-1}(\Fs)) + \sum_{z\in \Es^0_{S,i}(\Fs)} d_{z} + 2\lambda_{r}\cdot\sum_{x \in W_{j}} dd^c H_{S,j}(x,\Fs);\\
\vw(\Gamma^{\tot,\lin}_{S,i}(\Fs)) \le \vw(\Gamma^{\tot,\lin}_{S,i-1}(\Fs)) + \sum_{z\in \Es^0_{S,i}(\Fs)} d_{z} + 2\lambda_{r}\cdot\sum_{x \in W_{i}} dd^c H_{S,i}(x,\Fs);\\
\WC\Bigl(V_{1}(\Gamma^{\tot,\lin}_{S,i}(\Fs)) - V_{1}(\Gamma^{\tot,\lin}_{S,i-1}(\Fs))\Bigr) \le r\cdot\sum_{x \in  W_{i}} dd^c H_{S,i}(x,\Fs).
\end{cases}
\end{equation}
\end{proposition}
\begin{proof}
We may assume that~$K$ is algebraically closed. Let us first modify~$\Gamma^{\tot,\lin}_{S,i-1}(\Fs)$ by adding a vertex at each point of~$\Es^0_{S,i}(\Fs)$. The number of vertices of the resulting graph~$\Gamma^{\tot,\lin,0}_{S,i-1}(\Fs)$ is at most $e(\Gamma^{\tot,\lin}_{S,i-1}(\Fs)) + \Card(\Es^0_{S,i}(\Fs))$.

Let $x\in V(\Gamma^{\tot,\lin,0}_{S,i-1}(\Fs))$ and set $\sigma_{x}^d := dd^c_{\nsubseteq \Gamma^{\tot,\lin,0}_{S,i-1}(\Fs)} H_{S,i}(x,\Fs)$. Denote by~$C_{x}$ the preimage of~$x$ by the retraction map $X \to \Gamma^{\tot,\lin,0}_{S,i-1}(\Fs)$. For every~$x'$ such that~$[x,x']$ 
is an edge of~$\Gamma^{\tot,\lin,0}_{S,i-1}(\Fs)$, we denote by~$\sigma_{x,x'}$ the slope of $H_{S,i}(-,\Fs)$ at~$x$ in the direction of~$x'$. 

By Corollary~\ref{cor:edgesdisks2}, the number of vertices of $\Gamma^{\tot,\lin,0}_{S,i}(\Fs) \cap C_{x}$ is at most~$\lambda_{r}\sigma_{x}^d + 1$. By Corollary~\ref{cor:edgesJ}, for every~$x'$ such that $[x,x']$ is an edge 
of~$V(\Gamma^{\tot,\lin,0}_{S,i-1}(\Fs))$, the number of vertices of 
$\Gamma^{\tot,\lin}_{S,i}(\Fs)$ inside the open pseudo-annulus whose skeleton 
is~$\of{]}{x,x'}{[}$ is at most $2\lambda_{r}(\sigma_{x,x'}+\sigma_{x',x})$. Summing up, this proves the bound for the number of vertices with 
$W_{i}=V(\Gamma^{\tot,\lin}_{S,i-1}(\Fs))\cup\Es_{S,i}(\Fs)$. 
In addition, for all $y\in V(\Gamma^{\tot,\lin}_{S,i-1}(\Fs))-\Es_{S,i}(\Fs)$, one has
$dd^c H_{S,i}(y,\Fs)\leq 0$. 
This proves the bound in the general case.

Finally, the number of end-points of~$\Gamma^{\tot,\lin}_{S,i}(\Fs)$ that 
are not end-points of~$\Gamma^{\tot,\lin}_{S,i-1}(\Fs)$ is the sum of the 
number of end-points of the trees that grow out of each open edge 
and each vertex of~$\Gamma^{\tot,\lin}_{S,i-1}(\Fs)$. The result is now 
proved as in the previous paragraph.
\end{proof}

We now extend the result to non-compact curves with a condition of log-affinity of the radii at the boundary. The key ingredient is the following easy lemma.

\begin{lemma}\label{lem:opentocompact}
Assume that~$X$ is a finite curve. Let $i\in\{1,\dotsc,r\}$ and assume that the radii $\R_{S,1}(-,\Fs),\dotsc,\R_{S,i}(-,\Fs)$ are $\log$-affine along every germ of segment at infinity. 

For each $b\in\partial^o X$, let~$C_{b}'$ be a pseudo-annulus in $X-S$ whose skeleton~$\Gamma'_{b}$ is included in~$\Gamma_{S}$ and represents~$b$ and such that the radii $\R_{S,1}(-,\Fs),\dotsc,\R_{S,i}(-,\Fs)$ are $\log$-affine on~$\Gamma'_{b}$. Let~$C_{b}$ be a strict sub-pseudo-annulus of~$C'_{b}$ whose skeleton~$\Gamma_{b}$ is included in~$\Gamma'_{b}$ and still represents~$b$. Denote by~$\gamma_{b}$ the unique point of the topological boundary of~$C_{b}$ in~$C'_{b}$. Set $Y := X \setminus \bigcup_{b\in \partial^o X} C_{b}$.

Then, the following properties hold:
\begin{enumerate}
\item for each $b\in\partial^o X$, we have $\Gamma^{\tot,\lin}_{S,i}(\Fs) \cap C_{b} = \Gamma_{b}$;
\item $Y$ is compact and $\partial Y = \partial X \sqcup \{\gamma_{b},\ b\in\partial^o X\}$;
\item $S_{Y} := S  \sqcup \{\gamma_{b},\ b\in\partial^o X\}$ is a pseudo-triangulation of~$Y$ and $\Gamma_{S_{Y}} = \Gamma_{S} \cap Y$;
\item for each $j\in \{1,\dotsc,i\}$, $\Es_{S,j}(\Fs) \subseteq \Es_{S_{Y},j}(\Fs_{\vert Y}) \subseteq \Es_{S,j}(\Fs) \sqcup \{\gamma_{b},\ b\in\partial^o X\}$;
\item for each $j\in \{1,\dotsc,i\}$, $dd^c \R_{S_{Y},j}(\gamma_{b},\Fs_{|Y}) = \partial_{b} \R_{S,j}(-,\Fs)$;
\item for each $b\in\partial^o X$ and each $j\in \{1,\dotsc,i\}$, $dd^c \R_{S_{Y},j}(\gamma_{b},\Fs_{|Y}) = \partial_{b} \R_{S,j}(-,\Fs)$;
\item for each $j\in \{1,\dotsc,i\}$, the graph $\Gamma^{\tot,\lin}_{S,j}(\Fs)$ is obtained from $\Gamma^{\tot,\lin}_{S_{Y},j}(\Fs_{|Y})$ by adding the (non-relatively-compact) edges~$\Gamma_{b}$ and then removing the (bivalent) vertices~$\gamma_{b}$.
\end{enumerate}

In particular, for each $j\in \{1,\dotsc,i\}$, the graph~$\Gamma^{\tot,\lin}_{S,j}(\Fs)$ is a quasi-finite open graph and we have
\begin{equation}
\begin{cases}
\ew(\Gamma^{\tot,\lin}_{S,j}(\Fs)) = \ew(\Gamma^{\tot,\lin}_{S_{Y},j}(\Fs_{|Y}));\\
\vw(\Gamma^{\tot,\lin}_{S,j}(\Fs)) = \vw(\Gamma^{\tot,\lin}_{S_{Y},j}(\Fs_{|Y})) - N_{\infty}(X);\\
v_{1,\rw}(\Gamma^{\tot,\lin}_{S,j}(\Fs)) = v_{1,\rw}(\Gamma^{\tot,\lin}_{S_{Y},j}(\Fs_{|Y})) - N_{\infty}(X).
\end{cases}
\end{equation}
\end{lemma}
\begin{proof}
By assumption, the radii $\R_{S,1}(-,\Fs),\dotsc,\R_{S,i}(-,\Fs)$ are log-affine on~$\Gamma'_{b}$. By Corollary~\ref{cor:annulusloglinear}, we deduce that $\Gamma^{\tot,\lin}_{S,i}(\Fs) \cap C'_{b} = \Gamma'_{b}$. Property~i) follows. The others are clear. 
\end{proof}

\begin{corollary}\label{cor:branchesrayonsaffines}
Assume that~$X$ is a finite curve. Let $i\in\{1,\dotsc,r\}$ and assume that the radii $\R_{S,1}(-,\Fs),\dotsc,\R_{S,i}(-,\Fs)$ are $\log$-affine along every germ of segment at infinity. Set $\Es^0_{S,i}(\Fs) = \Es_{S,i}(\Fs) - V(\Gamma^{\tot,\lin}_{S,i-1}(\Fs))$ and, for every $z\in \Es^0_{S,i}(\Fs)$, let~$d_{z}$ be the degree of the edge of~$\Gamma^{\tot,\lin}_{S,i-1}(\Fs)$ containing~$z$. Let $W_{i}$ be a set such that $\Es_{S,i}(\Fs)\subseteq W_{i}\subseteq V(\Gamma^{\tot,\lin}_{S,i-1}(\Fs))\cup\Es_{S,i}(\Fs)$. Then, we have
\begin{equation}
\begin{cases}
\ew(\Gamma^{\tot,\lin}_{S,i}(\Fs)) \le \ew(\Gamma^{\tot,\lin}_{S,i-1}(\Fs)) +  \sum_{z\in \Es^0_{i}} d_{z}\\ 
\phantom{\ew(\Gamma^{\tot,\lin}_{S,i}(\Fs)) \le} \ 
+ 2\lambda_{r}\cdot \Bigl(\sum_{x \in W_{i}} dd^c H_{S,i}(x,\Fs)
+\sum_{b \in \partial^o X} \partial_{b} H_{S,i}(-,\Fs)\Bigr);\\
\vw(\Gamma^{\tot,\lin}_{S,i}(\Fs)) \le \vw(\Gamma^{\tot,\lin}_{S,i-1}(\Fs)) +  \sum_{z\in \Es^0_{i}} d_{z}\\ \phantom{\vw(\Gamma^{\tot,\lin}_{S,i}(\Fs)) \le} \ + 2\lambda_{r}\cdot \Bigl(\sum_{x \in W_{i}} dd^c H_{S,i}(x,\Fs)
+ \sum_{b \in \partial^o X} \partial_{b} H_{S,i}(-,\Fs)\Bigr);\\
\WC\Bigl(V_{1}(\Gamma^{\tot,\lin}_{S,i}(\Fs))-V_{1}(\Gamma^{\tot,\lin}_{S,i-1}(\Fs))\Bigr)\\ 
\quad \le r\cdot \Bigl(\sum_{x \in W_{j}} dd^c H_{S,i}(x,\Fs) + \sum_{b \in \partial^o X} \partial_{b} H_{S,i}(-,\Fs)\Bigr).
\end{cases}
\end{equation}
\end{corollary}
\begin{proof}
We use Lemma~\ref{lem:opentocompact} in order to get back to the case of a compact curve~$Y$. In the process, the set~$W_{i}$ gets changed into $W_{i} \sqcup \{\gamma_{b},\ b\in\partial^o X\}$ (with the notations of the lemma). The result now follows from Proposition~\ref{prop:branchesXcompact}.
\end{proof}

We would now like to remove the dependency on~$\Es_{S,i}(\Fs)$. We first handle the case where super-harmonicity of the partial height $H_{S,i}(-,\Fs)$ holds everywhere outside~$S$. For instance, this is always true when~$i=1$. In this case, we can choose $W_{i} = S$. We deduce the following result.

\begin{corollary}\label{cor:W=ScupEj=1}
Assume that~$X$ is a finite curve. Let $i\in\{1,\dotsc,r\}$ and assume that the radii $\R_{S,1}(-,\Fs),\dotsc,\R_{S,i}(-,\Fs)$ are log-affine along every germ of segment at infinity. Moreover, assume that $H_{S,i}(-,\Fs)$ is super-harmonic outside~$S$, \emph{i.e.} $\Es_{S,i}(\Fs) \subseteq S$. Then, we have
\begin{equation}
\begin{cases}
\ew(\Gamma^{\tot,\lin}_{S,i}(\Fs)) \le \ew(\Gamma^{\tot,\lin}_{S,i-1}(\Fs)) \\ 
\phantom{\ew(\Gamma^{\tot,\lin}_{S,i}(\Fs)) \le} \ + 2\lambda_{r}\cdot \Bigl(\sum_{x \in S} dd^c H_{S,i}(x,\Fs) +\sum_{b \in \partial^o X} \partial_{b} H_{S,i}(-,\Fs)\Bigr);\\
\vw(\Gamma^{\tot,\lin}_{S,i}(\Fs)) \le \vw(\Gamma^{\tot,\lin}_{S,i-1}(\Fs)) \\ 
\phantom{\vw(\Gamma^{\tot,\lin}_{S,i}(\Fs)) \le} \ + 2\lambda_{r}\cdot \Bigl(\sum_{x \in S} dd^c H_{S,i}(x,\Fs) +\sum_{b \in \partial^o X} \partial_{b} H_{S,i}(-,\Fs)\Bigr);\\
\WC(V_{1}(\Gamma^{\tot,\lin}_{S,i}(\Fs))-V_{1}(\Gamma^{\tot,\lin}_{S,i-1}(\Fs)))\\ 
\quad \le r\cdot \Bigl(\sum_{x \in S} dd^c H_{S,i}(x,\Fs) +\sum_{b \in \partial^o X} \partial_{b} H_{S,i}(-,\Fs)\Bigr).
\end{cases}
\end{equation}
\qed
\end{corollary}

Let us now deal with the general case, taking the sets~$\Cs_{S,i}(\Fs)$ into account.

\begin{corollary}\label{cor:W=ScupEjge2}
Assume that~$X$ is a finite curve. Let $i\in\{2,\dotsc,r\}$ and assume that the radii $\R_{S,1}(-,\Fs),\dotsc,\R_{S,i}(-,\Fs)$ are log-affine along every germ of segment at infinity. 
Then, we have
\begin{equation}
\begin{cases}
\ew(\Gamma^{\tot,\lin}_{S,i}(\Fs)) \le \ew(\Gamma^{\tot,\lin}_{S,i-1}(\Fs)) +  \WC(\Cs_{i-1}(\Fs))+ v_{1,\rw}(\Gamma^{\tot,\lin}_{S,i-1}(\Fs))\\ 
\phantom{\ew(\Gamma^{\tot,\lin}_{S,i}(\Fs)) \le} \ + 2\lambda_{r}\cdot \Bigl((i-1)\bigl(\WC(\Cs_{i-1}(\Fs))+ v_{1,\rw}(\Gamma^{\tot,\lin}_{S,i-1}(\Fs))\bigr)\\ 
\phantom{\ew(\Gamma^{\tot,\lin}_{S,i}(\Fs)) \le \quad + 2\lambda_{r}\cdot} \  + \sum_{x \in S} dd^c H_{S,i}(x,\Fs) +\sum_{b \in \partial^o X} \partial_{b} H_{S,i}(-,\Fs)\Bigr);\\
\vw(\Gamma^{\tot,\lin}_{S,i}(\Fs)) \le \vw(\Gamma^{\tot,\lin}_{S,i-1}(\Fs)) +  \WC(\Cs_{i-1}(\Fs))+ v_{1,\rw}(\Gamma^{\tot,\lin}_{S,i-1}(\Fs))\\ 
\phantom{\vw(\Gamma^{\tot,\lin}_{S,i}(\Fs)) \le} \ + 2\lambda_{r}\cdot \Bigl((i-1)\bigl(\WC(\Cs_{i-1}(\Fs))+ v_{1,\rw}(\Gamma^{\tot,\lin}_{S,i-1}(\Fs))\bigr)\\ 
\phantom{\vw(\Gamma^{\tot,\lin}_{S,i}(\Fs)) \le \ + 2\lambda_{r}\cdot} \quad + \sum_{x \in S} dd^c H_{S,i}(x,\Fs) +\sum_{b \in \partial^o X} \partial_{b} H_{S,i}(-,\Fs)\Bigr);\\
\WC(V_{1}(\Gamma^{\tot,\lin}_{S,i}(\Fs))-V_{1}(\Gamma^{\tot,\lin}_{S,i-1}(\Fs)))\\ 
\quad \le r\cdot \Bigl((i-1)\bigl(\WC(\Cs_{i-1}(\Fs))+ v_{1,\rw}(\Gamma^{\tot,\lin}_{S,i-1}(\Fs))\bigr)\\ 
\quad \phantom{\le \lambda_{r}\cdot } \quad + \sum_{x \in S} dd^c H_{S,i}(x,\Fs) +\sum_{b \in \partial^o X} \partial_{b} H_{S,i}(-,\Fs)\Bigr).
\end{cases}
\end{equation}
\end{corollary}
\begin{proof}
We may assume that~$K$ is algebraically closed. By Remark~\ref{rem:cardinalCi}, the cardinal of~$\aleph_{S,i}(\Fs)$ 
is at most~$v_{1}(\Gamma^{\tot,\lin}_{S,i-1}(\Fs))$, hence the cardinal of~$\C_{S,i}(\Fs)$ is at most $\Card(\Cs_{S,i-1}(\Fs))+v_{1}(\Gamma^{\tot,\lin}_{S,i-1}(\Fs))$. 

Set $W_{i} = S\cup\Es_{S,i}(\Fs)$. We have $\Es_{S,i}(\Fs) \subseteq W_{i} \subseteq V(\Gamma^{\tot,\lin}_{S,i-1})(\Fs) \cup \Es_{S,i}(\Fs)$. By Corollary~\ref{cor:R1sh}, $\Es_{S,i}(\Fs)-S$ and $\Es^0_{S,i}(\Fs) := \Es_{S,i}(\Fs)-V(\Gamma^{\tot,\lin}_{S,i-1}(\Fs))$ are contained in~$\C_{S,i}(\Fs)$.

By Proposition~\ref{prop:HiCi}, for every $x\in 
\C_{S,i}(\Fs)$, we have $dd^c H_{S,i}(x,\Fs) \le i-1$, hence
\begin{align}
\sum_{x\in W_{i}} dd^c H_{S,i}(x,\Fs) &= \sum_{x\in S} dd^c H_{S,i}(x,\Fs) + \sum_{x\in \Es_{S,i}(\Fs)-S} dd^c H_{S,i}(x,\Fs)\\
&\le \sum_{x\in S} dd^c H_{S,i}(x,\Fs) + (i-1)\bigl(\Card(\Cs_{S,i-1}(\Fs))+ v_{1}(\Gamma^{\tot,\lin}_{S,i-1}(\Fs))\bigr).
\end{align}
The result now follows from Corollary~\ref{cor:branchesrayonsaffines}.
\end{proof}

We finally derive upper bounds that depend only on the topology of the curve and of the behavior of the partial heights at its boundary. We will use the invariants~$\chi(x,S)$ introduced in Notation~\ref{Def; chi(x,S)}.


\begin{remark}\label{rk:gcnegatif}
Let $x\in X$. We have 
\begin{equation}
\chi(x,S) = 2 \deg(x) - 2g(x) - N_S(x) >0
\end{equation}
if, and only if, $g(x)=0$ and 
$N_S(x) \in \{0,\deg(x)\}$\footnote{Note that~$N_{S}(x)$ is a multiple of~$\deg(x)$.}.


Assume that $g(x)=0$. If $N_S(x)=0$, then~$X$ is isomorphic to a projective curve of genus~0. If $N_S(x)=\deg(x)$ and $x\notin \partial X$, then 
the skeleton~$\Gamma_{S}$ contains a branch~$b$ that ends at the interior point~$x$ of genus~0 with $\deg(b)=\deg(x)$, hence~$S$ is not minimal, by \cite[5.2.2.3]{Duc}.
\end{remark}

In what follows, we will assume that, for each $x\in S - \partial X$, we have $\chi(x,S) \le 0$. By Remark~\ref{rk:gcnegatif}, we can always achieve this condition by considering a minimal pseudo-triangulation, provided that $X$~is not a projective curve of genus~0. On the other hand, by GAGA, differential modules on projective curves of genus~0 are algebraic, hence trivial, and it is quite harmless to exclude them from our study. 


\bigbreak

We first consider the case where super-harmonicity of the partial height $H_{S,i}(-,\Fs)$ holds everywhere outside~$S$. This is alway the case when~$i=1$. 

\begin{notation}
Let $i\in \{1,\dotsc,r\}$. For $x\in S$, set
\begin{equation}
\Delta_{S,i}(x,\Fs) = dd^c H_{S,i}(x,\Fs) + \chi(x,S)\, i.
\end{equation}

Assume that~$X$ is a finite curve and that the map $H_{S,i}(-,\Fs)$ is $\log$-affine along every germ of segment at infinity in~$X$. Set 
\begin{equation}
\Irr_{S,i}(\Fs) := \sum_{x \in \partial X} \Delta_{i}(x,\Fs) + \sum_{b \in  \partial^o X} \partial_{b} H_{S,i}(-,\Fs).
\end{equation}
\end{notation}

The reason for the notation~$\Irr_{S,i}(\Fs)$ will be clarified in our subsequent papers. There, the quantity~$\Irr_{S,r}(\Fs)$ will appear as the irregularity of~$\Fs$.


\begin{corollary}\label{cor:RS1edge}
Assume that~$X$ is a finite curve. Assume that~$S$ meets every connected component of~$X$ and that, for each $x\in S- \partial X$, we have $\chi(x,S) \le 0$. Let $i\in\{1,\dotsc,r\}$. Assume moreover that the radii $\R_{S,1}(-,\Fs),\dotsc,\R_{S,i}(-,\Fs)$ are $\log$-affine along every germ of segment at infinity in~$X$ and that $H_{S,i}(-,\Fs)$ is super-harmonic outside~$S$, \emph{i.e.} $\Es_{S,i}(\Fs) \subseteq S$. Then, we have
\begin{equation}
\begin{cases}
\ew(\Gamma^{\tot,\lin}_{S,i}(\Fs)) \;\le\; \ew(\Gamma^{\tot,\lin}_{S,i-1}(\Fs)) +2 \lambda_{r} (-\chi_{c}(X) i + \Irr_{S,i}(\Fs));\bigskip\\
\vw(\Gamma^{\tot,\lin}_{S,i}(\Fs)) \;\le\; \vw(\Gamma^{\tot,\lin}_{S,i-1}(\Fs)) +2 \lambda_{r} (-\chi_{c}(X) i + \Irr_{S,i}(\Fs));\bigskip\\
\WC\bigl(V_{1}(\Gamma^{\tot,\lin}_{S,i}(\Fs))-V_{1}(\Gamma^{\tot,\lin}_{S,i-1})\bigr) \;\le\; r (-\chi_{c}(X) i + \Irr_{S,i}(\Fs)).
\end{cases}
\end{equation}
\end{corollary}
\begin{proof}
By Theorem~\ref{thm:curveTR} and Remark~\ref{rem:noTR}, for each $x\in S - \partial X$ and each~$i\in \{1,\dotsc,r\}$, we have 
\begin{equation}\label{eq: dddddc H_SI}
dd^c H _{S,i}(x,\Fs) \;\le\; 
-\chi(x,S)\, \min(i,i^\sp_{x}) \;\le\; 
-\chi(x,S)\, i,
\end{equation} 
since $\chi(x,S) \le 0$ by assumption. We deduce that 
\begin{equation}\label{eq:ddcHj}
\renewcommand{\arraystretch}{1.3}
\begin{array}{rl}
\sum_{x \in S} dd^c H_{S,i}(x,\Fs) & = \sum_{x \in S-\partial X} dd^c H_{S,i}(x,\Fs) + \sum_{x \in \partial X} dd^c H_{S,i}(x,\Fs) \medskip\\
& \le \sum_{x \in S-\partial X} -\chi(x,S)i  + \sum_{x \in \partial X} dd^c H_{S,i}(x,\Fs)\medskip\\
& \le \sum_{x \in S} -\chi(x,S)i  +  \sum_{x \in \partial X} \Delta_{S,i}(x,\Fs)\medskip\\
& \le -\chi_{c}(X)i + \sum_{x \in \partial X} \Delta_{S,i}(x,\Fs)
\end{array}
\end{equation}
by Lemma~\ref{Lemma : chi(X)=sum chi(x,S)}. The result now follows from Corollary~\ref{cor:W=ScupEj=1}.
\end{proof}

\begin{corollary}\label{cor:HSjedgetotal}
Assume that~$X$ is a finite curve. Assume that~$S$ meets every connected component of~$X$ and that, for each $x\in S - \partial X$, we have $\chi(x,S) \le 0$. Let $i\in\{1,\dotsc,r\}$. Assume moreover that the radii $\R_{S,1}(-,\Fs),\dotsc,\R_{S,i}(-,\Fs)$ are $\log$-affine along every germ of segment at infinity in~$X$ and that $H_{S,2}(-,\Fs),\dotsc,H_{S,i}(-,\Fs)$ are super-harmonic outside~$S$, \emph{i.e.} $\bigcup_{2\le j\le i}\Es_{S,j}(\Fs) \subseteq S$. Then, we have
\begin{equation}
\begin{cases}
\ew(\Gamma^{\tot,\lin}_{S,i}(\Fs)) \le \ew(\Gamma_{S}) + \lambda_{r}\bigl( - \chi_{c}(X)i(i+1) + 2 \sum_{j=1}^i \Irr_{S,j}(\Fs)\bigr);\\
\vw(\Gamma^{\tot,\lin}_{S,i}(\Fs)) \le  \vw(\Gamma_{S}) + \lambda_{r}\bigl( - \chi_{c}(X)i(i+1) + 2 \sum_{j=1}^i \Irr_{S,j}(\Fs)\bigr);\\
v_{1,\rw}(\Gamma^{\tot,\lin}_{S,i}(\Fs)) \le v_{1,\rw}(\Gamma_{S}) +r \bigl(-\chi_{c}(X) \frac{i(i+1)}{2} 
+ \sum_{j=1}^i \Irr_{S,j}(\Fs)\bigr).
\end{cases}
\end{equation}
\qed
\end{corollary}

\begin{remark}
We will prove in the next sections that the condition 
$\Es_{S,i}\subseteq S$ is related to the presence of Liouville numbers 
in the equation at the points of $\C_{S,i}$. 
Moreover we will prove the super-harmonicity under the 
further assumption that the radii of $\Fs$ equals those if its dual.
\end{remark}


Without the super-harmonicity assumption, we still find bounds that depend only on~$r$, $\Gamma_{S}$ and the slopes of the partial heights at infinity, although they are more complicated.

\begin{corollary}\label{cor:HSjedgetotalCS}
Assume that~$X$ is a finite curve. Assume that~$S$ meets every connected component of~$X$ and that, for each $x\in S - \partial X$, we have $\chi(x,S) \le 0$. Let $i\in\{1,\dotsc,r\}$. Assume moreover that the radii $\R_{S,1}(-,\Fs),\dotsc,\R_{S,r}(-,\Fs)$ are $\log$-affine along every germ of segment at infinity in~$X$. 

Define sequences $(f_{n})_{0\le n\le r}$, $(v_{n})_{0\le n\le r}$ and $(e_{n})_{0\le n\le r}$ by $f_{0}=v_{1,\rw}(\Gamma_{S})$, $v_{0}=\vw(\Gamma_{S})$, $e_{0}=\ew(\Gamma_{S})$ and, for each $n\in\{0,\dotsc,r-1\}$,
\begin{equation}
\begin{cases}
f_{n+1} = f_{n} -r\chi_{c}(X)(n+1) + r n \sum_{m=1}^n f_{m} +r \Irr_{S,n+1};\\
v_{n+1} = v_{n} -2\lambda_{r}\chi_{c}(X)(n+1) + (2\lambda_{r}n+1) \sum_{m=1}^n f_{m} + 2\lambda_{r} \Irr_{S,n+1};\\
e_{n+1} = e_{n} -2\lambda_{r}\chi_{c}(X)(n+1) + (2\lambda_{r}n+1) \sum_{m=1}^n f_{m} + 2\lambda_{r} \Irr_{S,n+1}.
\end{cases}
\end{equation}
Then, for each $i\in\{1,\dotsc,r\}$, we have
\begin{equation}
\begin{cases}
v_{1,\rw}(\Gamma^{\tot,\lin}_{S,i}(\Fs)) \le v_{1,\rw}(\Gamma_{S}) + f_{i};\\
\vw(\Gamma^{\tot,\lin}_{S,i}(\Fs)) \le v_{i};\\
\ew(\Gamma^{\tot,\lin}_{S,i}(\Fs)) \le e_{i}.
\end{cases}
\end{equation}
\end{corollary}
\begin{proof}
We may assume that~$K$ is algebraically closed. Let~$n\in\{1,\dotsc,r\}$. By Definition~\ref{Def. C_i}, the cardinal of~$\C_{S,n}(\Fs)$ is at most the sum of the cardinals of the $\aleph_{S,m}(\Fs)$'s, for $m\in\{2,\dotsc,n\}$. By Remark~\ref{rem:cardinalCi}, the cardinal of such an~$\aleph_{S,m}(\Fs)$ is itself at most the number of end-points of~$\Gamma^{\tot}_{S,m-1}(\Fs)$ that do not belong to~$\Gamma_{S}$. We deduce that
\begin{equation}
\Card(\C_{S,n}(\Fs)) \le \sum_{m=1}^{n-1} v_{1,\rw}(\Gamma^{\tot,\lin}_{S,m}(\Fs)).
\end{equation}

The result now follows from Corollary~\ref{cor:RS1edge} (with $i=1$) and Corollary~\ref{cor:W=ScupEjge2}, using \eqref{eq:ddcHj} again.
\end{proof}

To finish this section, let us consider the specific case of elliptic curves. 

\begin{corollary}\label{Cor Ell curves}
Assume that~$X$ is an elliptic curve and that the pseudo-triangulation~$S$ is minimal. Then, for each $i\in \{1,\dotsc,r\}$, the map $\Rc_{S,i}(-,\Fs)$ is constant on~$X$ and we have $\Gamma_{S,i}^{\tot,\lin}(\Fs) = \Gamma_{S}$. 

Moreover, the module with connection~$(\Fs,\nabla)$ admits a direct sum decomposition
\begin{equation}
\Fs = \bigoplus_{\rho \in ]0,1]} \Fs^\rho
\end{equation}
with the property that, for each $\rho \in ]0,1]$ such that $\Fs^\rho \ne 0$ and each $j\in \{1,\dotsc,\mathrm{rank}(\Fs^\rho)\}$, we have $\Rc_{S,j}(-,\Fs^\rho) = \rho$.
\end{corollary}
\begin{proof}
Let us first prove that, for each $i\in \{1,\dotsc,r\}$, we have $\Gamma_{S,i}^{\tot,\lin}(\Fs) = \Gamma_{S}$. We will proceed by induction on~$i$.

Let us begin with $i=1$. By Remark~\ref{rk:gcnegatif}, for each~$x\in S$, we then have $\chi(x,S)\le 0$. We have $\chi_{c}(X) =0$ and $\partial X = \partial^o X = \emptyset$, hence, by Corollary~\ref{cor:RS1edge} (whose assumptions are satisfied by Corollary~\ref{cor:R1sh}), we have $\Gamma^{\tot,\lin}_{S,1}(\Fs)=\Gamma_{S}$. 

Let $i\in \{1,\dotsc,r-1\}$ and assume that we have $\Gamma_{S,i}^{\tot,\lin}(\Fs) = \Gamma_{S}$. By Corollary~\ref{cor:R1sh} and Remark~\ref{Rk : C_Sis position}, 
we have $\Es_{S,i+1}(\Fs) \subseteq S$. It then follows from Corollary~\ref{cor:RS1edge} again that $\Gamma^{\tot,\lin}_{S,i+1}(\Fs)=\Gamma_{S}$. This proves the result.

\medbreak

Let us now prove that, for each $i\in \{1,\dotsc,r\}$, the map $\Rc_{S,i}(-,\Fs)$ is constant on~$X$. 

Let $i\in \{1,\dotsc,r\}$. It is enough to prove that the map $\Rc_{S_{\wKa},i}(-,\Fs_{\wKa})$ is constant on~$X_{\wKa}$. Note that the pseudo-triangulation~$S_{\wKa}$ may no longer be minimal. However, its associated skeleton $\Gamma_{S_{\wKa}} = \pi^{-1}_{\wKa}(\Gamma_{X})$ is still the analytic skeleton of~$X_{\wKa}$ (which is the smallest skeleton). Since the radii only depend on the skeleton and not on the triangulation itself, we may replace~$S_{\wKa}$ by a minimal triangulation~$S'$ of~$X_{\wKa}$. As a result, up to replacing~$X$ by~$X_{\wKa}$ and~$S$ by~$S'$, we may assume that~$K$ is algebraically closed.

By Remark~\ref{rem:genus}, we have 
\begin{equation}
\sum_{s\in S} g(s) + 1 - \chi_{\topo}(\Gamma_S) = g(X) = 1.
\end{equation}

We distinguish two cases.

$\bullet$ Assume that $\chi_{\topo}(\Gamma_S)=1$. (This is the good reduction case.)

Then there exists a point $x_0 \in X$ such that $g(x_0)=1$ and every other point has genus~0. We have $\Gamma_{S} = S = \{x_{0}\}$, so $\Rc_{S,i}(-,\Fs)$ is constant on~$\Gamma_{S}$, hence on~$X$, since $\Gamma_{S,i}^{\tot,\lin}(\Fs) = \Gamma_{S}$, by the first part of the proof.

$\bullet$ Assume that $\chi_{\topo}(\Gamma_S)=0$. (This is the bad reduction case, also known as \emph{Tate curve} case.)

Then $\Gamma_S$ is a circle and $X$ contains no points with positive genus. By definition of~$\Gamma^{\tot,\lin}_{S,i}(\Fs)$, the map $\Rc_{S,i}(-,\Fs)$ is $\log$-affine and continuous on~$\Gamma_{S}$. It follows that it is constant on~$\Gamma_{S}$, hence on~$X$.

\medbreak

The last part of the theorem follows from \cite{NP-III}.

\end{proof}

\begin{remark}
The same proof applied to $\mathbb{P}^{1,\mathrm{an}}_K$ endowed with 
a minimal triangulation~$S$ (consisting of a point) shows 
that, in this case, the map $\R_{S,1}(-,\Fs)$ is constant, equal to~1.
This confirms the fact that every differential equation over 
$\mathbb{P}^{1,\mathrm{an}}_K$ is trivial (because it is algebraic by GAGA).
\end{remark}

\subsection{Controlling graphs of partial heights}

In this section, we prove bounds for the controlling graphs $\Gamma^{\lin}_{S}(H_{S,j}(-,\Fs))$ of the partial heights. The proofs are very similar to that of Sections~\ref{sec:boundsondisks} and~\ref{sec:globalbounds} and we only explain which modifications are required. We will need an additional assumption of super-harmonicity.


Let us start with Corollary~\ref{cor:edgesdisks}. In this result, we show that $H_{\emptyset,j}(-,\Fs)$ behaves as a first radius of convergence on a virtual open disk~$D$, so as to be able to use Proposition~\ref{prop:edgesdisksR1}. Going back to the proof of the latter, we observe that we actually use a single property of~$H_{\emptyset,j}(-,\Fs)$, namely its super-harmonicity on~$D$. Assuming that it holds, by copying the proof of Proposition~\ref{prop:edgesdisksR1}, we directly get the following result.



\begin{proposition}\label{prop:Hjdisk}
Let~$D$ be a virtual open disk endowed with the empty pseudo-triangulation. Let $\Fs$ be a module with connection of rank~$r$ on~$D$. Let $i\in\{1,\dotsc,r\}$. Assume that the map $H_{\emptyset,i}(-,\Fs)$ is $\log$-affine on the germ of segment at infinity~$b_{\infty}$ of~$D$ and let~$\sigma$ be its slope. Assume, moreover, that $H_{\emptyset,i}(-,\Fs)$ is super-harmonic on~$D$. Then, we have
\begin{equation}
\frac{1}{r} \,v_{1,\rw}(\Gamma^{\lin}_{\emptyset}(H_{\emptyset,i}(-,\Fs))) + \frac{1}{\kappa_{r}} \, v_{2,\rw}(\Gamma^{\lin}_{\emptyset}(H_{\emptyset,i}(-,\Fs))) \le \deg(b_{\infty})\, \sigma
\end{equation}
and 
\begin{equation}
v_{\rw}(\Gamma^{\lin}_{\emptyset}(H_{\emptyset,i}(-,\Fs))) = e_{\rw}(\Gamma^{\lin}_{\emptyset}(H_{\emptyset,j}(-,\Fs))) \le \deg(b_{\infty}) \lambda_{r} \sigma.
\end{equation}
In particular, $\sigma \geq 0$ and we have $\sigma=0$ if, and only if, $\Gamma^{\lin}_{\emptyset}(H_{\emptyset,i}(-,\Fs)) = \emptyset$. 
%
%
\qed
\end{proposition}

Corollary~\ref{cor:edgesdisks2} may now be easily transfered to this setting. We skip it and consider directly the analogue of Corollary~\ref{cor:edgesJ}.

\begin{corollary}\label{cor:Hjannulus}
Let~$C$ be an open pseudo-annulus endowed with the empty pseudo-triangulation. Let $\Fs$ be a module with connection of rank~$r$ on~$C$. Let $i\in\{1,\dotsc,r\}$. Assume that $H_{\emptyset,i}(-,\Fs)$ is $\log$-affine on the two germs of segment~$b_{-}$ and~$b_{+}$ at infinity in~$C$. Let~$\sigma_{-}$ and~$\sigma_{+}$ be the associated slopes. Assume, moreover, that $H_{\emptyset,i}(-,\Fs)$ is super-harmonic on~$C$. Then, we have
\begin{equation}\label{eq:GammaC1}
 \frac{1}{\kappa_{r}} \, v_{\rw}(\Gamma^{\lin}_{S}(H_{\emptyset,i}(-,\Fs)) \cap \Gamma_{C}) =  \frac{1}{\kappa_{r}} \, \big(e_{\rw}(\Gamma^{\lin}_{S}(H_{\emptyset,i}(-,\Fs)) \cap \Gamma_{C}) -1\big)  \le \deg(C)\,(\sigma_{-} + \sigma_{+}),
\end{equation}
\begin{equation}\label{eq:GammaC2}
\frac{1}{r} \,v_{1,\rw}(\Gamma^{\lin}_{S}(H_{\emptyset,i}(-,\Fs)) \cap (C -\Gamma_{C}) ) + \frac{1}{\kappa_{r}} \, v_{2,\rw}(\Gamma^{\lin}_{S}(H_{\emptyset,i}(-,\Fs)) \cap (C -\Gamma_{C}))\le \deg(C)\,(\sigma_{-} + \sigma_{+})
\end{equation}
and 
\begin{equation}\label{eq:GammaC3}
\begin{cases}
v_{\rw}(\Gamma^{\lin}_{S}(H_{\emptyset,i}(-,\Fs))) \le \deg(C)\, (2\lambda_{r} (\sigma_{-} + \sigma_{+})) ;\\ e_{\rw}(\Gamma^{\lin}_{S}(H_{\emptyset,i}(-,\Fs))) \le \deg(C)\, (2\lambda_{r} (\sigma_{-} + \sigma_{+})+1).
\end{cases}
\end{equation}
\qed
\end{corollary}

We now pass to more general curves and adapt Corollary~\ref{cor:W=ScupEj=1}.

\begin{corollary}\label{cor:Hjcurve}
Let~$X$ is a finite curve endowed with a pseudo-triangulation~$S$. Let~$\Fs$ be a module with connection of rank~$r$ on~$\Fs$. Let $i\in\{1,\dotsc,r\}$. Assume that $H_{S,i}(-,\Fs)$ is $\log$-affine along every germ of segment at infinity in~$X$. Assume, moreover, that $H_{S,i}(-,\Fs)$ is super-harmonic on~$X-S$. Then, we have
\begin{equation}
\begin{cases}
\ew(\Gamma^{\lin}_{S}(H_{S,i}(-,\Fs))) \le \ew(\Gamma_{S}) + 2\lambda_{r}\cdot \Bigl(\sum_{x \in S} dd^c H_{S,i}(x,\Fs) +\sum_{b \in \partial^o X} \partial_{b} H_{S,i}(-,\Fs)\Bigr);\\
\vw(\Gamma^{\lin}_{S}(H_{S,i}(-,\Fs)) ) \le \vw(\Gamma_{S}) + 2\lambda_{r}\cdot \Bigl(\sum_{x \in S} dd^c H_{S,i}(x,\Fs) +\sum_{b \in \partial^o X} \partial_{b} H_{S,i}(-,\Fs)\Bigr);\\
\WC(V_{1}(\Gamma^{\lin}_{S}(H_{S,i}(-,\Fs))) -V_{1}(\Gamma_{S})) \le r\cdot \Bigl(\sum_{x \in S} dd^c H_{S,i}(x,\Fs) +\sum_{b \in \partial^o X} \partial_{b} H_{S,i}(-,\Fs)\Bigr).
\end{cases}
\end{equation}
\qed
\end{corollary}

%

Finally, we state the analogue of Corollary~\ref{cor:RS1edge}.

\begin{corollary}\label{cor:RS1edgeH}
Let~$X$ be a finite curve endowed with a pseudo-triangulation~$S$. Assume that~$S$ meets every connected component of~$X$ and that, for each $x\in S- \partial X$, we have $\chi(x,S) \le 0$. Let~$\Fs$ be a module with connection of rank~$r$ on~$\Fs$. Let $i\in\{1,\dotsc,r\}$. Assume that $H_{S,i}(-,\Fs)$ is $\log$-affine along every germ of segment at infinity in~$X$. Assume, moreover, that $H_{S,i}(-,\Fs)$ is super-harmonic on~$X-S$. Then, we have
\begin{equation}
\begin{cases}
\ew(\Gamma^{\lin}_{S}(H_{S,i}(-,\Fs))) \le \ew(\Gamma_{S}) +2 \lambda_{r} (-\chi_{c}(X) i + \Irr_{S,i}(\Fs));\\
\vw(\Gamma^{\lin}_{S}(H_{S,i}(-,\Fs))) \le \vw(\Gamma_{S}) +2 \lambda_{r} (-\chi_{c}(X) i + \Irr_{S,i}(\Fs));\\
\WC(V_{1}(\Gamma^{\lin}_{S}(H_{S,i}(-,\Fs)))-V_{1}(\Gamma_{S})) \le r (-\chi_{c}(X) i + \Irr_{S,i}(\Fs)).
\end{cases}
\end{equation}
\qed
\end{corollary}

In a forthcoming paper, we will show that the irregularity~$\Irr_{S,r}(\Fs)$ defined above has a cohomological interpretation and appears in a Grothendieck-Ogg-Shafarevich-type formula. 


\begin{theorem}\label{thm:GOS}
Let~$X$ be a finite curve endowed with a pseudo-triangulation~$S$. Let~$\Fs$ be a module with connection of rank~$r$ on~$X$. Assume that the radii of~$\Fs$ are spectral non-solvable at the points of~$\partial X$. Assume that~$\Fs$ is strongly free of Liouville numbers along the singular germs of segment at the points of~$\partial X$ and the germs of segment at infinity in~$X$. Assume that $H_{S,r}(-,\Fs)$ is $\log$-affine along every germ of segment at infinity in~$X$. Then, the de Rham cohomology spaces of~$X$ with coefficients in~$\Fs$ are finite-dimensional and we have
\begin{equation}
\chidr(X,\Fs) = r \cdot \chi_{c}(X) - \Irr_{S,r}(\Fs),
\end{equation}
where $\chidr(X,\Fs)$ denotes the Euler characteristic of~$X$ with respect to the de Rham cohomology with coefficients in~$\Fs$.
\qed
\end{theorem}


\begin{corollary}\label{cor:HrIndex}
Retain the hypotheses of Theorem~\ref{thm:GOS}. Assume, moreover, that $S$ meets every connected component of~$X$ and that, for each $x\in S- \partial X$, we have $\chi(x,S) \le 0$. Assume also $H_{S,r}(-,\Fs)$ is super-harmonic on~$X-S$. 

Then, we have
\begin{equation}
\begin{cases}
\ew(\Gamma^{\lin}_{S}(H_{S,r}(-,\Fs))) \le \ew(\Gamma_{S}) - 2\lambda_{r}\chidr(\Fs);\\ 
\vw(\Gamma^{\lin}_{S}(H_{S,r}(-,\Fs))) \le \vw(\Gamma_{S}) - 2\lambda_{r}\chidr(\Fs);\\ 
v_{1,\rw}(\Gamma^{\lin}_{S}(H_{S,r}(-,\Fs)) \le v_{1,\rw}(\Gamma_{S}) - r\chidr(\Fs).
\end{cases}
\end{equation}
\qed
\end{corollary}

\section{Some counterexamples.} 
\label{An explicit counterexample.}


In this section we provide the following counterexamples:
\begin{enumerate}
\item Non compatibility of solvable and over-solvable 
radii with duals; 
\item Uncontrolled behavior of solvable and over-solvable 
radii by exact sequences; 
\item An explicit example of differential module over a 
disk for which $\Fs_{\geq i}$ is not a direct summand 
\item Some basic relations between the 
Grothendieck-Ogg-Shafarevich formula and 
super-harmonicity of partial heights 
(cf. Remark \ref{rk : super-harmonicity}). 
\end{enumerate}
All the examples involve a differential module over an open disk with 
empty triangulation. Indeed 
any possible counterexample to the above situations is reduced 
to this case by localizing to a maximal disk $D(x,S)$. 
This is because on $\Gamma_S$ all the radii are spectral and are 
compatible with duality and spectral sequences. 
Concerning iv) the potential failure of super-harmonicity at some 
points was one of the crucial difficulties of \cite{NP-I}, we here relate 
this to the presence of Liouville numbers.

\subsection{Setting.}\label{section : setting examples ..htfhg}
Let $D:=D^-(0,1)$ be the unit open disk, and let $T$ be its coordinate.
In this section all differential module will be defined over the ring 
$\O^\dag(D):=\cup_{\varepsilon>0}\O(D_\varepsilon)$, where 
$D_\varepsilon:=D^-(0,1+\varepsilon)$. Namely $\O^\dag(D)$ 
is formed by power series $f(T)=\sum_{i\geq 0}a_iT^i$ satisfying 
$\lim_i|a_i|\rho^i=0$, for some unspecified $\rho>1$.
The data of a differential module $\M$ 
over $\O^\dag(D)$ is equivalent to 
that of a differential module $\M_\varepsilon$ over 
$\O(D_\varepsilon)$ for some 
unspecified $\varepsilon>0$. The triangulation on 
$D_\varepsilon$ will always be the empty one. Moreover the 
radii are assumed to be \emph{all solvable or over-solvable} at the 
boundary $x_{0,1}$ of $D$. By Lemma 
\ref{Lemma : solvable radii at boundary}, this implies $\M$ is trivial on  
all sub-disks of $D_\varepsilon$ with boundary $x_{0,1}$. 

All differential modules will have the property $\bs{\mathrm{NL}}$ 
of non Liouville exponents. 
By definition this means that the exponents, and their 
differences, are not Liouville numbers (cf. \cite{Ch-Me-I}, 
\cite{Ch-Me-II}, \cite{Ch-Me-III}, \cite{Ch-Me-IV}).
Under this condition 
$H^1(\M,\O^\dag(D)):=\textrm{Coker}(\nabla:\M\to\M)$ 
has finite dimension.


The $\bs{\mathrm{NL}}$ condition is quite implicit. 
The effective way to ensure it is to assume either that the restriction of 
$\M$ to the annulus $C_\varepsilon=\{|T(x)|\in]1,1+\varepsilon[\}$ 
has a Frobenius structure, or that 
none of the radii $\R_{\emptyset,i}(-,\M)$ of 
$\M$ verifies the Robba condition along 
$b=[x_{0,1},x_{0,1+\varepsilon}[$:
\begin{equation}\label{eq : negative break of R_i NL}
\R_{\emptyset,i}(-,\M)\textrm{ is solvable at $x_{0,1}$, 
and $\partial_b\R_{\emptyset,i}(x_{0,1},\M)=1$}\;.
\end{equation} 
\begin{remark}\label{Rk : NL condition restr}
Condition \eqref{eq : negative break of R_i NL} differs from 
the same condition for the radii of $\M_{|C_\varepsilon}$. 
Indeed over-solvable radii over $D_\varepsilon$ are truncated by 
localization to $C_\varepsilon$ and become solvable. 
This corresponds to the existence of trivial 
submodules of $\M_\varepsilon$ that of course satisfy the condition 
$\bs{\mathrm{NL}}$. So the $\bs{\mathrm{NL}}$ condition really 
arises from the presence of a break of some 
$\R_{\emptyset,i}(-,\M)$ at $x_{0,1}$, as in 
\eqref{eq : negative break of R_i NL}, before localization.
\end{remark}
\if{
will have an \emph{unspecified} Frobenius 
structure, i.e. an isomorphism of differential modules 
$(\varphi^n)^*\M\simto\M$, for an unspecified $n\geq 1$, 
where the first module is the pull-back by the Frobenius map 
$\varphi^n:\O^\dag(D)\to\O^\dag(D)$ associating to 
$\sum_{i\geq 0}a_iT^i$ the series 
$\sum_{i\geq 0}\sigma^n(a_i)T^{ip^n}$, where 
$\sigma:K\to K$ is a lifting of the $p$-th power map 
$x\mapsto x^p$ of the residual field of $K$.
}\fi

\subsection{Grothendieck-Ogg-Shafarevich formula and 
super-harmonicity}
\label{A remark on the Grothendieck-Ogg-Shafarevich formula}
In  this section $\M$ is a differential module over $\O^\dag(D)$ of 
rank $r=\mathrm{rank}(\M)$, and $b=]x_{0,1},x_{0,1+\varepsilon}[$ 
is a germ of segment oriented as out of $x_{0,1}$.
\begin{definition}[\cite{Ch-Me-III},\cite{Ch-Me-IV}]
One defines the $p$-adic irregularity of $\M$ 
at $\infty$ as
\begin{equation}\label{Eq : def of irririr}
\mathrm{Irr}_{\infty}\M\;:=\;
-\partial_bH_{\emptyset,r}(x_{0,1},\M_{|C_\varepsilon})\;=\; 
-\sum_{i=0}^r\partial_b\R_{\emptyset,i}(x_{0,1},
\M_{|C_\varepsilon})\;.
\end{equation}
\end{definition}
Over-solvable radii of $\M$ corresponds to solutions of $\M$ on some 
$D_\varepsilon$. These radii are truncated by localization to 
$C_\varepsilon$. As a result their slope remains zero, but the 
localization to $C_{\varepsilon}$ adds $-1$ to the slope of each other 
radius. 
One finds $\partial_bH_{\emptyset,r}(x_{0,1},\M_{|C_\varepsilon})=
\partial_bH_{\emptyset,r}(x_{0,1},\M)-r+h^0(\M)$:
\begin{equation}\label{eq :IIR as intrinsic def}
\mathrm{Irr}_\infty(\M)\;=\;
\mathrm{rank}(\M)-\partial_bH_{\emptyset,r}(x_{0,1},\M_\varepsilon)- 
h^0(\M)\;.
\end{equation}
Assume now that $\M$ has the $\bs{\mathrm{NL}}$ property, and 
that $K$ is spherically complete. Then
one has the Grothendieck-Ogg-Shafarevich formula 
(often called Euler-PoincarÃ© formula):
\begin{equation}\label{eq : GOS}
h^0(\M)-h^1(\M)\;=\;
\mathrm{rank}(\M)-\mathrm{Irr}_\infty(\M)\;.
\end{equation}
By \eqref{eq :IIR as intrinsic def} this formula can be written as :
\begin{equation}\label{GOS-changed}
h^1(\M)\;=\;-\partial_bH_{\emptyset,r}(x_{0,1},\M_\varepsilon) \;.
\end{equation}

\begin{proposition}\label{Prop : super-harm everywhere}
Let $\M_\varepsilon$ be a differential module over a disk 
$D_\varepsilon$, 
such that $\R_{\emptyset,1}(-,\M_\varepsilon)$ 
is solvable or over-solvable at 
$x_{0,1}$, and such that $\M_\varepsilon$ 
has the $\bs{\mathrm{NL}}$ property.
Then for all $i=1,\ldots,r=\mathrm{rank}(\M_\varepsilon)$ 
the $i$-th partial height 
$H_{\emptyset,i}(-,\M_\varepsilon)$ is super-harmonic at $x_{0,1}$: 
$dd^cH_{\emptyset,i}(-,\M_\varepsilon)\leq0$.
\end{proposition}
\begin{proof}
Super-harmonicity is insensitive to scalar extensions of $K$. 
So we can assume that  $K$ is spherically complete.  
If $\R_{\emptyset,1}(-,\M_\varepsilon)$ 
is over-solvable at $x_{0,1}$, then $\M_\varepsilon$ 
is trivial over $D_\varepsilon$ for some $\varepsilon$, and the 
statement is trivial. Over-solvable radii do not contribute to 
super-harmonicity, 
so 
we can assume that $\M_\varepsilon$ has no 
trivial submodules over $D_\varepsilon$, i.e. 
$\R_{\emptyset,i}(x_{0,1},\M_\varepsilon)=1$ for all $i$.
For small values of $\varepsilon$, the radii are all $\log$-affine on 
$[x_{0,1},x_{0,1+\varepsilon}[$. 
So Lemma \ref{Lemma : solvable radii at boundary}, together with 
Cor. \ref{Cor : effective computation of Gamma},
and Remark \ref{Remark : segments inside Gamma}
prove that $\Gamma_{S,i}(\M_\varepsilon)$ is either equal to  
$[x_{0,1},x_{0,1+\varepsilon}[$, or empty if 
$\R_{\emptyset,i}(-,\M_\varepsilon)$ is constant.  
By \eqref{GOS-changed}, $H_{\emptyset,r}(-,\M_\varepsilon)$ is 
concave at $x_{0,1}$, and constant outside 
$[x_{0,1},x_{0,1+\varepsilon}[$, 
hence $H_{\emptyset,r}(-,\M_\varepsilon)$ is super-harmonic on $D_\varepsilon$.
The assertion for $H_{\emptyset,i}(-,\M_\varepsilon)$ 
is deduced from that of $H_{\emptyset,r}(-,\M_\varepsilon)$ by 
interpolation. Namely by convexity of the 
convergence Newton polygon one has 
$H_{\emptyset,i}(-,\M_\varepsilon)\leq 
\frac{i}{r}H_{\emptyset,r}(-,\M_\varepsilon)$, 
moreover these two functions coincide at $x_{0,1}$. This proves that 
$H_{\emptyset,i}(-,\M_\varepsilon)$ is super-harmonic. 
\end{proof}
\begin{remark}\label{rk : super-harmonicity}
In \cite[Thm.4.7]{NP-I} one proves that 
$H_{\emptyset,i}(-,\M_\varepsilon)$ 
are super-harmonic over $D_\varepsilon$ 
outside the \emph{finite} set $\mathscr{C}_i$ 
of Def. \ref{Def. C_i}. However 
$H_{\emptyset,1}(-,\M_\varepsilon)$ 
is super-harmonic on the whole disk.  
Proposition \ref{Prop : super-harm everywhere} 
shows that if $H_{\emptyset,i}(-,\M_\varepsilon)$ 
is not super-harmonic
then $\M_\varepsilon$ has some Liouville exponent, 
or possibly some non solvable radii in order to be outside 
the range of validity of G.O.S. formula. 
\if{
\framebox{
\begin{minipage}{430pt}
QUI COMPLETARE :

1) decomporre con Dwork-Robba separare il caso solubile dal caso non 
solubile. 

2) Guardare la formula di GOS in Robba IV nel caso spettrale non 
soluble. E dimostrare che la formula vale nel caso non solubile.

3) Mostrare che il caso solubile Ã¨ la formula di GOS fatta sopra 
(risciverla con un numero arbitrario di buchi surconvergenti). 

4) Dedurre che quando GOS funziona abbiamo super-harmonicitÃ .

....

....

The above example then proves then that the super-harmonicity holds 
if the Grothendieck-Ogg-Shafarevich formula holds. Unfortunately the 
assumptions on $\M$ for this formula are Non Liouville exponents ... 
completare !!!! .......

.........

........

\end{minipage}}
}\fi

\end{remark}

\begin{remark}
Since $\mathrm{Irr}_\infty(\M)$ involves spectral radii, these are 
stable by duality and one has 
$\mathrm{Irr}_\infty(\M)=\mathrm{Irr}_\infty(\M^*)$. 
Hence $h^0(\M)-h^1(\M)=h^0(\M^*)-h^1(\M^*)$ as soon as Grothendieck-Ogg-Shafarevich formula holds.
\end{remark}

\bibliographystyle{amsalpha}
\bibliography{NP}

\end{document}